\documentclass[10pt,a4paper,reqno]{article}
\usepackage[margin=2cm]{geometry}
\usepackage[T1]{fontenc}
\ifdefined\pdfgentounicode
  \input{glyphtounicode}
\fi
\usepackage{charter}
\usepackage[english]{babel} %for biblatex
\usepackage{csquotes} %for biblatex
\usepackage{filecontents} %for biblatex
\usepackage[dvipsnames]{xcolor}
\usepackage{color}

\usepackage{amsmath,amssymb,amsthm,mathtools,mathrsfs}
\usepackage[expansion=false]{microtype}

 \usepackage[
 backend=biber,
 hyperref=true,
 backref=true,
 isbn=false,
 doi=true,
 url=false,
 natbib=true,
 eprint=true,
 useprefix=true,
 maxcitenames=99,
 maxbibnames=99,  
 maxalphanames=99, 
 minalphanames=99,
 safeinputenc,
 style=alphabetic,
 citestyle=alphabetic,
 block=space,
 datamodel=ext-eprint,
 ]{biblatex}
 \usepackage[
 hypertexnames = false,
 colorlinks    = true,
 citecolor     = teal,
 linkcolor     = blue,
 urlcolor      = blue,
 linktocpage = true,
 hyperfootnotes=true,
 breaklinks
 ]{hyperref} 
 
 \renewbibmacro{in:}{} %%% for removing In : in biblatex %%%%%
 
\DeclareFieldFormat[article,unpublished,misc]{title}{\mkbibemph{#1}}
\DeclareFieldFormat[article]{journaltitle}{#1}
\DeclareFieldFormat[article]{volume}{\mkbibbold{#1}}
\AtBeginBibliography{\raggedright}

\newcommand{\ip}{\mathbin{\lrcorner}}
\newcommand{\tr}{\operatorname{tr}}

\newcommand{\Lie}{\mathcal L}
\newcommand{\ad}{\operatorname{ad}}
\newcommand{\DD}{\mathcal D}
\newcommand{\PP}{\mathcal P}
\newcommand{\QQ}{\mathcal Q}
\newcommand{\UU}{\mathcal U}
\newcommand{\BB}{\mathcal B}

\newcommand{\vol}{\operatorname{vol}}
\numberwithin{equation}{section}
\newtheorem{theorem}{Theorem}[section]
\newtheorem{proposition}[theorem]{Proposition}
\newtheorem{lemma}[theorem]{Lemma}
\newtheorem{corollary}[theorem]{Corollary}
\newtheorem{definition}[theorem]{Definition}
\theoremstyle{remark}
\newtheorem{remark}[theorem]{Remark}
\theoremstyle{plain}
\allowdisplaybreaks[2]
\AtBeginDocument{%
  \setlength{\abovedisplayskip}{16pt plus 3pt minus 2pt}
  \setlength{\belowdisplayskip}{16pt plus 3pt minus 2pt}
  \setlength{\abovedisplayshortskip}{16pt plus 3pt minus 2pt}
  \setlength{\belowdisplayshortskip}{16pt plus 3pt minus 2pt}
}
\title{Coupled flows in $G_2$ geometry: preservation, existence and stability}

\author{Shubham Dwivedi and Siva Tej Pindiprolu}
\date{\today}

\begin{document}

\maketitle

\begin{abstract}
We study coupled flows on compact manifolds with $G_2$ structures that preserve coassociative submanifolds and $G_2$ instantons as the ambient $G_2$ structure evolves.

For coassociatives, we give an explicit example on a torus showing that the Laplacian flow of closed $G_2$ structures coupled to mean curvature flow can fail to preserve coassociativity, even arbitrarily close to a flat structure. We circumvent this by adding a Lie derivative to the ambient equation, using an elliptic correction of minimum $L^2$ norm along the submanifold. We prove short-time existence and uniqueness of the  corrected coupled flow from coassociative initial data, without any rigidity assumption. We also prove global existence and smooth exponential convergence for coassociative initial pairs sufficiently close to a torsion-free
coassociative pair, with the ambient cohomology class fixed.

For instantons, we propose a similar corrected coupled system by adding a torsion term and a nonlocal elliptic correction to Yang–Mills flow along prescribed coclosed ambient $4$-forms. We prove short-time existence and uniqueness of the corrected coupled system for connections nearby an acyclic $G_2$ instanton. We obtain global existence and smooth exponential convergence for the coupled corrected connection flow and the modified coflow of coclosed $G_2$ structures with parameter zero, near a torsion-free structure with an acyclic instanton. We prove a new result on the global existence and smooth exponential convergence near a torsion-free $4$-form in a fixed cohomology class for the modified coflow with $c=0$ which is a result of independent interest.

We also give continuation criteria for both the coupled systems and  construct compact examples with nonzero correction. Finally, for every positive coflow parameter, we exhibit a flat torus with a growing linear mode which shows why $c=0$ is the meaningful one.

\end{abstract}

\tableofcontents

\section{Introduction and main results}\label{sec:intro}

A geometric condition used in defining a submanifold of a smooth manifold or a connection on a principal bundle over a smooth manifold depends on the ambient structure. Allowing that structure to evolve by a geometric flow therefore changes the equation that the submanifold or connection must satisfy.  It is natural to study coupled geometric flows, pairing the ambient evolution with a flow of the submanifold or connection and asking whether a geometrically preferred property of the submanifold (for instance minimality or being calibrated) or of the connection (for instance the connection being an instanton or Yang-Mills) is preserved along the coupled system of flows. Such a study can have consequences in a combination of manners. Depending on what properties are preserved, the coupled system can tell us about the preferred geometric flows for the ambient structures and the sub-objects. It can also tells us that what geometric property is better behaved and worth looking into in detail. This is the motivating philosophy behind the present paper and for us the setting is that of a seven dimensional smooth manifold with a $G_2$-structure where the underlying geometry is related to the nonassociative algebra of the octonions $\mathbb{O}$. A natural first choice in the coupled system is to use a familiar geometric flow for the ambient structure as well as for the special structure which will be submanifolds or instantons in our setting. However, \emph{a priori} there is no reason for the two velocities to preserve the equation relating them. We must differentiate that equation and check whether
the contributions cancel. This question is motivated by some concrete affirmative answers in literature where such a cancellation do occur: the preservation of Lagrangian submanifolds under mean curvature flow coupled to K\"ahler--Ricci flow following Smoczyk \cite{Smo, Smo-habil} and the work of Lotay and Pacini in \cite{LP} where among other things, they prove preservation of the condition of being totally real in an almost K\"ahler manifold under the coupling of the \emph{Maslow flow} and the symplectic curvature flow of Streets and Tian \cite{SCF}. Lotay and Pacini also proposed studying analogous couplings for
coassociatives and $G_2$ instantons \cite[Section~1]{LP}. We study this question by first identifying the failure of preservation for the unmodified equations and then constructing
elliptic corrections that preserve the required conditions.

\medskip

Let $M$ be a smooth seven dimensional manifold. A positive three form $\varphi$ is a differential form which is pointwise linearly equivalent to the
model in \eqref{eq:flat} on $\mathbb{R}^7$. It defines a $G_2$ structure, a metric
$g=g_\varphi$, and an orientation, with the last two depending on $\varphi$ in a nonlinear fashion. The existence of such a structure on a smooth manifold is a topological condition. A $G_2$ structure on $M$ exists if and only if the structure group of the frame bundle reduces from the group $\mathrm{GL}(7, \mathbb R)$ to the Lie group $G_2$. This condition is also equivalent to the existence of a positive $3$-form $\varphi$ as above which s equivalent to the vanishing of the first and second Stiefel-Whitney classes which in turn, are equivalent to $M$ being orientable and spinnable, respetively. We use Karigiannis' conventions \cite{Kar} and write $\psi=*\varphi$. We fix the orientation and consider positive forms with that orientation. The $G_2$-structure is closed when $d\varphi=0$, co-closed when $d\psi=0$, and torsion-free when $\nabla^{g_\varphi}\varphi=0$ where $\nabla^{g_{\varphi}}$ is the Levi-Civita connection of the metric induced by $\varphi$. Torsion-freeness of  $\varphi$ is equivalent to closedness and co-closedness together. We will keep these assumptions separate, since the two constructions use different ones. See \cite{Joyce-book, Salamon-book, Bryant-rmks, Kar} for more on manifolds with $G_2$ structures.

\medskip

Let $M$ be a manifold with a $G_2$ structure $\varphi$. Let $i:C^4\hookrightarrow M^7$ be an embedding, and let $A$ be a
connection on a principal $G$-bundle $P\to M$. The curvature $F_A$
is a two form with values in the adjoint bundle $\operatorname{ad}P$.
The defining equations for coassociativity and the $G_2$ instanton
condition are, respectively,
\begin{align}\label{eq:definingeqns}
	i^*\varphi=0,
	\qquad
	F_A\wedge\psi=0.
\end{align}
We orient a coassociative so that
$i^*\psi=\operatorname{vol}_{i^*g}$.

As is evident from their defining equations, both conditions depend on the ambient $G_2$ structure. A submanifold
can be coassociative for one structure and fail to be coassociative
for a nearby structure. Similarly, a connection which is an
instanton for one structure need not remain an instanton when the
structure changes. Thus, when we evolve the ambient structure, we
must also choose how the submanifold or connection moves if we want to preserved the special conditions in \eqref{eq:definingeqns}.

\medskip

If $d\psi=0$, then $\psi$ is a calibration form in the sense of Harvey--Lawson \cite{HL} and every coassociative
is minimal. In the first coupled system which we study in this paper, we only assume $d\varphi=0$ and not $d\psi=0$ hence the initial coassociative (or the evolving ones) need not be minimal. If a connection is an instanton and the underlying $G_2$ structure is coclosed then the connection satisfies a version of the Yang-Mills equation with an extra term consisting of the torsion of the $G_2$ structure. However, it need not satisfy the ordinary Yang--Mills equation. More about coassociatives and calibrations can be found in the seminal paper on calibrated geometry by Harvey--Lawson \cite{HL} and the lecture notes by Lotay \cite{Lot}. In both the coupled systems, our choice of assumptions of $d\varphi=0$ and $d\psi=0$ is motivated by the deformation theory of the two geometric constraints in \eqref{eq:definingeqns}. When $d\varphi = 0$, the normal linearization of the coassociative equation identifies with the exterior derivative on anti-self-dual two-forms, and the nonlinear deformation map is surjective on the space of exact three-forms. These properties underlie the unobstructed deformation theory of compact coassociatives see \cite{McL} and \cite[Theorem 5.5]{Lot}. Thus, the natural background in which to study possible evolutions of coasscitatives is when the underlying $G_2$ structure is close. Similarly, for $G_2$-instantons, the condition $d\psi=0$, together with the Bianchi identity places the linearization of the instanton equation $F_{A}\wedge\psi=0$ in an elliptic deformation complex, whose first cohomology describes infinitesimal deformations modulo gauge see for instance \cite{Wal17} or \cite[\textsection 8.4]{Lot}. Although the instanton deformation problem may have obstructions, this provides its standard elliptic framework. This is the reason for working with coclosed $G_2$ structures in the second coupled system in the paper. For the system of coupled flows in this paper, we consider Bryant's Laplacian flow of closed $G_2$ structures \cite{Bryant-rmks, BX} coupled to mean curvature flow of an initial coassociative embedding, and
Yang--Mills flow of connections coupled to a prescribed evolution of co-closed $G_2$ structures, a particular case of which is Grigorian's modified Laplacian coflow \cite{Gri}. In both the cases we first calculate the error in preservation of the coassociative or instanton condition. We then
solve a linear elliptic equation to introduce a correction which cancels that error. However, this procedure is not straight-forward and there are three separate issues: 1) the corrected system of PDEs
must be defined on a neighborhood of the initial data, 2) it must have
solutions for short-time, 3) and those solutions must preserve the required condition.
Long-time existence and convergence of solutions is an entirely different issue and require further analytic estimates.

\medskip

We write $\delta_g=d_g^*$ and use the Hodge Laplacian
$\Delta_g=d\delta_g+\delta_gd$. Subscripts on operators indicate the
current metric. For a closed $G_2$ structure, if $\tau$ denotes the only surviving \emph{intrinsic torsion form} of $\varphi$ then $\tau=\delta_\varphi\varphi$ \cite{LW}
and Laplacian flow is the evolution equation $\partial_t\varphi=\Delta_\varphi\varphi=d\tau$.
This flow was introduced by Bryant \cite{Bryant-rmks} as a way to obtain holonomy $G_2$ manifolds using parabolic methods. The short-time existence and uniqueness of solutions of this flow on closed manifolds was proved by Bryant and Xu \cite{BX}. Further analytic and geometric properties of the Laplacian flow were proved by Lotay and Wei who
established derivative estimates, uniqueness and compactness \cite{LW},
real analyticity \cite{LW19a}, and dynamical stability of torsion-free $G_2$ structures \cite{LW19b}. There is an analogous proposal for flows of coclosed $G_2$-structures called the Laplacian coflow of coclosed $G_2$ structures and that flow and its modification are studied in\cite{KMT,Gri}. For connections on principal bundles, Yang--Mills flow on torsion-free special holonomy manifolds is studied by Oliveira and Waldron \cite{OW}.

Motivated by the results on the preservation of the Lagrangian condition in an ambient K\"ahler-Einstein manifold or when the mean curvature flow of an initial Lagrangian in a K\"ahler manifold is coupled with the K\"ahler-Ricci flow of the ambient K\"ahler structures in \cite{Smo, Smo-habil} as well as the results on the Maslow flow of totally real submanifolds in an ambient symplectic curvature flow background in \cite{LP, LP-persist}, one might naively want to couple the mean curvature flow of an initial coassociative with the Laplacian flow of the ambient closed $G_2$-structures and hope to preserve the coassociative condition. Trying to attempt the same, for a coassociative moving by mean curvature in an ordinary
Laplacian-flow background, our first calculation in Proposition~\ref{prop:unmodified-pullback} gives
\begin{align*}
 \frac{\partial}{\partial t}(i^*\varphi)=d_C(*_Ci^*\tau).
\end{align*}
\emph{A priori}, there is no reason for this expression to vanish. In fact, in Proposition~\ref{prop:torus-example} we give an explicit example on the compact torus where we show that this term is non-zero even for arbitrarily small exact perturbations of a flat structure and in
the same cohomology class. Thus, the naive expectation of coupling the mean curvature flow of coassociatives with the Laplacian flow of closed $G_2$ structures cannot work. Nevertheless, we circumvent this problem as we now describe. We cancel the error introduced by the previous equation by adding a term $\Lie_X\varphi$ to the ambient Laplacian flow. The restriction of the vector field $X$ to the submanifold is a  vector field in the normal bundle of the coassociative determined by the deformation operator of the coassociative equation. There would be many such normal vector fields and we choose the vector field with the smallest
$L^2$ norm. We then extend it to the ambient manifold $M$ using fixed data on a tubular neighbourhood of the coassociative and different choices of this extension give pairs which are related by 
diffeomorphisms of $M$. Equivalently, we may keep the evolution of the $G_2$ structures by the Laplacian flow and in that case we
modify the normal velocity of the coassociatives. The Laplacian flow equation
determines one part of this velocity. The remaining part lies in
the tangent space to the local family of coassociatives and is the
negative $L^2$ gradient of volume on that family. Together these
two parts give the compatible velocity closest to the mean curvature vector
in $L^2$. Thus, our construction selects an ordinary differential
equation on a moving finite-dimensional deformation space of the coassociatives. We
explain this interpretation in Proposition~\ref{prop-coassoc-volume-gradient}.

Having established the coupled system, we look at its analytic properties. The coassociative existence proof uses the closed-structure version of Mclean's deformation theory \cite{McL}, including the varying ambient results of Kovalev and Lotay \cite[Theorem 2.7-2.8]{KL}. These give a finite dimensional family of nearby coassociatives. An ordinary differential equation in its harmonic parameters selects the path with the prescribed normal velocity. Quantitative control of this family then allows us to lift a stable ambient flow for all time, and we use the ambient Laplacian flow stability theorem from \cite{LW19b} for this argument. We mention that Lambert and Lotay \cite{LL21} also relate Laplacian flow and mean curvature flow in the setting of semi-flat coassociative fibrations, where mean curvature flow acts on a space-like three-dimensional base. Here, we move a chosen four-dimensional coassociative in a seven-dimensional manifold, without a fibration or symmetry assumptions.
\medskip

For connections, we prescribe a smooth family of closed positive
four forms $\psi_t$. The correction in this case has two parts. A torsion term
removes the difference between ordinary Yang--Mills flow and the
negative gradient of the instanton defect energy. An elliptic
correction then accounts for $\partial_t\psi_t$. We explain that the theory is most well-behaved near an \emph{acyclic instanton} (see page 5 or \cite[Def. 3.6]{Wal}). Near an acyclic
instanton, the inverse defining this correction exists even at
connections which are not instantons. The energy identity proves
that an initially zero instanton defect remains zero, thus proving the preservation of the instanton condition along the prescribed flow of $\psi_t$.

We will see that decay of the instanton defect alone does not give convergence of the connection. To obtain convergence from connections which are not initially
instantons, we also need to control the connection itself. Acyclicity
provides this control after fixing gauge: the linearized operator has
a positive spectral gap, and the nonlinear terms can be absorbed in
the corresponding Sobolev energy estimate. The same equation therefore
attracts sufficiently close connections to an instanton. We prove this in both the cases, when a coclosed $G_2$ structure is fixed and  also along an explicitly controlled flow of coclosed $G_2$ structures. For the ambient flow of coclosed $G_2$ structures we follow  our own philosophy of using a well-studied flow of coclosed $G_2$ structures. The flow we look at is the modified Laplacian coclow of coclosed $G_2$ structures which was introduced by Grigorian \cite{Gri}.
We obtain global convergence for the coupling system when the constant parameter in the modified coflow is $0$ and the evolution takes place near a torsion-free $G_2$ structure. We explain both the settings in detail below.

Throughout the paper, the manifolds and submanifolds are compact and
have no boundary, and initial data and solutions are smooth.
Sobolev and H\"older norms use fixed smooth reference metrics and
connections. We fix $0<\alpha<1$ when using $C^{k,\alpha}$ norms.

\subsection{Coupling with coassociatives}\label{sec:first-statement}
Let $C^4 \hookrightarrow M^7$ be a coassociative submanifold. If $g$ is a Riemannian metric induced by the underlying $G_2$ structure, we use on $C$ the induced metric $g_C=i^*g$ and we denote its normal bundle by $\nu_i$. With $\perp$ denoting normal projection, the second fundamental
form and the mean curvature vector of $C$ in $M$ are
\begin{align*}
 \mathrm{II}(U,V)=(\nabla_UV)^\perp,\qquad
 H=\tr_{g_C}\mathrm{II},\qquad U,V\in TC.
\end{align*}
To describe the change of the coassociative equation under a normal
motion, we associate a two form to a normal field $V$ and then take
its exterior derivative. We set
\begin{align*}
 J_iV=i^*(V\ip\varphi),\qquad D_i:\Gamma(\nu_i)\longrightarrow d_C\Omega^2(C)\ \text{with}\ D_iV=d_CJ_iV,\qquad \theta=i^*\tau.
\end{align*}
At a coassociative, $J_i$ identifies $\nu_i$ with the anti-self-dual
two forms on $C$, for the orientation determined by
$i^*\psi=\operatorname{vol}_{g_C}$. We use the following notations:
$\pi_-=(\mathrm{id}-*_C)/2$, $\delta_C=d_C^*$ and
$\Delta_C=d_C\delta_C+\delta_Cd_C$. The Green's operator $G_C$ of $\Delta_C$ is zero
on harmonic forms, and its degree will be suppressed. The harmonic
anti-self-dual two forms are denoted by $\mathcal H^2_-(C)$ and their
dimension by $b_2^-(C)$.

Near a coassociative pair $(\varphi, C)$, $D_i$ has a minimum norm right inverse $R_i$
on exact three forms, with $D_iR_iq=q$ and $R_iq\perp\ker D_i$.
Choose the diffeomorphism-equivariant extension operator
$E_{\varphi,i}$ defined in \eqref{eq:extension}, fixing its tubular
radius and cutoff. It takes a normal vector field along $i(C)$ and
extends it to a vector field on $M$. Put (see \textsection~\ref{sec:first} for more details)
\begin{align}\label{eq:introW}
 W_{\varphi,i}=-R_i(d_C\theta+D_iH),\qquad
 X_{\varphi,i}=E_{\varphi,i}W_{\varphi,i}.
\end{align}
Here $W_{\varphi, i}$ is a normal vector field and $X_{\varphi, i}$ is its extension to the manifold $M$. Even though the subscripts $\varphi, i$ are cluttering the notations, we still keep them to emphasize the dependence of these objects on the pair $(\varphi, C)$. The term $d_C\theta$ is caused by the ambient Laplacian flow, and
$D_iH$ is caused by mean curvature flow. Thus $W$ is chosen to
cancel their sum. The equation $D_iR_iq=q$ is a right-inverse
identity on exact three forms and we remark that it does not say that $D_i$ has zero kernel. Orthogonality to that kernel chooses one solution. We
construct these operators in Section~\ref{sec:first}. Having discussed the notations, our first main theorem is the following.\footnote{Similar results for coupling of the Laplacian flow and mean curvature flow of coassociatives have been independently obtained by Jason Lotay (private communications).}

\begin{theorem}\label{thm:corrected}
Let $M^7$ be compact with a closed $G_2$ structure $\varphi_0$, and let
$i_0:C^4\hookrightarrow M$ be a compact embedded coassociative. Fix the
extension data above in \eqref{eq:introW}. The coupled system
\begin{align}\label{eq:correctedMCF}
 \partial_t\varphi_t=\Delta_{\varphi_t}\varphi_t+
       \Lie_{X_{\varphi_t,i_t}}\varphi_t,\qquad
 d\varphi_t=0,\qquad \partial_ti_t=H_t
\end{align}
has a unique smooth solution from $(\varphi_0,i_0)$ for a short time.
Every smooth solution from these initial data satisfies
$i_t^*\varphi_t=0$ and $[\varphi_t]=[\varphi_0]$ on each interval
where the pair remains in a neighborhood on which $R_i$ and the fixed
tubular extension are defined. At a coassociative pair the normal vector field
$W$ along $C$ is explicitly
\begin{align}\label{eq:introexplicitW}
 W=-2J_i^{-1}\pi_-\delta_CG_Cd_C(*_C\theta).
\end{align}
\end{theorem}

\begin{remark}
No condition on $b_2^-(C)$ is required in Theorem~\ref{thm:corrected}.
\end{remark}

The normal field $W$ has the least $L^2$ norm among those cancelling
the pullback derivative. A version of the theorem for a fixed embedding is given in Corollary~\ref{cor:fixedC}. The system \eqref{eq:correctedMCF} is a modified coupled system of the Laplacian flow and the mean curvature flow and for comparison, the unmodified system is
\begin{align}\label{eq:MCFsystem}
 \partial_t\varphi_t=\Delta_{g_t}\varphi_t,\qquad d\varphi_t=0,
 \qquad g_t=g_{\varphi_t},\qquad \partial_ti_t=H_t.
\end{align}
In this description the coassociative has normal velocity
\begin{align*}
 V=H+W=\Pi_iH-R_i d_C\theta,
\end{align*}
where $\Pi_i$ is the $L^2$ orthogonal projection onto $\ker D_i$.
The term $-R_i d_C\theta$ is the compatible normal velocity of least
$L^2$ norm and its geometric motivation is explained in Proposition~\ref{prop-coassoc-volume-gradient}. If $b_2^-(C)=0$,
the ambient form determines a unique nearby coassociative image
at each sufficiently small time.

For mean curvature flow with the ambient metric fixed, Tsai and Wang
proved stability for certain calibrated zero sections~\cite{TW18}.
They also proved stability for minimal submanifolds satisfying their
strong stability condition~\cite{TW20}.

We fix the initial parametrizations and measure closeness of embeddings
in a normal graph chart at the reference embedding. The assumption
$[\varphi_0]=[\bar\varphi]$ allows us to apply the stability theorem
for Laplacian flow.

\begin{theorem}\label{thm:first-global}
Let $(M,\bar\varphi)$ be compact with $\bar \varphi$ torsion-free and let $\bar i$ be a compact coassociative embedding. Fix admissible tubular extension data used in Theorem~\ref{thm:corrected}. There exists an integer $N$ and a $C^{N,\alpha}$ neighborhood of
$(\bar\varphi,\bar i)$ such that every smooth initial pair in this
neighborhood satisfying
\begin{align*}
 d\varphi_0=0,\qquad [\varphi_0]=[\bar\varphi],\qquad
 i_0^*\varphi_0=0
\end{align*}
has a unique global solution of \eqref{eq:correctedMCF}. The solution
converges smoothly and exponentially to a torsion-free coassociative
pair. The limiting ambient $G_2$ structure belongs to
$\mathrm{Diff}_0(M)\cdot\bar\varphi$. The limiting coassociative may
vary in its local deformation space.
\end{theorem}

\subsection{Coupling with instantons}\label{sec:second-statement}
Let $P\to M$ be a principal $G$-bundle with $G$ compact. Fix an
$\operatorname{Ad}(G)$-invariant positive definite inner product on
$\mathfrak g$, and set $E=\ad P$. We write $\Omega^k(M,E)$ for the
$E$-valued $k$-forms and use the fixed inner product in the adjoints
and norms of these forms. Let $\psi$ be a closed positive four-form,
with associated metric $g_\psi$ and dual positive three-form
$\varphi=*_\psi\psi$. For a connection $A$ on $P$ with curvature
$F_A$ and for $(f,a)\in\Omega^0(M,E)\oplus\Omega^1(M,E)$, set
\begin{align}
 \chi_{A,\psi}&=*(F_A\wedge\psi),&
 \PP_{A,\psi}a&=*(d_Aa\wedge\psi),\label{eq:introdefect}\\
 \DD_{A,\psi}(f,a)&=(d_A^*a,d_Af+\PP_{A,\psi}a).
 \label{eq:introD}
\end{align}

\begin{definition}
A connection $A$ on $P\rightarrow M$ is called a $G_2$-instanton if its curvature $F_A$ satisfies
\begin{align*}
F_A\wedge \psi=0.
\end{align*}
\end{definition}

We call $\chi_{A,\psi}$ the instanton defect. With $\psi$ fixed, its
derivative in the direction $a$ is $\PP_{A,\psi}a$. At an instanton,
\begin{align*}
 \PP_{A,\psi}(d_Aq)=*([F_A,q]\wedge\psi)=0
 \qquad\text{for }q\in\Omega^0(M,E).
\end{align*}
We impose the Coulomb condition $d_A^*a=0$ to fix this gauge freedom
and add the zero-form $f$ to obtain the elliptic operator
$\DD_{A,\psi}$. An instanton is \textbf{acyclic} if $\DD_{A,\psi}$ is
invertible~\cite[Definition 3.6]{Wal}. We identify the kernel of the operator $\mathcal D$ at an instanton with the
direct sum of the parallel sections of $E$ and the infinitesimal
instanton deformations in Coulomb gauge in Lemma~\ref{lem:D}.

We use the inverse operator $\mathcal D^{-1}$ in the correction below but for instance this does not apply on the full adjoint bundle of a $U(r)$-bundle as this inverse never exists simply because $\sqrt{-1}\operatorname{Id}$ is parallel for every connection $A$. In particular,
\begin{align*}
 \DD_{A,\psi}(\sqrt{-1}\operatorname{Id},0)=0.
\end{align*}
Thus, a fixed-determinant version would assume acyclicity on the traceless adjoint bundle. The fixed determinant connection must also satisfy the instanton equation for each $\psi_t$.

For a closed four-form $Q$ and a pair $(A,\psi)$ for which
$\DD_{A,\psi}$ is invertible, set
\begin{align}\label{eq:introU}
 (f,\UU_{A,\psi}(Q))=-\DD_{A,\psi}^{-1}(0,*(F_A\wedge Q)).
\end{align}
The operator stays invertible near an acyclic instanton, so this
formula also applies to nearby connections with nonzero defect.
The Bianchi identity and $d\psi=0$ give $d_A(F_A\wedge\psi)=0$,
and hence $d_A^*\chi_{A,\psi}=0$. On a compact manifold without
boundary, integration by parts gives
\begin{align*}
 \langle\chi_{A,\psi},d_Af\rangle_{L^2}
 =\langle d_A^*\chi_{A,\psi},f\rangle_{L^2}=0.
\end{align*}

\begin{theorem}\label{thm:transport}
Let $M^7$ be compact without boundary and let $\psi_t$ be a prescribed
smooth family of closed positive four-forms with
$\varphi_t=*_{t}\psi_t$ and $Q_t=\partial_t\psi_t$.
If $A_0$ is an acyclic instanton for $\psi_0$, then
\begin{align}\label{eq:coupledYM}
 \partial_tA=-d_A^{*_{t}}F_A
              +*_{t}(F_A\wedge d\varphi_t)+\UU_{A,\psi_t}(Q_t)
\end{align}
has a unique smooth solution for short time with initial value $A_0$.
Every smooth solution from $A_0$ satisfies $F_{A_t}\wedge\psi_t=0$
as long as $\DD_{A_t,\psi_t}$ remains invertible. Along any such
solution, the section $f$ in \eqref{eq:introU} vanishes and
\begin{align}\label{eq:horizontaltransport}
 d_A^*\partial_tA=0,\qquad
 \PP_{A,\psi}\partial_tA=-*(F_A\wedge Q).
\end{align}
\end{theorem}

On a family of nearby instantons which varies with time, the Yang--Mills and torsion terms in
\eqref{eq:coupledYM} cancel. The remaining equation
$\partial_tA=\UU_{A,\psi_t}(Q_t)$ gives the Coulomb lift of the branch.
Section~\ref{sec:connection-examples} gives examples in which the gauge class of instantons is moving, i.e, the velocity in \eqref{eq:coupledYM} is nonzero.

If the coclosed $G_2$ structure is not moving, we have the following 

\begin{theorem}\label{thm:static-stability}
Let $A_\star$ be an acyclic instanton for a closed positive four form
$\psi_\star$, and put $\varphi_\star=*_{\psi_\star}\psi_\star$.
Fix any integer $s\geq8$, every sufficiently close smooth initial
connection in the Sobolev space $H^s$ has a unique global solution of
\begin{align*}
 \partial_tA=-d_A^*F_A+*(F_A\wedge d\varphi_\star).
\end{align*}
It converges smoothly and exponentially to a gauge transform of
$A_\star$. The initial connection need not be an instanton.
\end{theorem}

%%%%%%%%%%%%%%%%%%%%%%%%%%%%%%%%%%%%%%%%

On a fixed torsion-free background, the equation in the preceding
theorem is the ordinary Yang--Mills flow. Global existence and
convergence for this flow near a local minimum of the Yang--Mills
energy were proved by Feehan~\cite[Thm. 6]{Feehan16} on closed
manifolds of arbitrary dimension. On a fixed coclosed background, the equation is the negative
gradient flow of the defect energy described in Section~\ref{sec:second}.
Its relation to the $\Omega$--Yang--Mills functional, with
$\Omega=-\varphi_\star$, is explained there; see also \cite{CW}.
Acyclicity makes the Hessian positive in Coulomb gauge, so the
convergence mechanism is the usual stability mechanism near a
nondegenerate minimum modulo gauge. We include this case as a
consequence of Theorem~\ref{thm:moving-stability}, which also
treats the additional terms produced by a moving ambient form and
a moving reference instanton.

%Our statement also allows coclosedbackgrounds with torsion, for which the torsion term in the equationneed not vanish. A moving-background version of the preceding theorem is given in Theorem~\ref{thm:moving-stability} with the uniform smallness hypotheses stated there.

\medskip

We now related our result to a well-studied flow of coclosed $G_2$-structures called the modified Laplacian coflow of coclosed $G_2$-structures introduced by Grigorian \cite{Gri}. We recall that the full torsion tensor of $\varphi$ is defined by
$\nabla_i\varphi_{abc}=T_i{}^q\psi_{qabc}$ and
$\tr T=g^{ij}T_{ij}$. For a constant $c\in\mathbb R$, which we call the coflow parameter,
we take the ambient velocity in Theorem~\ref{thm:transport} to be 
\begin{align}\label{eq:modifiedcoflow}
 \QQ_c(\psi)=\Delta_\psi\psi-2d\big((c+\tr T_\psi)\varphi\big).
\end{align}
This is precisely Grigorian's modified coflow in our convention; the conversion
and the parameter-zero case are discussed in Section~\ref{sec:coflow}. We have the following

\begin{corollary}\label{cor:coflow}
Under the initial hypotheses of Theorem~\ref{thm:transport}, fix
$c\in \mathbb R$ and couple \eqref{eq:coupledYM} to
$\partial_t\psi=\QQ_c(\psi)$, with $Q=\QQ_c(\psi)$ and $d\psi=0$.
This system has a unique smooth solution from $(A_0,\psi_0)$ for a
short time. Every smooth solution from these initial data preserves
$F_A\wedge\psi=0$ on each interval where $\DD_{A_t,\psi_t}$ is
invertible.
\end{corollary}

Note that the parameter $c$ in Grigorian's modified coflow is a constant. Since torsion-free $G_2$-structures are stationary points of the flow, one can ask whether they are dynamically stable along the flow. Such a result was proved by Bedulli and Vezzoni in \cite{BVStab20} for the choice $c=0$ which in a way suggests that $c=0$ might be a suitbale choice in the flow. The global theorem below uses $c=0$. The local result above does not
by itself imply global stability for positive $c$. In fact, Proposition~\ref{prop:positive-c-instability} gives backgrounds with
positive linear growth rates when $c>0$ which shows that $c>0$ can introduce unstable modes. It does not assert
instability for every positive parameter on every background.

\begin{theorem}\label{thm:coflow-global}
Let $(M,\bar\psi)$ be compact and torsion-free, and let $A_\star$ be an acyclic $\bar\psi$-instanton. Fix $0<\alpha<1$. There are a finite integer $N$ and a constant $\delta>0$ such that the following holds. Suppose that $\psi_0$ and $A_0$ are smooth,
that $\psi_0$ is closed and positive, that $[\psi_0]=[\bar\psi]$,
and that
\begin{align*}
 \|\psi_0-\bar\psi\|_{C^{N,\alpha}}
 +\|A_0-A_\star\|_{H^8}<\delta.
\end{align*}
These norms use the fixed metric of $\bar\psi$ and the reference
connection $A_\star$ and the connection difference is measured in a
fixed initial gauge. Then the system
\begin{align*}
 \partial_t\psi&=Q_0(\psi),\\
 \partial_tA&=-d_A^{*t}F_A+*_t(F_A\wedge d\varphi_t)
                  +\UU_{A,\psi_t}(\partial_t\psi_t)
\end{align*}
has a unique global smooth solution. The operator
$\DD_{A_t,\psi_t}$ remains invertible. The solution converges
smoothly and exponentially to a torsion-free structure together
with an acyclic instanton. The initial connection may have nonzero
instanton defect. If it is initially an instanton, it remains one.
\end{theorem}
In this theorem, one may take $N=r+8$, where $r\geq\max\{r_0,10\}$ and $r_0$
is an ambient differentiability order sufficient for the transport
estimates in the proof of Theorem~\ref{thm:global-instanton-transport}
to control the transported connection in $C^9$. These orders are
fixed before $\delta$ is chosen.

\subsection{Strategy of proofs}

We now explain the strategy of proving the above mentioned theorems.

\medskip\noindent\textbf{Preservation and existence.}
For coassociatives, the error takes values in exact three forms.
Hodge theory gives a right inverse of $D_i$ on this space. Its
kernel consists of the harmonic normal deformations; projecting
away from that kernel selects the solution of smallest norm.
If $R_0$ is the right inverse at the initial coassociative then because $D_iR_0$ is close to the identity operator on the space of exact $3$-form, it is invertible and lets us define a right inverse for $D_i$ and hence $W$ at nearby pairs $(\varphi, i)$ even when $i*\varphi \neq 0$. We therefore
define the correction before assuming that the evolving pair is
coassociative. The coassociative equation gives
$\partial_t(i^*\varphi)=d_C\theta+D_iH+D_iW=0$ for every 
solution in the chosen neighborhood.

We construct a solution by first solving the unmodified Laplacian flow.
The varying coassociative deformation theorem \cite{McL,KL} gives
a local family of coassociatives for the evolving $G_2$ structures.
The parameters in the kernel satisfy an ODEs on a finite-dimensional space. Its solution selects the required normal velocity. A diffeomorphism of $M$ then moves the correction into
the ambient Laplacian flow equation, leaving the submanifold velocity equal to mean curvature. Reversing these steps proves uniqueness of the solutions.

\medskip

For instantons, we differentiate the equation for the defect $\chi$ and use $d_A^*\chi=0$. When we differentiate $\chi= *_\psi(F_A\wedge\psi)$, the correction
cancels $*_\psi(F_A\wedge Q)$ and leaves a term $-d_Af$ in the
defect equation. The Bianchi identity and $d\psi=0$ give
$d_A^*\chi=0$, so the $L^2$ pairing of $d_Af$ with $\chi$ vanishes
by integration by parts. The energy estimate and Gronwall's inequality
then give $\chi(t)=0$ if $\chi(0)=0$.

More precisely, to construct nearby instantons in Lemma~\ref{lem:IFT}, we solve
Walpuski's augmented equation~\cite[eq.~(2.4)]{Wal17} in a fixed
Coulomb slice:
\begin{align*}
 *_\psi(F_A\wedge\psi)+d_A\xi=0,\qquad
 d_{A_0}^{*_{g_{\psi_0}}}(A-A_0)=0.
\end{align*}
Here $\xi$ is an adjoint-valued zero-form. The linearization in
$(\xi,A)$ at $(0,A_0,\psi_0)$ is $\DD_{A_0,\psi_0}$, which is
invertible by the acyclicity assumption. The implicit function theorem gives solutions for nearby closed positive forms $\psi$, by the argument of \cite[Proposition~2.11]{Wal17}. Applying $d_A^*$ to the first equation gives $d_A^*d_A\xi=0$. Integration by parts yields
$d_A\xi=0$, so $A$ satisfies the instanton equation.

For initial connections with nonzero defect, we must solve the full
connection equation. Following \cite[Section~2.1]{Hong}, we first add
a gauge term so that its highest-order part becomes a heat operator.
We also check that the nonlocal correction $\UU$ preserves the
Sobolev regularity of the connection. Indeed, forming the curvature
uses one derivative, while applying $\DD^{-1}$ restores that
derivative. In particular, $\UU$ maps the $H^6$ neighborhood used in
the existence proof into $H^6$, with local Lipschitz estimates.
These estimates allow us to solve the modified equation for a short
time. We then apply a time-dependent gauge transformation that
cancels the added gauge term and gives a solution of the original
equation.

\medskip\noindent\textbf{Long-time behavior.}
For understanding the long time behaviour of the coupled system for coassociatives, we control the graph charts and right inverses
uniformly near a fixed pair. We show that
the evolving pair $(\varphi_t, i_t)$ stays in the neighborhood where the local
construction applies. We describe nearby coassociatives as normal
graphs over the reference submanifold and obtain uniform estimates
for these parametrizations and the right inverses $R_i$. 
%%The velocities of the deformation parameters are bounded by suitable
%norms of the ambient velocity and torsion. Integrating this estimate
%bounds the total change of the parameters. If this bound is
%sufficiently small, the parameters cannot leave the region where
%the construction is defined. A small velocity at each time would
%not be enough, since even a small constant velocity can eventually
%produce a large displacement.
Near a torsion-free structure, the
stability theorem of \cite{LW19b} keeps the ambient form close to
the reference torsion-free $G_2$ structure, and we extract the required integral bounds from its proof.

For connections, we compare the solution with a nearby instanton
which moves along the evolving ambient $4$-forms. After choosing
gauge, we estimate the difference between the two connections.
Acyclicity gives a positive spectral gap for the second-order
operator in the linearized equation, so the linearized flow
decreases this difference. We show that sufficiently small
perturbations and the explicit bounds on the moving $\psi_t$
allow us to control the additional terms. An estimate in one
fixed Sobolev norm then gives exponential decay and keeps the
solution in the neighborhood where the correction is defined.
Parabolic smoothing gives exponential decay of every higher
spatial derivative after positive time. The moving-background
assumptions are stated in Theorem~\ref{thm:moving-stability}.
If the evolving instanton also converges, the original
connection converges to its limit up to a gauge transformation.

For the ordinary Laplacian flow,
Lemma~\ref{lemma-ambient-laplacian} uses the stability proof of
Lotay and Wei \cite{LW19b} to bound the distance between the evolving $G_2$ structures and the reference torsion-free one as well as the the time integrals of its velocity
and torsion. These bounds keep the coassociative in its
deformation chart for all time.

For the modified coflow with $c=0$, Bedulli and Vezzoni
\cite[Theorem~3.4]{BVStab20} prove stability near a torsion-free
structure under smallness in every $C^k$ norm. In
Lemma~\ref{lemma-ambient-coflow}, we prove the finite-order
estimates needed here by a direct energy argument. We show that for each fixed
spatial order $r$, it is enough to assume smallness in
$C^{r+8,\alpha}$. 
%The initial form is smooth, closed, positive, and cohomologous to the torsion-free reference form. Its higher derivatives need not be small.

The same initial norm controls the distance from the reference
throughout the flow and gives an exponentially decaying bound for
the ambient velocity in $C^{r,\alpha}$. In particular, it bounds
\begin{align*}
 \int_0^\infty\|\partial_t\psi_t\|_{C^{r,\alpha}}\,dt.
\end{align*}
These quantitative estimates keep the transported instanton close
to the reference connection and verify the moving-background
hypotheses of Theorem~\ref{thm:moving-stability}. This gives the
global coupled result in Theorem~\ref{thm:coflow-global}.

\medskip\noindent\textbf{Organization of the paper.}
We give basics of $G_2$ geometry and fix our notations and conventions in Section~\ref{sec:basics}.
Section~\ref{sec:first} treats the coassociative system. We prove our short-time and long-time existence results in that section which includes the proofs of Theorems~\ref{thm:corrected} and \ref{thm:first-global}. We also give an explicit example on the flat torus in Proposition~\ref{prop:torus-example} which shows that a straightforward coupling of the Laplacian flow and the mean curvature flow does not preserve the coassociative condition. We discuss the coupled flow system of connections and the $4$-form in  Section~\ref{sec:second}. This sections contains the local
theory and nonlinear stability of the coupled system, and the application to modified coflow. It also gives examples with nonzero correction and the
compatibility condition at a nonacyclic instanton.
We finish the paper with some questions which we find interesting but which we do not treat in this paper in Section~\ref{sec:questions}.

\medskip

\noindent
\textbf{Acknowledgments.}
SP would like to thank Vivek Shende for his immense support and
encouragement throughout the project. SP was supported by VILLUMFONDEN, VILLUM Investigator grant 37814 and Danish National Research Foundation grant DNRF157. SD would like to thank Jason Lotay for several discussions on the topics of this paper and for sharing his ideas on coupled flows. He would also like to thank Thomas Walpuski for discussions on gauge theory on $G_2$ manifolds and the Yang-Mills flow. SD acknowledges support by the Deutsche Forschungsgemeinschaft (DFG, German Research Foundation) under Germany’s Excellence Strategy - EXC 2121 "Quantum Universe" - 390833306. This project started when SP was a visiting PhD student of SD at the University of Hamburg in 2024-25. SP would like to thank the Department of Mathematics at the University of Hamburg for support and for providing excellent working conditions.

\medspace

\noindent
\textbf{AI usage.}
AI tools were used for language editing, polishing and for literature review. They were also used as an assistant for proof-checking. The authors take full responsibility for correctness and final contents of the paper.

\section{Basics on \texorpdfstring{$G_2$}{G2} structures}\label{sec:basics}
We first fix the algebraic signs used throughout the paper. The
metric, the Hodge star, and the decomposition of forms all depend on
the $G_2$ structure. In particular, differentiating a Hodge star means
differentiating the metric as well. We use the model
\begin{align}\label{eq:flat}
 \varphi_{\mathrm{flat}}=dy^{123}+dy^1\wedge(dx^{12}+dx^{34})
          +dy^2\wedge(dx^{13}-dx^{24})-dy^3\wedge(dx^{14}+dx^{23}).
\end{align}
The ambient orientation is $-dy^{123}\wedge dx^{1234}$, and the
coassociative $x$ plane is oriented by $-dx^{1234}$. We use the sign conventions in \cite{Kar}. We use the
irreducible summands
$\Lambda^2_7=\{V\ip\varphi\}$ and
$\Lambda^2_{14}=\{\beta:\beta\wedge\psi=0\}$ of the two form bundle,
and write $\Omega^2_{14}(M)$ for smooth sections of the latter.
The corresponding decomposition of three forms is
$\Lambda^3=\Lambda^3_1\oplus\Lambda^3_7\oplus\Lambda^3_{27}$,
where $\Lambda^3_1=\mathbb R\varphi$,
$\Lambda^3_7=\{V\ip\psi\}$ and $\Lambda^3_{27}$ is their
orthogonal complement. Put $\Lambda^4_j=*\Lambda^3_j$ and let
$\pi_j$ denote the corresponding projection on four forms when it
is applied to a four form. We will use the algebraic variation
formula \cite{Kar}
\begin{align}\label{eq:dual-variation}
 \dot\varphi
 =*\bigl(\tfrac34\pi_1\dot\psi+\pi_7\dot\psi-\pi_{27}\dot\psi\bigr).
\end{align}
In particular, a variation in $\Lambda^4_{27}$ gives
$\dot\varphi=-*\dot\psi$. We will use this sign when linearizing
the coflow. The map from $\psi$ to its dual three form is nonlinear and \eqref{eq:dual-variation} is its derivative, including the change
of the Hodge star.

For tangent vectors $u,v$ and a two form $\beta$, our conventions give
\begin{align}\label{eq:mainconv}
 (u\ip\varphi)\wedge(v\ip\varphi)\wedge\varphi
       =-6g(u,v)\vol_g,\qquad
 *(\varphi\wedge\beta)=\begin{cases}-2\beta,&\beta\in\Lambda^2_7,\\
 \beta,&\beta\in\Lambda^2_{14}.
 \end{cases}
\end{align}
The torsion of a $G_2$-structure is a $2$-tensor $T$ which is essentially $\nabla \varphi$. The explicit expression for $T$ in a local orthonormal frame is 
\begin{align*}
T_{ij}=\frac{1}{24}\nabla_i\varphi_{abc}\psi_j^{\, \, abc}.
\end{align*}
For a closed $G_2$ structure, the torsion form $\tau$ satisfies
\begin{align}\label{eq:maintau}
 d\psi=-\tau\wedge\varphi,\qquad \tau\in\Omega^2_{14}(M).
\end{align}
In the notation of \cite{Kar}, $\tau=-\tau_2$ for a closed structure.
The identity $d\psi=-\tau\wedge\varphi$ will turn the first variation
of coassociative volume into a formula involving the restriction of
$\tau$. 

The cross product is defined by $g(U\times V,W)=\varphi(U,V,W)$,
and our convention for pointwise norms is
$|\omega|^2=\frac1{k!}\omega_{i_1\cdots i_k}\omega^{i_1\cdots i_k}$.

\medskip\noindent\textbf{Coassociatives and instantons.}
At a coassociative, the map $J_i$ introduced in
Section~\ref{sec:first-statement} takes values in $\Lambda^2_-T^*C$ with norm $|J_iV|^2=2|V|^2$, where the anti-self-dual notation follows
from our chosen orientation. For nearby embeddings, we retain the
orientation of $C$ induced by the initial coassociative embedding, and
Lemma~\ref{lem:H} computes $J_i(H)$ when $d\varphi=0$.

For the bundle $E$ in Section~\ref{sec:second-statement}, the instanton
condition is equivalently $F_A\in\Omega^2_{14}(M,E)$.
Under \eqref{eq:mainconv} this implies $*(\varphi\wedge F_A)=F_A$.
Thus an instanton can have nonzero curvature and the instanton equation only requires the $\Lambda^2_7$ part of the curvature to vanish. This
will matter when interpreting the defect energy: it controls that
component, not the full curvature.

\medskip\noindent\textbf{Differentiating a pullback.}
A dot or $\partial_t$ denotes time differentiation on the fixed
domain $C$. For smooth $\varphi_t$ and embeddings $i_t$ with velocity
$V_t=\partial_ti_t$, let $\gamma_t=i_t^*\varphi_t$. Then
\begin{align}\label{eq:basicpullback}
 \dot\gamma=i^*\dot\varphi+d_Ci^*(V\ip\varphi)
                         +i^*(V\ip d\varphi).
\end{align}
There are two contributions in this expression. The $3$-form changes with time, and the embedding moves as well. To compute the second term, extend $V$ locally to an ambient vector field and use Cartan's magical formula $\Lie_V=d(V\ip\,\cdot)+V\ip d$. After pullback, the result depends only on $V$ along the embedding, so it is independent of the chosen extension. If $d\varphi=0$, the last term vanishes. The two remaining terms are the ones we must cancel in order to preserve the coassociative construction.

For local coordinate one forms we use
$dx^{i_1\cdots i_k}=dx^{i_1}\wedge\cdots\wedge dx^{i_k}$, and
similarly $e^{i_1\cdots i_k}$ for a coframe. 

\section{The first coupled system: coassociatives}\label{sec:first}
\subsection{Correction and preservation}\label{sec:ambientcorrection}
We first calculate how the coassociative equation changes when the
ambient $G_2$ structure evolves by the Laplacian flow and the submanifold moves by the mean curvature flow. This calculation will identify the term that we need to
cancel. We then construct a right inverse that produces the correction.
The existence of a solution to the corrected equation is a separate
question, which we address in Section~\ref{sec:first-existence}.

\medskip

We begin with the relation between mean curvature and torsion. This
relation is part of the general theory developed by Ball and
Madnick~\cite{BM}. For closed $G_2$ structures, Pacini and
Raffero~\cite[Lemma~2.3]{PR24} give the first variation of
coassociative volume in terms of the ambient torsion. We include the
calculation in our conventions because the signs and normalization
enter the correction below. 

\begin{lemma}\label{lem:H}
For a compact coassociative $i:C^4\hookrightarrow(M^7,\varphi)$ with $d\varphi=0$, let $\theta=i^*\tau$ where $\tau \in \Omega^2_{14}$ is the torsion form. The map
$J_i:\nu C\to\Lambda^2_-T^*C$, $J_i(V)=i^*(V\ip\varphi)$, satisfies
\begin{align}\label{eq:JH}
 \langle J_i(V),J_i(W)\rangle=2g(V,W),\qquad J_i(H)=-2\pi_-\theta.
\end{align}
\end{lemma}
\begin{proof}
We first check the pointwise statement. Choose an adapted orthonormal
coframe and write $\varphi$ as in \eqref{eq:flat}, with the orientations
fixed in Section~\ref{sec:basics}. Contraction with the three unit normal
vectors gives three mutually orthogonal anti-self-dual two forms, each
of squared norm two. By linearity, their pairings give the first identity
in \eqref{eq:JH}. In particular, $J_i$ is an isomorphism from the normal
bundle to the bundle of anti-self-dual two forms and the factor $2$ records
its normalization.

We next vary the embedding while keeping the ambient forms and metric
fixed. Let $V$ be any smooth normal field and choose a variation $i_s$
with initial velocity $V$. The comass inequality says pointwise that
$i_s^*\psi\leq\vol_{i_s^*g}$. At the coassociative embedding $i_0=i$,
the two sides agree. Their difference is therefore a nonnegative
function of $s$ that vanishes at zero, so its derivative at zero is
zero. We can thus compute the variation of the volume by differentiating
$\int_Ci_s^*\psi$.
The first variation of volume gives the first equality below. Cartan's
formula and Stokes' theorem give the second with the exact term integrating
to zero because $C$ is compact and has no boundary:
\begin{align*}
 -\int_Cg(H,V)\,d\mu
   =\left.\frac{d}{ds}\right|_0\int_Ci_s^*\psi
   =\int_Ci^*(V\ip d\psi).
\end{align*}
We now use \eqref{eq:maintau}. Since $i^*\varphi=0$, contraction of
$d\psi=-\tau\wedge\varphi$ gives
$i^*(V\ip d\psi)=-\theta\wedge J_i(V)$. For an anti-self-dual two form
$\alpha$, we have $\theta\wedge\alpha=-\langle\pi_-\theta,\alpha\rangle d\mu$.
Consequently the integrand on the right is
$\langle\pi_-\theta,J_i(V)\rangle d\mu$.
We may choose $V$ supported in an arbitrarily small neighborhood, so the
integral identity gives
$-g(H,V)=\langle\pi_-\theta,J_i(V)\rangle$ pointwise.
%The first identity in \eqref{eq:JH} rewrites the left side as
%$-\tfrac12\langle J_i(H),J_i(V)\rangle$. 
As $J_i(V)$ ranges over all
anti-self-dual two forms, we obtain $J_i(H)=-2\pi_-\theta$. This completes the proof of the lemma.
\end{proof}

\begin{proposition}\label{prop:unmodified-pullback}
For a smooth solution of the unmodified coupled system \eqref{eq:MCFsystem}, at every time
\begin{align}
 \partial_t(i_t^*\varphi_t)
      =d_C\bigl(i_t^*\tau_t+i_t^*(H_t\ip\varphi_t)\bigr).
\end{align}
At a time when $i_t$ is coassociative, this becomes
\begin{align}
 \partial_t(i_t^*\varphi_t)=d_C(*_Ci_t^*\tau_t).
\end{align}
At such a time, the coassociative preservation condition
$\partial_t(i_t^*\varphi_t)=0$ is equivalent to
$\delta_C(i_t^*\tau_t)=0$.
\end{proposition}

\begin{proof}
We differentiate both the ambient form and the embedding. Since the $G_2$-structures remain closed along the Laplacian flow, we get that the final contraction term in \eqref{eq:basicpullback} vanishes, and the Laplacian flow equation
$\partial_t\varphi=d\tau$ gives
\begin{align*}
 \partial_t(i_t^*\varphi_t)
 =d_C(i_t^*\tau_t)+d_C\bigl(i_t^*(H_t\ip\varphi_t)\bigr).
\end{align*}
If $i_t$ is coassociative at the time under consideration, we can use
\eqref{eq:JH}. Since $2\pi_-=I-*_C$, we obtain
$\theta+J_i(H)=\theta-2\pi_-\theta=*_C\theta$.
This proves the second formula. Finally, $*_C$ is an isomorphism and
$\delta_C=-*_Cd_C*_C$ on two forms. Thus $d_C(*_C\theta)=0$ is equivalent
to $\delta_C\theta=0$.
\end{proof}
The proposition gives condtions for a preservation condition. It is
not automatically satisfied, as the compact example below will show.
To obtain preservation, we add a term that cancels this exact three
form. We will do this by changing the ambient form through a Lie
derivative whose vector field is prescribed along the submanifold.

We first specify the spaces on which the cancellation equation is
solved. We suppress the
dependence on the ambient $G_2$ structure $\varphi$ in the notation $J_i$, $D_i$, and $R_i$.

\begin{lemma}\label{lem:subinverse}
If $i$ is coassociative, $D_i$ is surjective onto exact three forms.
Its kernel is $J_i^{-1}\mathcal H^2_-(C)$, and its minimum norm right
inverse is
\begin{align}\label{eq:subexplicit}
 R_iq=2J_i^{-1}\pi_-\delta_CG_Cq,\qquad q\in d_C\Omega^2(C),
\end{align}
where $G_C$ is the Green operator on three forms. In particular,
$D_iR_iq=q$ and $R_iq\perp\ker D_i$.
\end{lemma}
\begin{proof}
We construct an anti-self-dual primitive of $q$. Set
$\eta=\delta_CG_Cq$. Since $q$ is exact it has zero harmonic
projection, and also that $d_Cq=0$. Since the Green operator commutes
with $d_C$, we have $d_CG_Cq=G_Cd_Cq=0$. The identity
$\Delta_CG_Cq=q$ now reduces to $d_C\eta=q$.
Also $\delta_C\eta=\delta_C^2G_Cq=0$, which is equivalent to
$d_C(*_C\eta)=0$. Subtracting $*_C\eta$ therefore does not change the
exterior derivative, and gives an anti-self-dual form
\begin{align*}
 d_C(2\pi_-\eta)=d_C(\eta-*_C\eta)=q.
\end{align*}
Applying $J_i^{-1}$ proves that the formulas is a right
inverse. To determine the kernel, suppose $D_iV=0$. Then $J_iV$ is
closed and anti-self-dual. Anti-self-duality makes a closed two form
co-closed as well, so $J_iV$ is harmonic. The converse is immediate.

We next check that the right inverse takes values orthogonal to this
kernel. If $h$ is harmonic and anti-self-dual, integration by parts gives
$\langle2\pi_-\delta_CG_Cq,h\rangle_{L^2}
=2\langle G_Cq,d_Ch\rangle_{L^2}=0$.
The metric identity in \eqref{eq:JH} transfers this orthogonality to
normal fields. Every other solution of $D_iV=q$ differs from $R_iq$ by
an element $k$ of the kernel. Orthogonality gives
$\|R_iq+k\|_2^2=\|R_iq\|_2^2+\|k\|_2^2$, proving both minimality and
uniqueness of the minimum norm solution. In particular, a nonzero
$b_2^-(C)$ is the degree of freedom for the solution.

For later use, we also express this right inverse through an elliptic
operator. The normalization of $J_i$ gives $J_iJ_i^*=2\pi_-$, so on
exact three forms
\begin{align}
 D_iD_i^*q=2d_C\pi_-\delta_Cq=d_C\delta_Cq=\Delta_Cq.
\end{align}
Here $d_C*_C\delta_Cq=0$ follows from $\delta_C^2q=0$, and $d_Cq=0$
gives the last equality. There are no nonzero exact harmonic forms.
Thus $D_iD_i^*$ is invertible on the exact forms, and the same minimum
norm inverse is $R_i=D_i^*(D_iD_i^*)^{-1}$.
\end{proof}

We must define the correction before assuming that a solution preserves
coassociativity. We therefore extend the inverse $R_i$ in \eqref{eq:subexplicit} to an open neighborhood of the initial pair $(\varphi_0, i_0)$, including embeddings that are not necessarily coassociative.
The next lemma supplies this extension and the estimates needed later.
Its proof first constructs some right inverse, then removes its kernel
component to obtain the minimum norm inverse.

\begin{lemma}\label{lem:subopen}
For smooth pairs $(\varphi,i)$ sufficiently close to a coassociative
pair in the $C^2\times C^3$ topology,
with $d\varphi=0$, the operator $D_i$ remains surjective onto exact
three forms and has a finite-dimensional kernel of dimension $b_2^-(C)$.
Its minimum norm right inverse
\begin{align}\label{eq:subRnear}
 R_i=D_i^*(D_iD_i^*)^{-1}
\end{align}
is defined on exact three forms and depends smoothly on smooth
finite-dimensional families of such pairs.
\end{lemma}
\begin{proof}
Fix a reference coassociative pair $(\varphi_0,i_0)$, and write
$\nu_0=\nu_{\varphi_0,i_0}$ and $D_0=D_{\varphi_0,i_0}$. Let $R_0$ be
the right inverse in Lemma~\ref{lem:subinverse}. We will compare nearby
operators on fixed spaces, so that smallness of their difference has a
precise meaning.

We first check the symbol. We use Hamilton's idea of computing the symbol of an operator by replacing the derivatives by the Fourier transform variable. We use the convention
$\partial_j\mapsto\sqrt{-1}\xi_j$. Thus, the symbol of the operator $D_i$ is 
$\sigma_\xi(D_i)V=\sqrt{-1}\xi\wedge J_iV$ at a coassociative pair.
Take $\xi=|\xi|e^1\ne0$. If an anti-self-dual two form $\alpha$
satisfies $\xi\wedge\alpha=0$, then $\alpha=e^1\wedge\lambda$.
But $*_C\alpha$ has no $e^1$ component, whereas
$\alpha=-*_C\alpha$, hence $\alpha=0$.
The symbol is therefore injective. Its image is contained in
$\ker(\xi\wedge:\Lambda^3T^*C\to\Lambda^4T^*C)$, which has real
dimension three, the same as the normal bundle and hence it is an isomorphism onto that space.
On the compact unit cosphere, the inverse symbol has a uniform bound.
The same property therefore holds for sufficiently nearby pairs. In
particular, $D_i^*D_i$ is elliptic there.

Use reference norms on
\begin{align*}
 \mathscr X=H^1(C,\nu_0),\qquad
 Y=\overline{d_C\Omega^2(C)}^{\,L^2(g_{C,0})}.
\end{align*}
Hodge decomposition shows that $Y$ is the closed space of exact
$L^2$ three forms. We identify the normal bundles as follows. Transport a
vector from $i_0(p)$ to $i(p)$ along the short geodesic of a fixed
background metric, and then project it onto $\nu_{i,p}$. For nearby
pairs this gives a bundle isomorphism $S_i:\nu_0\to\nu_i$, depending
smoothly on the pair and satisfying $S_0=I$. Thus we can regard every
nearby operator as an operator with the same domain and target and we have
\begin{align*}
 \widehat D_i=D_iS_i:\mathscr X\to Y,\qquad
 \|\widehat D_i-D_0\|_{H^1\to L^2}
       \leq C\|J_iS_i-J_0\|_{C^1}.
\end{align*}
The estimate follows by differentiating $J_iS_i-J_0$ once, since
$\widehat D_i-D_0=d_C(J_iS_i-J_0)$. Make the pair close enough that
$\|(\widehat D_i-D_0)R_0\|_{Y\to Y}<1$. We can then define
\begin{align*}
 \widehat R_i=R_0\bigl[I_Y+(\widehat D_i-D_0)R_0\bigr]^{-1},\qquad
 \widetilde R_i=S_i\widehat R_i,\qquad D_i\widetilde R_i=I_Y.
\end{align*}
The inverse in brackets is the convergent Neumann series for
$I_Y+(\widehat D_i-D_0)R_0$. Multiplying by $\widehat D_i$ shows directly
that $\widehat D_i\widehat R_i=I_Y$, proving surjectivity.

\medskip

The inverse found above need not yet have minimum norm. We now do that.
The fixed decomposition
$\mathscr X=\ker D_0\oplus\operatorname{im}R_0$ also describes the new
kernel. The restriction of $\widehat D_i$ to $\operatorname{im}R_0$ is
an isomorphism onto $Y$. Thus, for each $k\in\ker D_0$, there is a
unique component in $\operatorname{im}R_0$ that cancels
$\widehat D_i k$, namely $-\widehat R_i\widehat D_i k$. Consequently
\begin{align*}
 \ker\widehat D_i
 =\{k-\widehat R_i\widehat D_i k:k\in\ker D_0\}.
\end{align*}
The $\ker D_0$ component is exactly $k$, so
this parametrization is injective. It is also surjective by the fixed
decomposition. The kernel dimension therefore remains $b_2^-(C)$.
Let $\Pi_i$ be the current $L^2$ projection onto $\ker D_i$. Subtracting
$\Pi_i\widetilde R_iq$ does not change the image under $D_i$, and makes
the result orthogonal to its kernel. The minimum norm right inverse is
therefore $R_i=(I-\Pi_i)\widetilde R_i$.

To prove the regularity of the right inverse, we next prove the derivative estimates and smooth dependence. The operator $D_iD_i^*$ is used only on exact forms. To apply the usual
elliptic estimates, we place it inside an elliptic operator on all three
forms, by adding a term that vanishes on closed forms. Write
\begin{align*}
 Y_s=Y\cap H^s(C,\Lambda^3T^*C),\qquad
 L_i=D_iD_i^*+\delta_Cd_C.
\end{align*}
The underlying spaces are fixed, but the adjoints in $L_i$ use the
current metrics. Its symbol has quadratic form
$|\sigma_\xi(D_i)^*q|^2+|\xi\wedge q|^2$. If this is zero, then
$q\in\ker(\xi\wedge)=\operatorname{im}\sigma_\xi(D_i)$ and $q$ is
orthogonal to that same image. Hence $q=0$, proving ellipticity.
The operator is also nonnegative and self-adjoint.

We must check that restriction to exact forms is compatible with this
self-adjoint operator. The current $L^2$ orthogonal decomposition is
$Y\oplus Y^\perp$, where $Y^\perp=\ker\delta_C$, and Hodge theory
makes its projections bounded on $H^2$. If $q$ is exact, then $d_Cq=0$
and $L_iq=D_iD_i^*q$ is exact. If $q\in Y^\perp$, then
$D_i^*q=J_i^*\delta_Cq=0$, so
$L_iq=\delta_Cd_Cq\in\ker\delta_C$. Thus $L_i$ preserves both
orthogonal summands. Its restriction to $Y$ is consequently
self-adjoint as well. All these pairings use the current metrics and 
near the reference pair their norms are uniformly equivalent to the
fixed reference norms.

On $Y$, we have $L_i=D_iD_i^*$. If $L_iq=0$, pairing with $q$ gives
$D_i^*q=0$. Such a $q$ is orthogonal to the image of $D_i$, which is all
of $Y$, so $q=0$. The self-adjoint elliptic operator on this closed
subspace has compact resolvent and zero kernel. It is therefore
invertible, with $L_i^{-1}:Y_0\to Y_2$. Here and below, $L_i^{-1}$ means
the inverse of this restriction, not an inverse on all three forms.
The elliptic estimate at every integer order $s\geq0$ is
\begin{align*}
 \|q\|_{H^{s+2}}
 \leq C_{i,s}\bigl(\|L_iq\|_{H^s}+\|q\|_{L^2}\bigr).
\end{align*}
The low-order inverse bound is uniform in a sufficiently small initial
neighborhood. Indeed, the already constructed right inverse satisfies
$D_i\widetilde R_iq=q$. For $q\in Y_1$, integration by parts gives
$\|q\|_2^2=\langle\widetilde R_iq,D_i^*q\rangle$.
The uniform bound for $\widetilde R_i$ therefore gives
$\|q\|_2\leq C\|D_i^*q\|_2$. This bounds the inverse of $L_i$ in
$L^2$. Applying the preceding elliptic estimate to $L_i^{-1}b$ now gives
\begin{align*}
 \|L_i^{-1}b\|_{H^{s+2}}\leq C_{i,s}\|b\|_{H^s},
 \qquad
 R_i=D_i^*L_i^{-1}:Y_s\longrightarrow H^{s+1}(C,\nu_i).
\end{align*}
Thus $R_i$ gains one derivative, and its formula is exactly
\eqref{eq:subRnear}. These estimates hold on the same low-order
neighborhood at every order. Their constants may depend on higher
norms of the coefficients.
%we do not assert uniform high-order bounds without corresponding coefficient bounds.

For a smooth finite-dimensional family of pairs,
$L_i:Y_{s+2}\to Y_s$ depends smoothly on the parameters. Inversion of
bounded isomorphisms is smooth. Differentiating $L_iL_i^{-1}=I$ gives
$\partial L_i^{-1}=-L_i^{-1}(\partial L_i)L_i^{-1}$.
The coefficients of $S_i^{-1}D_i^*$ also depend smoothly on the pair.
Hence $S_i^{-1}R_i:Y_s\to H^{s+1}(C,\nu_0)$ depends smoothly on the
parameters. This is the asserted smooth dependence of $R_i$ after
identifying the moving normal bundles.
\end{proof}

The correction $W$ in \eqref{eq:introW} is a normal field along the
submanifold. To form the ambient Lie derivative, we need a vector field
on $M$ that restricts to $W$. We choose this extension once as part of
the equation.

Fix a radius $r>0$ smaller than the initial normal tubular radius, and a
smooth cutoff $\rho$ equal to one near zero and zero near one. Keep
$r$ and $\rho$ fixed for nearby pairs. By compactness and the inverse
function theorem, after shrinking the neighborhood, every normal
exponential map is a diffeomorphism on the normal disc bundle of radius
$r$, and its inverse depends smoothly on smooth parameter families.
For $x=\exp_p^{g_\varphi}(\nu)$, with $p\in i(C)$ and $|\nu|<r$,
transport $W(p)$ along the radial normal geodesic from $p$ to $x$.
Writing this transport as $\mathsf{PT}_{p\to x}$, define
\begin{align}\label{eq:extension}
 (E_{\varphi,i}W)(x)=\rho\left(\frac{|\nu|^2}{r^2}\right)\,\mathsf{PT}_{p\to x}W(p),
\end{align}
and extend by zero outside the tube. Since the cutoff vanishes near the
boundary of the tube, the resulting ambient field is smooth.

We record how this choice behaves under changes of parametrization and
ambient coordinates. Replacing $i,W$ by $i\circ a,W\circ a$, where $a$
is a diffeomorphism of $C$, leaves the image and its normal field
unchanged, hence leaves the extension unchanged. If $F$ is an ambient
diffeomorphism, then $F:(M,F^*g)\to(M,g)$ is an isometry. It preserves
normal exponential maps, radial parallel transport, and $|\nu|$.
The extension formed from $F^*\varphi$, $F^{-1}\circ i$, and
$dF^{-1}W$ is therefore $F^*(E_{\varphi,i}W)$. This is the equivariance
needed in Section~\ref{sec:first-existence}. We retain $r,\rho$ as our fixed
data and write $X_{\varphi,i}=E_{\varphi,i}W_{\varphi,i}$.

\begin{proposition}\label{prop:exactpullback}
Every smooth solution of \eqref{eq:correctedMCF} remaining in the
neighborhood of Lemma~\ref{lem:subopen}, where the fixed tubular
extension is defined, satisfies
\begin{align}\label{eq:constantpullback}
 \partial_t(i_t^*\varphi_t)=0.
\end{align}
The cohomology class $[\varphi_t]$ is constant.
At a coassociative pair, $W$ is given by \eqref{eq:introexplicitW}
and is the unique normal field along $i(C)$ of minimum $L^2(C,g_C)$ norm enforcing
\eqref{eq:constantpullback} with submanifold velocity $H$.
\end{proposition}
\begin{proof}
We use the pullback derivative, now including the ambient Lie derivative.
Since $d\varphi=0$, Cartan's formula gives
$\Lie_X\varphi=d(X\ip\varphi)$. Since $X$ restricts to $W$ on the
submanifold, its pullback contribution is $D_iW$. Thus
\begin{align*}
 \partial_t(i^*\varphi)
 =i^*(d\tau+\Lie_X\varphi)+d_Ci^*(H\ip\varphi)=d_C\theta+D_iW+D_iH=0.
\end{align*}
The last equality is the defining equation
$D_iW=-d_C\theta-D_iH$. Notice that no coassociativity assumption was used in this cancellation, because the inverse was defined on an open neighborhood.
Integrating in time proves \eqref{eq:constantpullback}. Also
$\partial_t\varphi=d(\tau+X\ip\varphi)$ is exact, so
$[\varphi_t]=[\varphi_0]$.

In the same way as in Lemma~\ref{lem:subinverse} we get that $W$ has the smallest $L^2$-norm.

We can also remove the Lie derivative from the ambient equation. This
will be useful for proving existence. Let $F_t$ solve
$\partial_tF_t=-X_t\circ F_t$ with $F_0=\mathrm{id}$. The minus sign is
chosen so that differentiation of the pullback cancels $\Lie_X\varphi$. Precisely, we have
\begin{align*}
 \partial_t(F_t^*\varphi_t)
 =F_t^*(\partial_t\varphi_t-\Lie_{X_t}\varphi_t)
 =\Delta_{F_t^*\varphi_t}(F_t^*\varphi_t).
\end{align*}
Set $\widetilde\varphi_t=F_t^*\varphi_t$ and
$\widetilde i_t=F_t^{-1}\circ i_t$. Then
$\widetilde i_t^*\widetilde\varphi_t=i_t^*\varphi_t$.
Differentiating $F_t^{-1}\circ i_t$ adds $W$ to the mean-curvature
velocity: the inverse diffeomorphism moves with the opposite sign to
$F_t$. Since $R_iD_i=I-\Pi_i$ on normal fields, the same velocity can
be written as $H+W=\Pi_iH-R_id_C\theta$. Pulling back gives
\begin{align}\label{eq:ordinary-lift}
 \partial_t\widetilde i_t
 =dF_t^{-1}(H_t+W_t)
 =\Pi_{\widetilde i_t}\widetilde H_t
   -R_{\widetilde i_t}d_C\widetilde\theta_t,
 \qquad \widetilde\theta_t=\widetilde i_t^*
                 (\delta_{\widetilde\varphi_t}\widetilde\varphi_t).
\end{align}
Here $\Pi_{\widetilde i_t}$ is the $L^2$ projection onto
$\ker D_{\widetilde i_t}$, and $\widetilde H_t$ uses the pulled back
metric. The normal operators commute with pullback because $F_t$ is an
isometry between the two metrics. It preserves their defining equations
and the $L^2$ orthogonality used to select the inverse.

Thus, when the ambient $G_2$ structure follows ordinary Laplacian flow, the normal velocity has two parts. The equation $D_iV=-d_C\theta$ forces the
component $-R_id_C\theta$ orthogonal to the kernel. The component
$\Pi_iH$ lies in the kernel and does not change that equation. Keeping
this component selects the admissible velocity closest to $H$ in
$L^2$. At a coassociative pair, $J_i$ identifies $\Pi_i$ with harmonic
projection on anti-self-dual two forms. In particular, $\Pi_iH=0$ if
$b_2^-(C)=0$.
\end{proof}
\begin{remark}
At a coassociative pair with torsion-free $\varphi$, both $\tau$ and $H$
vanish, so $W=0$.	
\end{remark}

We can now give a proof of part of our first main result.

\begin{proof}[Proof of Theorem~\ref{thm:corrected}]
We have from Proposition~\ref{prop:exactpullback} that
$i_t^*\varphi_t=i_0^*\varphi_0$ and $[\varphi_t]=[\varphi_0]$ for every
solution in the stated neighborhood. Since $i_0^*\varphi_0=0$ hence the
embeddings remain coassociative. The same proposition gives the formula
\eqref{eq:introexplicitW} and its minimum norm property. These are the
preservation assertions of the theorem. The short-time existence and
uniqueness assertion is proved separately in
Theorem~\ref{thm:first-existence}.
\end{proof}
The proof of the following corollary is similar to Proposition~\ref{prop:exactpullback} and the  existence and uniqueness part is given at the end of
Section~\ref{sec:first-existence}.

\begin{corollary}\label{cor:fixedC}
With the same initial hypotheses and fixed embedding $i=i_0$, let
$W^0=-R_i d_C(i^*\tau)$ and $X^0=E_{\varphi,i}W^0$. The equation
\begin{align}\label{eq:fixedCflow}
 \partial_t\varphi=\Delta_\varphi\varphi+\Lie_{X^0}\varphi
\end{align}
has a unique smooth closed solution from $\varphi_0$ for a short time,
and every smooth closed solution from $\varphi_0$ remaining in the same
neighborhood preserves $i^*\varphi=0$.
\end{corollary}

\subsection{An explicit compact example}
Before proving the short-time existence and uniqueness assertions in Theorem~\ref{thm:corrected} we give an explicit example which shows the need to modify the coupled system of the Laplacian flow of closed $G_2$-structures and the mean curvature of submanifolds. Precisely, we show that the error in
Proposition~\ref{prop:unmodified-pullback} can be nonzero for compact
initial data, even arbitrarily close to a torsion-free structure.
The example also lets us compute the correcting normal field directly.
%We compute at the initial time; no explicit formula for the full solution is needed to prove failure of preservation.

\begin{proposition}
\label{prop:torus-example}
There are closed $G_2$ structures on $T^7$, arbitrarily close and
cohomologous to a flat torsion-free $G_2$ structure, and compact coassociative
four-tori for which ordinary Laplacian flow coupled to mean curvature
flow fails to preserve the coassociative condition. 
%The normal correction in the first system is nonzero for these data.
\end{proposition}

\begin{proof}
Use coordinates $(y^1,y^2,y^3,x^1,x^2,x^3,x^4)$ of period $2\pi$ and
the model from the earlier mean-curvature calculation. Thus, we have
\begin{align*}
 \varphi_{\mathrm{flat}}
 =dy^{123}+dy^1\wedge(dx^{12}+dx^{34})
  +dy^2\wedge(dx^{13}-dx^{24})
  -dy^3\wedge(dx^{14}+dx^{23}).
\end{align*}
Let
\begin{align}
\label{eq:torus-phi}
 \varphi_\varepsilon
 =\varphi_{\mathrm{flat}}
   +\varepsilon\sin y^1\sin x^1\,dy^1\wedge dx^{12},
 \qquad 0<|\varepsilon|<1.
\end{align}
The added term is
$d(-\varepsilon\cos y^1\sin x^1\,dx^{12})$. It is therefore exact,
so the new form is closed and has the same cohomology class as the flat
form. To check that it is indeed a $G_2$ structure and identify the metric, we put
$a=1+\varepsilon\sin y^1\sin x^1>0$.
Multiply the coordinate coframe by $a^{1/3}$ in the $y^1,x^1,x^2$
directions and by $a^{-1/6}$ in the other four directions. Substituting
this coframe into the flat model gives \eqref{eq:torus-phi}. Thus this
is an orthonormal coframe for its induced metric, and
\begin{align}
\label{eq:torus-metric}
 g_\varepsilon
 =a^{\frac 23}\big((dy^1)^2+(dx^1)^2+(dx^2)^2\big)
  +a^{-1/3}\big((dy^2)^2+(dy^3)^2+(dx^3)^2+(dx^4)^2\big).
\end{align}
Every summand of $\varphi_\varepsilon$ contains a $dy$ factor, so its
pullback to $C=\{y^1=y^2=y^3=0\}$ vanishes. Hence $C$ is coassociative.
Since $a=1$ on $C$, the induced metric is flat, and its coassociative
orientation is $-dx^{1234}$. We remark that flatness of the induced metric does not force $H$ to vanish as $H$ also depends on normal derivatives of the
ambient metric, which need not vanish on $C$.

To calculate $\theta=i^*\tau$, we use $\tau=-*d\psi$. Using the facts that $a=1$ along $C$, a short computation gives
\begin{align}
\label{eq:torus-theta}
 \theta=i^*\tau
 =\frac\varepsilon3\sin x^1(2dx^{12}-dx^{34}).
\end{align}
We can also compute $H$ from the metric. The trace of $\partial_{y^1}g_\varepsilon$ over
$TC$ is
$(\frac23+\frac23-\frac13-\frac13)\varepsilon\sin x^1
=\frac23\varepsilon\sin x^1$. Thus
\begin{align*}
 H=-\tfrac12\tr_{TC}(\partial_{y^1}g_\varepsilon)
       \partial_{y^1}
   =-\frac\varepsilon3\sin x^1\,\partial_{y^1}.
\end{align*}
This agrees with $JH=-2\pi_-\theta$.

Since $*_Cdx^{12}=-dx^{34}$ and $*_Cdx^{34}=-dx^{12}$,
\begin{align}
\label{eq:torus-obstruction}
 d_C(*_C\theta)=-\frac{2\varepsilon}{3}\cos x^1\,dx^{134}\neq0.
\end{align}
Proposition~\ref{prop:unmodified-pullback} now gives a nonzero initial
derivative of $i_t^*\varphi_t$, although its initial value is zero.
The ordinary ambient Laplacian flow exists for a short time. Mean
curvature flow in its smooth time-dependent metric also has the usual
short-time parabolic existence. Thus these initial data give actual
solutions of the unmodified coupling for which coassociativity is not
preserved.

Finally,
\begin{align}
\label{eq:torus-W}
 W=\frac{2\varepsilon}{3}\sin x^1\,\partial_{y^1}
\end{align}
satisfies $D_iW=-d_C(*_C\theta)$. Indeed,
$J_i\partial_{y^1}=dx^{12}+dx^{34}$, and differentiating the sine
coefficient gives
$D_iW=\frac{2\varepsilon}{3}\cos x^1\,dx^{134}$, the negative of
\eqref{eq:torus-obstruction}. Since $C$ is flat, its harmonic
anti-self-dual forms are constant, so the kernel of $D_i$ consists of
constant normal fields. The sine coefficient in \eqref{eq:torus-W}
has zero mean. Hence $W$ is orthogonal to this kernel and is exactly
the minimum norm correction required to make the coupling work.
\end{proof}

\subsection{Existence from a coassociative}\label{sec:first-existence}
We now construct solutions to the coupled systems. Even though we have proved when the coassociatives are preserved this alone does not
prove existence, because the correction is part of the unknown coupled
equation. We instead solve the ordinary ambient flow of the $G_2$ structures first and choose a coassociative embedding above each ambient $G_2$ structure. We then choose the motion within this family to have the required normal velocity for preservation.

The precise family needed for these steps is provided by McLean's deformation theory \cite{McL} and by \cite[Thm. 2.8]{KL}. We use bars below to refer
to the prescribed ambient family. Along a graph $j(t,z)$, the operators
$\bar D$, $\bar R$ and $\bar\Pi$ are formed using
$(\bar\varphi_t,j(t,z))$, with $\bar\Pi$ projecting onto $\ker\bar D$,
and $D_0$ denotes the operator at the initial pair
$(\bar\varphi_0,i_0)$.

%\begin{samepage}
\begin{lemma}\label{lem:selection}
Let $\bar\varphi_t$ be a smooth family of closed $G_2$ structures
with constant cohomology class, and let $i_0$ be coassociative for
$\bar\varphi_0$. Suppose $\partial_t\bar\varphi_t=d\bar\tau_t$ with
$\bar\tau_t$ smooth. Let $B(t,j)\in\Gamma(\nu_j)$ assign a normal field
to each pair $(t,j)$ in the local finite-dimensional family, smoothly
in the family parameters and with $B(t,j\circ a)=B(t,j)\circ a$ for
diffeomorphisms $a$ of $C$. Then there is a unique short-time path of coassociative
embeddings starting at $i_0$ with normal velocity
\begin{align}
 \bar V=\bar\Pi B-\bar R\,d_C(\bar i^*\bar\tau).
\end{align}
%The parametrization is fixed by requiring $\partial_t\bar i=\bar V$. In particular, $B=\bar H$ and $B=0$ are allowed.
\end{lemma}
%\end{samepage}
\begin{proof}
We first describe the nearby coassociatives, without prescribing their
velocity. Use the initial metric to write nearby embedded images as
normal graphs $j_v(p)=\exp_{i_0(p)}v(p)$, where
$v\in C^{2,\alpha}(C,\nu_0)$ is small. Their coassociative equation is
\begin{align*}
 \mathcal F(t,v)=j_v^*\bar\varphi_t
 \quad\hbox{in}\quad
 Y^{1,\alpha}:=d_CC^{2,\alpha}(C,\Lambda^2T^*C)
                  \subset C^{1,\alpha}(C,\Lambda^3T^*C).
\end{align*}
We equip $Y^{1,\alpha}$ with the $C^{1,\alpha}$ norm. Hodge
decomposition makes this space of exact three forms a closed Banach
subspace. The map does take values there because the nearby graphs are
isotopic, the ambient cohomology class is constant, and therefore
$[j_v^*\bar\varphi_t]=[i_0^*\bar\varphi_0]=0$.
This is the point where the cohomology hypothesis enters the local
construction.

At $(t,v)=(0,0)$, the derivative in $v$ is $D_0$. Its kernel is
$J_0^{-1}\mathcal H^2_-(C)$, and the inverse formula
$2J_0^{-1}\pi_-\delta_CG_C$ maps exact $C^{1,\alpha}$ three forms to
$C^{2,\alpha}$ normal fields. Hence the restriction of $D_0$ to the
fixed $L^2$ complement of its kernel is a bounded isomorphism onto
$Y^{1,\alpha}$.

Split the normal field into its kernel component $z$ and its orthogonal
component $w$. The implicit function theorem solves for $w$ in terms
of $(t,z)$. This gives all nearby coassociative graphs $j(t,z)$, with
$z$ in an open set of the finite-dimensional space $\ker D_0$.
%We do not need this kernel to vanish and it becomes the parameter space in which we will solve an ordinary differential equation.

We also need a smooth family on one common time interval, rather than
a separate interval for each derivative order. This is essentially the regularity
argument in \cite[Theorems 2.7 and 2.8]{KL} and we give some details. The graph equation is first order elliptic in the normal variables near a coassociative,
because its linearization has injective symbol. Differentiating the
equation and applying elliptic estimates therefore improves the spatial
regularity of a $C^{2,\alpha}$ graph to every order when the ambient data
are smooth. At each such graph, the linearization on the fixed kernel
complement remains invertible by the original implicit function theorem.
Elliptic regularity gives its inverse in the higher spaces as well.
We may therefore apply the implicit function theorem at each higher
order. The graphs are consequently smooth in space and in
$(t,z)$ on the same parameter neighborhood. 
%Each resulting branch also solves the original low-order equation, whose uniqueness identifies it with the branch already constructed. 

\medskip

We now choose the path in this family. Let
$V_0=(\partial_tj)^\perp$ be the normal velocity obtained by varying
$t$ while holding $z$ fixed, and let
$Z_a=(\partial_{z^a}j)^\perp$ be the normal fields obtained by varying
the parameters. Differentiate $j(t,z)^*\bar\varphi_t=0$ in these two
ways to obtain
\begin{align*}
 \bar D V_0=-d_C(j^*\bar\tau),\qquad \bar D Z_a=0.
\end{align*}
Tangential components do not contribute as their pullback variation is
the Lie derivative of the zero form $j^*\bar\varphi_t$.
The graph parametrization identifies variations in $z$ with all normal
infinitesimal deformations. Hence the $Z_a$ form a basis of
$\ker\bar D$. Write
\begin{align}\label{eq:parametric-basis}
 \mathsf Z_{t,z}v=\sum_a v^aZ_a,
 \qquad \mathsf Z_{t,z}:\mathbb R^{b_2^-(C)}\longrightarrow\ker\bar D.
\end{align}
This map is an isomorphism. The desired velocity $\bar V$ and the
fixed-parameter velocity $V_0$ satisfy the same inhomogeneous equation,
$\bar D\bar V=\bar DV_0=-d_C(j^*\bar\tau)$. Their difference therefore
lies in $\ker\bar D$, and there are unique coefficients $b^a(t,z)$ with
\begin{align*}
 \bar V-V_0=\sum_a b^a(t,z)Z_a.
\end{align*}
Pair this equation with each $Z_b$ in $L^2$. The resulting matrix is
the Gram matrix of a basis, so it is positive definite. Its inverse
and the fields $Z_a,V_0,\bar V$ depend smoothly on $(t,z)$, which makes
the coefficients $b^a$ smooth. We can therefore solve the ODE
\begin{align}\label{eq:selectionode}
 \dot z^a=b^a(t,z),\qquad z(0)=0,
\end{align}
by the finite-dimensional existence and uniqueness theorem. The
resulting graph path has the required normal velocity, but its full
velocity may also have a tangential component. Write it as
$\bar V+dj(Z)$ and solve $\partial_tf_t=-Z_t\circ f_t$ with
$f_0=\mathrm{id}_C$. The derivative of $j_t\circ f_t$ is
$\bar V\circ f_t$: the two terms containing $dj(Z)$ cancel.
Equivariance of $B$ and of the minimum norm inverse ensures that this
is the prescribed normal field for the reparametrized embedding.

Finally, any other nearby coassociative path with the same initial
embedding and prescribed normal velocity has the same equation
\eqref{eq:selectionode} for its graph parameter. Uniqueness of that
ordinary differential equation gives the same moving image. The
reparametrization equation then fixes the labeling of its points
uniquely. This proves uniqueness of the parametrized path.
\end{proof}

In the following proposition, we explain the geometric meaning of the selected velocity.

\begin{proposition}
\label{prop-coassoc-volume-gradient}
Fix a closed $G_2$ structure $\varphi$ and a compact
coassociative embedding $i$. Let $\mathcal M_\varphi$ be its local
coassociative deformation space, with embeddings identified when
they differ only by reparametrization. Give its tangent space
$T_i\mathcal M_\varphi=\ker D_i$ the $L^2$ metric on normal
vector fields. If $\mathcal V_\varphi$ is the volume of the
embedded submanifold in the metric $g_\varphi$, then
\begin{align*}
 \operatorname{grad}_{\mathcal M_\varphi}\mathcal V_\varphi
 =-\Pi_iH.
\end{align*}

Suppose now that the ambient form varies with
$\partial_t\varphi_t=d\tau_t$, and set $\theta=i^*\tau_t$ at a
coassociative pair. The normal velocity
\begin{align*}
 V=\Pi_iH-R_i d_C\theta
   =-R_i d_C\theta
      -\operatorname{grad}_{\mathcal M_{\varphi_t}}
                         \mathcal V_{\varphi_t}
\end{align*}
is the unique velocity closest to $H$ in $L^2$ among the normal
velocities satisfying the differentiated coassociative equation.
Thus it is the sum of the compatible velocity of minimum norm and
the negative volume gradient within the current deformation space.
\end{proposition}

\begin{proof}
First keep the ambient structure fixed. Let $i_s$ be a path of
coassociatives through $i$ whose normal velocity at $s=0$ is $Z$.
Differentiating $i_s^*\varphi=0$ gives $D_iZ=0$. Conversely,
the deformation chart identifies every element of $\ker D_i$
with a tangent vector to $\mathcal M_\varphi$. The first variation
of volume gives
\begin{align*}
 d\mathcal V_\varphi|_i(Z)
 &=-\int_C g_\varphi(H,Z)\,d\mu_i=-\int_C g_\varphi(\Pi_iH,Z)\,d\mu_i,
\end{align*}
where the second equality follows because $Z\in\ker D_i$ and $\Pi_i$
is the orthogonal projection onto that space. This proves the
gradient formula.

If the ambient $G_2$ structure is also moving then a compatible normal velocity
$Y$ must satisfy
\begin{align*}
 0=i^*(\partial_t\varphi_t)+D_iY=d_C\theta+D_iY.
\end{align*}
Since $D_iR_i=I$ on exact three forms, every solution has the form
\begin{align*}
 Y=-R_i d_C\theta+Z,\qquad Z\in\ker D_i.
\end{align*}
The field $R_i d_C\theta$ is orthogonal to $\ker D_i$. Therefore
the orthogonal decomposition of $H$ gives
\begin{align*}
 \|Y-H\|_{L^2}^2
 =\|Z-\Pi_iH\|_{L^2}^2
   +\|R_i d_C\theta+(I-\Pi_i)H\|_{L^2}^2.
\end{align*}
The second term does not depend on $Z$. The first has its unique
minimum at $Z=\Pi_iH$, which proves the assertion. The same
decomposition gives
\begin{align*}
 \|Y\|_{L^2}^2=\|R_i d_C\theta\|_{L^2}^2+\|Z\|_{L^2}^2,
\end{align*}
so the compatible velocity of minimum norm is $-R_i d_C\theta$.
\end{proof}

We now complete the proof of Theorem~\ref{thm:corrected}.

\begin{theorem}\label{thm:first-existence}
Let $\varphi_0$ be a smooth closed $G_2$ structure on compact $M$, and
let $i_0:C\hookrightarrow M$ be a compact coassociative. Fix the tubular
extension data in \eqref{eq:introW}. Then the coupled system \eqref{eq:correctedMCF} has a unique smooth
short-time solution with this initial pair and initial parametrization,
remaining in the neighborhood where the correction is defined. It
preserves $i_t^*\varphi_t=0$ and $[\varphi_t]=[\varphi_0]$.
\end{theorem}
\begin{proof}
We first use the short-time existence theorem for ordinary Laplacian
flow \cite{BX,LW} to obtain $\bar\varphi_t$ from $\varphi_0$.
Set $\bar\tau=\delta_{\bar\varphi}\bar\varphi$.
Its velocity is $d\bar\tau$, so Lemma~\ref{lem:selection} applies with
$B=\bar H$. Since $\bar R\bar D=I-\bar\Pi$, the selected normal
velocity, as described before, is
\begin{align*}
 \bar V=\bar\Pi\bar H-\bar R d_C\bar\theta
       =\bar H+\bar W.
\end{align*}
We next transfer the additional submanifold velocity $\bar W$ into the
ambient equation. Extend it to
$\bar X=E_{\bar\varphi,\bar i}\bar W$ and solve
$\partial_tF_t=\bar X_t\circ F_t$, $F_0=\mathrm{id}_M$.
%This time the sign is positive, because we are reconstructing the corrected system from the ordinary one.
Define
\begin{align}\label{eq:reconstruction}
 \varphi_t=F_t^*\bar\varphi_t,\qquad
 i_t=F_t^{-1}\circ\bar i_t.
\end{align}
The normal operators and our fixed extension rule are equivariant under
ambient diffeomorphisms, as shown in Section~\ref{sec:ambientcorrection}.
Thus $X_{\varphi,i}=F_t^*\bar X$. Differentiating both expressions in
\eqref{eq:reconstruction} gives
\begin{align*}
 \partial_t\varphi=\Delta_\varphi\varphi+\Lie_X\varphi,
 \qquad
 \partial_ti=dF_t^{-1}(\bar H+\bar W-\bar X|_{\bar i(C)})=H.
\end{align*}
For the second identity we used
$\bar X|_{\bar i(C)}=\bar W$, which cancels the extra normal velocity.
After shortening the interval, the pair remains in the neighborhood
where both the inverse and the fixed tubular extension are defined.
This gives the required solution. Closedness of the $G_2$ structures and the coassociative equation follow directly from pullback, and the cohomology statement also follows from the exact ambient velocity.

To prove uniqueness, take another solution with the same initial data.
The preservation identity already proved shows that it remains
coassociative. Solve $\partial_tq_t=-X_t\circ q_t$ and pull back by
$q_t$. The ambient $G_2$ structure $q_t^*\varphi_t$ then solves ordinary Laplacian flow, while $q_t^{-1}\circ i_t$ has normal velocity $\bar H+\bar W$.
Uniqueness of ordinary Laplacian flow fixes the ambient path, and
Lemma~\ref{lem:selection} fixes the lifted embedding. Thus the barred
pair agrees with the one constructed above. Its reconstruction map
$q_t^{-1}$ satisfies $\partial_tF_t=\bar X_t\circ F_t$ with the same
initial value. Uniqueness of this last ordinary differential equation
recovers the same corrected pair.
\end{proof}

\begin{proof}[Existence and uniqueness in Corollary~\ref{cor:fixedC}]
Apply Lemma~\ref{lem:selection} with $B=0$. 
\end{proof}

\subsection{Choice of extension and continuation}
The minimum norm condition determines the normal field $W$ uniquely.
Its extension away from the submanifold is an additional choice. We
now show that changing this extension changes the ambient coordinate
description, while leaving the ordinary-flow lift unchanged.

\begin{proposition}
\label{prop:extension-independence}
Choose two smooth diffeomorphism-equivariant tubular extensions, both equal
to $W$ along the submanifold. The resulting solutions with the same
coassociative initial pair are equivalent by time-dependent ambient
diffeomorphisms for as long as both solutions are defined. In the
description with ordinary Laplacian flow they give the same normal-velocity
embedding $j_t$.
\end{proposition}

\begin{proof}
After removing the ambient Lie derivative, both solutions have the same
ordinary ambient flow and the same prescribed normal velocity.  The uniqueness argument in Theorem~\ref{thm:first-existence} therefore gives the same lifted pair $(\widetilde\varphi_t,j_t)$. Let $G_t^1,G_t^2$ be the two maps that
reconstruct the corrected solutions from this pair. Setting
$F_t=(G_t^1)^{-1}\circ G_t^2$ and using the reconstruction formulas gives
\begin{align*}
 F_t^*\varphi_t^1=\varphi_t^2,\qquad
 F_t^{-1}\circ i_t^1=i_t^2,
\end{align*}
and hence we get the required equivalence of the ambient structure. 
\end{proof}

\begin{proposition}
\label{prop:coassoc-continuation}
Suppose an ordinary Laplacian flow $\widetilde\varphi_t$ extends smoothly
across a finite time $T$. A solution $j_t$ of
\eqref{eq:ordinary-lift} extends across $T$ if its images remain
in a compact subset of the space of smooth embedded coassociatives,
measured in a fixed sufficiently high finite regularity. More explicitly,
it is sufficient that, after reparametrization, the embeddings have
uniform $C^{m+2,\alpha}$ bounds for some $m\geq4$, their induced metrics
are uniformly equivalent to a fixed metric on $C$, and their normal
tubular radii have a positive lower bound.
\end{proposition}

\begin{proof}
The finite-order bounds give a subsequence of images converging in
$C^{m+1,\beta}$ for $0<\beta<\alpha$. The metric lower bound prevents
loss of immersion, and the positive lower bound for the tubular radius
prevents different parts of the image from meeting in the limit. The
limit is therefore an embedded image. Passing to the limit in the
pullback equation shows that it is coassociative for
$\widetilde\varphi_T$. In a normal graph chart, elliptic regularity for
the coassociative equation improves this limit to a smooth image, by
the same argument as in Lemma~\ref{lem:selection}.

At every limiting coassociative, the explicit formula in
Lemma~\ref{lem:subinverse} gives a right inverse for $D_i$.
The compact set of possible limits can therefore be covered by finitely
many deformation charts, with uniform inverse bounds after shrinking
them. In each chart, the moving image is described by
\eqref{eq:selectionode}, a smooth finite-dimensional equation.
Its coefficients are bounded on a smaller chart, and the local
existence time can be chosen uniformly over this finite collection.
Starting sufficiently close to $T$ therefore continues the image
beyond $T$.

The tangential reparametrization fields have the corresponding bounds
in these charts, so their ordinary differential equations also continue.
This recovers the parametrized embedding $j_t$. 
\end{proof}

\begin{remark}
The argument in the previous proposition gives a sufficient continuation
criterion. It does not claim that smoothness of the ambient flow alone
prevents degeneration of the submanifold.
\end{remark}

\subsection{Long-time persistence and convergence}
\label{sec:coassoc-stability}

We now look at the long-time existence of the coupled system. For long-time existence, we need to keep the lifted embedding inside a region where the preceding construction remains valid. We do this by bounding how far its deformation parameter can travel. A small velocity at each time would not be enough, since a small constant velocity can accumulate a large displacement. We will require a small time integral of the relevant velocities.

\medskip

\begin{lemma}
\label{lem:uniform-chart}
Fix a compact coassociative $(\bar\varphi,\bar i)$ and an integer $m\geq4$.
There is a neighborhood of $\bar\varphi$ in $C^{m+3,\alpha}$ such that the
construction of Lemma~\ref{lem:selection} is uniform for closed
forms $\varphi$ in that neighborhood with $[\bar i^*\varphi]=0$.
On a smaller deformation ball, the following estimates hold:
\begin{align}
\label{eq:uniform-chart-bounds}
 |\dot z|\leq C\big(\|\dot\varphi\|_{C^{m+3,\alpha}}
                       +\|\tau\|_{C^{m+1,\alpha}}\big), \quad \|V\|_{C^{m,\alpha}}+\|W\|_{C^{m,\alpha}}
 \leq C\|\tau\|_{C^{m+1,\alpha}},
\end{align}
with $V$ and $W$ the normal velocities as before.
\end{lemma}

\begin{proof}
We first obtain the graph coordinates in exactly the same way as in Lemma~\ref{lem:selection}. We can
then estimate how the graphs move along the ambient motion by differentiating
$\mathcal F(\varphi,z,w(\varphi,z))=0$ while holding $z$ fixed.
For an admissible ambient variation $\eta$, the direct variation of
$\mathcal F$ is $j^*\eta$, so
\begin{align*}
 \mathrm D_w\mathcal F\,[\mathrm D_\varphi w[\eta]]
 =-j^*\eta.
\end{align*}
Applying the uniformly bounded inverse of the $w$ derivative, and then
the uniform estimate for the normal graph parametrization, gives
\begin{align}\label{eq:chart-ambient-derivative}
 \|\mathrm D_\varphi j[\eta]\|_{C^{m+2,\alpha}}
 \leq C\|\eta\|_{C^{m+3,\alpha}}.
\end{align}
In particular the velocity $V_0=(\partial_tj)^\perp$, with $z$
held fixed, satisfies this estimate with $\eta=\dot\varphi$.

We next estimate $V$ and $W$. In the smaller graph neighborhood, the
induced metrics and the maps $J_i,J_i^{-1}$ have uniform coefficient
bounds. Lemma~\ref{lem:subopen} gives a uniform low-order inverse bound
for the Laplacian on exact three forms. The Schauder estimate upgrades
this bound to the derivative orders needed here. Combining it with
\eqref{eq:subexplicit}, we obtain
\begin{align*}
 \|R_iq\|_{C^{m,\alpha}}\leq C\|q\|_{C^{m-1,\alpha}}.
\end{align*}
The formula $J_iH=-2\pi_-\theta$ controls $H$ directly by the torsion
restriction. In $W=-R_id_C(*_C\theta)$, the derivative taken by $d_C$
is recovered by the right inverse. Thus
\begin{align*}
 \|H\|_{C^{m,\alpha}}+\|W\|_{C^{m,\alpha}}
 &\leq C\|\theta\|_{C^{m,\alpha}}\\
 &\leq C\|\tau\|_{C^{m+1,\alpha}}.
\end{align*}
The same bound holds for $V=H+W$ and we get the second estimates in \eqref{eq:uniform-chart-bounds}.

To estimate the deformation parameter, recall that
$Z_a=(\partial_{z^a}j)^\perp$ is a basis of $\ker D_i$.
Its Gram matrix $G_{ab}=\langle Z_a,Z_b\rangle_{L^2(C)}$ is positive
definite at the reference pair. After shrinking the graph neighborhood,
its inverse is uniformly bounded. Pairing
$V-V_0=\sum_a\dot z^a Z_a$ with the basis gives the explicit equation
\begin{align}\label{eq:selection-gram}
 \dot z^a=\sum_b(G^{-1})^{ab}\langle V-V_0,Z_b\rangle_{L^2(C)}.
\end{align}
Combining this formula with \eqref{eq:chart-ambient-derivative}
and the bound for $V$ proves the first estimate of
\eqref{eq:uniform-chart-bounds}.
\end{proof}

We can now prove the following.

\begin{theorem}[An integrable-velocity criterion]
\label{thm:integrable-coassoc}
Let $\bar\varphi$ be torsion-free on compact $M$, and let $\bar i$ be a
compact coassociative embedding.
Fix the tubular extension data with radius smaller than the tubular
radius of $\bar i$. There are $m\geq4$ and $\epsilon>0$ with the
following property. Let $\varphi_t$, $t\geq0$, be a smooth ordinary
Laplacian flow such that
\begin{align}
\label{eq:ambient-smallness}
 \sup_{t\geq0}\|\varphi_t-\bar\varphi\|_{C^{m+3,\alpha}}<\epsilon,
 \qquad [\bar i^*\varphi_t]=0,
\end{align}
and
\begin{align}
\label{eq:ambient-integrability}
 \int_0^\infty
 \big(\|\partial_t\varphi_t\|_{C^{m+3,\alpha}}
       +\|\tau_t\|_{C^{m+1,\alpha}}\big)\,dt<\epsilon.
\end{align}
Assume also that the ambient forms are uniformly bounded in every higher
norm and that the corresponding higher-norm integrals are finite.
For every sufficiently close $\varphi_0$-coassociative $i_0$, the lift
\eqref{eq:ordinary-lift} and the reconstructed corrected system
exist for all time. Both converge smoothly to coassociative pairs with
torsion-free ambient structure. If the ambient velocities and torsion
decay exponentially in every norm, the convergence of the pairs and the
reconstruction diffeomorphisms is exponential in every norm.
\end{theorem}

\begin{proof}
Choose the uniform chart from Lemma~\ref{lem:uniform-chart} with a
parameter ball of radius $2r$, and take the initial graph parameter in
the concentric ball of radius $\frac r2$. On the interval on which the lift
stays in this chart we integrate the bound for $\dot z$ in \eqref{eq:uniform-chart-bounds} to obtain
\begin{align*}
 |z(t)-z(0)|
 \leq C\int_0^t
 \big(\|\partial_s\varphi_s\|_{C^{m+3,\alpha}}
       +\|\tau_s\|_{C^{m+1,\alpha}}\big)\,ds.
\end{align*}
We choose $\epsilon>0$ so that the bound on the right is less than $\frac r2$.
The parameter then remains within distance $r$ of the center. If there
were a first exit from this smaller ball, the same estimate up to the
exit time would give a contradiction. The finite-dimensional equation
therefore stays in a fixed interior region and can be continued for
all time. This is where smallness of the total time integral \eqref{eq:ambient-integrability} is used.

By Lemma~\ref{lem:uniform-chart}, the normal velocities $V_t$ are
integrable in every norm. Integrating from one time to a later time
shows that the embeddings are Cauchy in each of these norms. The limit
$j_\infty$ remains embedded, because the low-order smallness keeps the
images in a graph neighborhood with a positive tubular radius.

The same estimates give a small value of
$\int_0^\infty\|\widetilde X_t\|_{C^1}\,dt$ and finite integrals at
every higher order. Hence the reconstruction maps $F_t$ exist globally
and converge. We now check that their limit is actually a diffeomorphism, rather
than a map with a degenerate derivative. We differentiate their equation
\begin{align*}
 \partial_t(dF_t)
 =(\nabla\widetilde X_t)\circ F_t\,dF_t.
\end{align*}
Thus, we can again use Gronwall's inequality to obtain uniform bounds for $dF_t$ and its inverse. The right-hand side is then integrable, so $dF_t$ is Cauchy. The maps
$F_t$ themselves are Cauchy by the integrability of $\widetilde X_t$.
Repeated differentiation gives the same conclusion at higher orders.
Thus $F_t$ converges smoothly to a diffeomorphism $F_\infty$.
Taking the low-order integral sufficiently small also keeps the
reconstructed pair inside the domain of the chosen fixed extension.

\medskip

The time integral of $\partial_t\varphi_t$, together with its higher
order versions, gives a smooth limit $\varphi_\infty$. It is a $G_2$ structure because the $3$-forms stay in the chosen small neighborhood of
$\bar\varphi$ (see \cite[Thm. 5.1]{LW} and \cite[Thm. 5.1]{Dwi-RHF} for more details on this argument). Since the torsion $\tau_t$ depends smoothly on the form and its first derivative, so $\tau_t$ also converges. If its limit were nonzero, its
norm would be bounded below by a positive constant for all sufficiently
large times, contradicting its finite time integral. Hence
$\tau_\infty=0$.  Passing to the limit in $j_t^*\varphi_t=0$ gives a coassociative limit, and \eqref{eq:reconstruction} gives the limiting corrected pair.

We can also estimate the distances to these limits. At any fixed order
$k$, the construction uses only finitely many derivatives of the ambient $G_2$ structure. Integrating the corresponding velocity estimates from $t$ to
infinity therefore gives an integer $r_k$ for which
\begin{align}\label{eq:coassoc-tail}
 \|\varphi_t-\varphi_\infty\|_{C^k}
  +\|j_t-j_\infty\|_{C^k}+\|F_t-F_\infty\|_{C^k}
 \leq C_k\int_t^\infty
 \left(\|\partial_q\varphi_q\|_{C^{r_k}}
              +\|\tau_q\|_{C^{r_k}}\right)\,dq.
\end{align}
%The constants include the uniform composition bounds for the domain reparametrizations and the ambient reconstruction maps proved above.
If the integrand is bounded by $Ce^{-\lambda q}$, its tail is bounded
by $\left(\dfrac{C}{\lambda}\right)e^{-\lambda t}$. Thus the lift and the reconstruction maps converge exponentially. Finally, composing the limits in
\eqref{eq:reconstruction} gives the same conclusion for the corrected
pair.
\end{proof}

We will need the following lemma which prove quantitative estimates for the ambient Laplacian flow.

\begin{lemma}
	\label{lemma-ambient-laplacian}
	Let $\bar\varphi$ be a torsion-free $G_2$ structure on a compact
	manifold $M$. Fix an integer $r\geq0$ and $0<\alpha<1$, and put
	$k=\max\{5,r+3\}$ and $N=k+4$. Then there are constants
	$\delta,C,\sigma>0$ with the following property. Suppose that
	$\varphi_0$ is smooth, closed and positive, that
	$[\varphi_0]=[\bar\varphi]$, and that
	$\|\varphi_0-\bar\varphi\|_{C^{N,\alpha}}<\delta$.
	Then the ordinary Laplacian flow from $\varphi_0$ exists for all time
	and satisfies
	\begin{align}
		\label{eq:LW-input}
		&\sup_{t\geq0}\|\varphi_t-\bar\varphi\|_{C^{r,\alpha}}
		+\int_0^\infty
		\bigl(\|\partial_t\varphi_t\|_{C^{r,\alpha}}
		+\|\tau_t\|_{C^{r,\alpha}}\bigr)\,dt
		\leq C\|\varphi_0-\bar\varphi\|_{C^{N,\alpha}}.
	\end{align}
	In fact,
	\begin{align*}
		\|\partial_t\varphi_t\|_{C^{r,\alpha}}
		+\|\tau_t\|_{C^{r,\alpha}}
		\leq C e^{-\sigma t}
		\|\varphi_0-\bar\varphi\|_{C^{N,\alpha}}.
	\end{align*}
	The forms $\varphi_t$ are uniformly bounded in every higher spatial
	norm. For each integer $j\geq0$, there is a finite constant $C_j$ such
	that
	\begin{align*}
		\|\partial_t\varphi_t\|_{C^{j,\alpha}}
		+\|\tau_t\|_{C^{j,\alpha}}
		\leq C_j e^{-\sigma t/2}.
	\end{align*}
	These higher-order constants may depend on higher norms of the smooth
	initial form. The flow converges
	smoothly and exponentially to a torsion-free form in
	$\mathrm{Diff}_0(M)\cdot\bar\varphi$.
\end{lemma}

%%% I give all the details below
\begin{proof}
We first use the Laplacian--DeTurck flow of \cite{LW19b} and
write its solution as $\widetilde\varphi_t=\bar\varphi+\eta_t$.
The form $\eta_t$ is exact because the cohomology class is fixed.
Fix $k=\max\{5,r+3\}$ and $N=k+4$, as in the statement.
All norms in the next estimates use the fixed reference metric.

Choose the neighbourhood size $\varepsilon$ in
\cite[Theorem~5.1]{LW19b}, and let $C_1$ be the constant in
their equation~(4.45). That theorem requires
\begin{align*}
 \|\eta_0\|_{C^{k+4}}<\delta_{\rm LW},\qquad
 \|\eta_0\|_{W^{k+4,2}}\leq\frac{\varepsilon}{2C_1}.
\end{align*}
Here $\delta_{\rm LW}$ depends on the fixed background, $k$,
and $\varepsilon$. On the compact reference manifold, both
initial norms are bounded by $C\|\eta_0\|_{C^{N,\alpha}}$.
Reducing our $\delta$ therefore ensures both hypotheses and
gives a global DeTurck solution in that fixed neighbourhood.

We obtain from \cite[eq. (4.45)]{LW19b} 
\begin{align*}
 \sup_{t\geq0}\|\eta_t\|_{C^k}
 \leq C_1\|\eta_0\|_{W^{k+4,2}}
 \leq C\|\eta_0\|_{C^{N,\alpha}}.
\end{align*}
The time-cutoff argument in the proof of \cite[Theorem~5.1]{LW19b},
ending in their \cite[eq. (5.6)]{LW19b}, gives, after a fixed time
$t_0=t_0(k)>0$,
\begin{align*}
 \|\bar\nabla^{\ell}\eta_t\|_{L^2}^{,2}
 \leq C e^{-\lambda_1t/2}\|\eta_0\|_{L^2}^{,2},
 \qquad 0\leq\ell\leq k+4,
\end{align*}
where $\lambda_1>0$ is the reference spectral gap used there.
Their derivative condition is $\lfloor\ell/2\rfloor\leq k-1$ and
it includes $\ell\leq k+4$ because $k\geq5$.
The constants at these fixed orders depend on the background,
$k$, and the chosen neighbourhood, and not on higher initial
derivatives. Summing these inequalities, taking square roots,
and using $W^{k+4,2}\hookrightarrow C^k$ in dimension seven gives
\begin{align*}
 \|\eta_t\|_{C^k}
 \leq C e^{-\lambda_1t/4}\|\eta_0\|_{C^{N,\alpha}},
 \qquad t\geq t_0.
\end{align*}
On $0\leq t\leq t_0$, the uniform estimate above gives the same
bound after increasing $C$ by a fixed factor. Thus the estimate
holds for all $t\geq0$ with $\sigma=\frac{\lambda_1}{4}$.
In particular, the DeTurck limit is $\bar\varphi$ itself.

	Let $V_t$ be the DeTurck vector field. It depends smoothly on the
	form and its first derivative, and it vanishes at $\bar\varphi$.
	Consequently,
	\begin{align*}
		\|V_t\|_{C^{k-1}}
		\leq C\|\eta_t\|_{C^k}
		\leq C e^{-\sigma t}\|\eta_0\|_{C^{N,\alpha}}.
	\end{align*}
	Let $F_t$ solve $\partial_tF_t=-V_t\circ F_t$, with $F_0$ the
	identity. 

Computations similar to the ones in the proof of Theorem~\ref{thm:integrable-coassoc} gives
	\begin{align*}
		\sup_{t\geq0}\|F_t-\mathrm{id}\|_{C^{k-1}}
		\leq C\int_0^\infty\|V_t\|_{C^{k-1}}\,dt
		\leq C\|\eta_0\|_{C^{N,\alpha}}.
	\end{align*}
	Here differences of maps are measured in fixed charts near the
	identity. These estimates give uniform bounds for pullback in the
	norms used below.
	
	The ordinary solution is $\varphi_t=F_t^*\widetilde\varphi_t$.
	Naturality of the Laplacian and of the torsion form gives
	\begin{align*}
		\partial_t\varphi_t
		&=F_t^*(\Delta_{\widetilde\varphi_t}\widetilde\varphi_t),
		&\tau_t&=F_t^*\tau(\widetilde\varphi_t).
	\end{align*}
	The first expression inside the pullback is second order in
	$\eta_t$, and the second is first order and they both vanish when
	$\eta_t=0$. On the fixed small neighbourhood,
	\begin{align*}
		\|\Delta_{\widetilde\varphi_t}\widetilde\varphi_t\|_{C^{r,\alpha}}
		+\|\tau(\widetilde\varphi_t)\|_{C^{r,\alpha}}
		\leq C\|\eta_t\|_{C^{r+2,\alpha}}
		\leq C\|\eta_t\|_{C^k}.
	\end{align*}
	The last inequality follows from $k\geq r+3$. Pullback of these
	forms uses $C^{r+1,\alpha}$ bounds for $F_t$, which are supplied by
	the $C^{k-1}$ bounds because $k-1\geq r+2$. We obtain the stated
	exponential estimate for $\partial_t\varphi_t$ and $\tau_t$.
	Furthermore,
	\begin{align*}
		\|\varphi_t-\bar\varphi\|_{C^{r,\alpha}}
		&\leq
		\|F_t^*\eta_t\|_{C^{r,\alpha}}
		+\|F_t^*\bar\varphi-\bar\varphi\|_{C^{r,\alpha}}\\
		&\leq C\bigl(\|\eta_t\|_{C^k}
		+\|F_t-\mathrm{id}\|_{C^{k-1}}\bigr)
		\leq C\|\eta_0\|_{C^{N,\alpha}}.
	\end{align*}
	Integration of the velocity and torsion bounds proves
	\eqref{eq:LW-input}.
	
	It remains to explain the higher-order assertion without assuming
	smallness of every initial derivative. Theorem~1.3 of \cite{LW19b},
	whose smooth-convergence conclusion is proved after Theorem~5.1,
	gives smooth convergence of this ordinary solution from the
	finite-order neighbourhood just used. Hence $\varphi_t$ has uniform
	bounds in every fixed spatial norm. Since its velocity is a smooth
	second-order expression in $\varphi_t$ and its torsion is a smooth
	first-order expression, both also have uniform bounds in every norm.
	These bounds are finite on the initial time interval because the
	initial form is smooth.
	
	We combine these bounds with the exponential decay already proved
	in $C^0$. For a tensor $u$ and any integer $L>j+\alpha$, the
	interpolation inequality in the fixed reference metric gives
	\begin{align*}
		\|u\|_{C^{j,\alpha}}
		\leq C\|u\|_{C^0}^{1-(j+\alpha)/L}
		\|u\|_{C^L}^{(j+\alpha)/L}.
	\end{align*}
	Take $L=2j+2$ and apply this to
	$u=\partial_t\varphi_t$ and $u=\tau_t$.
	The exponent of the decaying $C^0$ norm is greater than $1/2$,
	so both tensors decay at least as fast as $C_j e^{-\sigma t/2}$
	in $C^{j,\alpha}$. Their time integrals are therefore finite at
	every order. Integrating the velocity from $t$ to infinity proves
	exponential convergence of $\varphi_t$ in every norm. The limit is
	torsion-free and belongs to $\mathrm{Diff}_0(M)\cdot\bar\varphi$. The ordinary limit can differ from
	$\bar\varphi$, which is why the first term in
	\eqref{eq:LW-input} is a uniform smallness bound rather than a
	decay estimate towards $\bar\varphi$.
\end{proof}

We can now give a proof of another of our main theoresm.

\begin{proof}[Proof of Theorem~\ref{thm:first-global}]
	Let $m$ and $\epsilon$ be supplied by
	Theorem~\ref{thm:integrable-coassoc}, with the fixed tubular
	extension data in that theorem. Apply
	Lemma~\ref{lemma-ambient-laplacian} with $r=m+3$ so its finite initial
	regularity can be taken to be $N=m+10$. We reduce the initial
	neighbourhood so that the right-hand side of
	\eqref{eq:LW-input} is smaller than $\epsilon$, after including
	the fixed norm-comparison constants. Then
	\begin{align*}
		\sup_{t\geq0}\|\varphi_t-\bar\varphi\|_{C^{m+3,\alpha}}
		&<\epsilon,\\
		\int_0^\infty
		\bigl(\|\partial_t\varphi_t\|_{C^{m+3,\alpha}}
		+\|\tau_t\|_{C^{m+1,\alpha}}\bigr)\,dt
		&<\epsilon.
	\end{align*}
	The second estimate uses the continuous inclusion
	$C^{m+3,\alpha}\subset C^{m+1,\alpha}$ for the torsion term.
	Since the Laplacian flow preserves the cohomology class and
	$\bar i$ is $\bar\varphi$-coassociative, we also have
	\begin{align*}
		[\bar i^*\varphi_t]
		=\bar i^*[\bar\varphi]=0.
	\end{align*}
	These are exactly
	\eqref{eq:ambient-smallness}--\eqref{eq:ambient-integrability}.
	The higher-order bounds and integrals required by
	Theorem~\ref{thm:integrable-coassoc} follow from the same lemma.
	Taking the initial embedding in the neighbourhood specified in that
	theorem therefore gives the global lift and the global corrected
	solution. The exponential bounds in every spatial norm give their
	smooth exponential convergence.
	
	The ordinary ambient limit belongs to
	$\mathrm{Diff}_0(M)\cdot\bar\varphi$ by
	Lemma~\ref{lemma-ambient-laplacian}. The reconstruction
	diffeomorphisms supplied by
	Theorem~\ref{thm:integrable-coassoc} start at the identity and
	converge smoothly to a diffeomorphism. Their limit is also isotopic
	to the identity as it is $C^1$-close to one of these diffeomorphisms
	for sufficiently large time, and sufficiently close
	diffeomorphisms are isotopic. Thus the corrected ambient limit is
	again in $\mathrm{Diff}_0(M)\cdot\bar\varphi$.
	Finally, the short-time uniqueness in
	Theorem~\ref{thm:first-existence} applies at every time and gives
	uniqueness on the whole interval.
\end{proof}

\begin{remark}
\label{rem:no-coassoc-attraction}
Theorem~\ref{thm:first-global} concerns initially coassociative pairs. For the correction used in this paper, the pullback $i_t^*\varphi_t$ is constant even when it is initially nonzero. Thus this equation does not make arbitrary
nearby embeddings coassociative. The convergence assertion above also
allows the neutral directions coming from harmonic anti-self-dual
two forms but it does not assert convergence to the initial embedding.
\end{remark}

\section{The second coupled system: instantons}\label{sec:second}

We now study our second coupled system of flows in $G_2$ geometry in this section. We are interested in the Yang-Mills flow of connections on a principal $G$-bundle over a manifold $M$ with a $G_2$ structure. We will proceed in the same way as in \textsection~\ref{sec:first} and along the way, we will figure out the correct flow of coclosed $G_2$ structures which preserves the condition of a family of connections being a $G_2$ instanton. More on gauge theory on manifolds with $G_2$ structures can be found in \cite{Wal, Wal17}.

\subsection{The correction and its energy identity}
We now choose the connection velocity so that an instanton remains an
instanton when the ambient four form $\psi$ evolves along a geometric flow. We first calculate the error caused by the ambient motion. We then construct an elliptic
correction and prove an energy identity. As was the case with coassociatives in \textsection~\ref{sec:first} the existence of the resulting flow is a separate question which is addressed in \textsection~\ref{sec:second-existence}.

\medskip

Throughout this section, $d\psi=0$. All adjoints, Hodge stars and inner
products use $g_\psi$ which is the metric which is obtained from $g_{\varphi}, \varphi=*\psi$ and the fixed fibre metric on $E=\ad P$. A difference
of connections, and hence a connection velocity, is an $E$ valued one
form. Our curvature convention is $F_A=dA+\frac12[A\wedge A]$.
The bracket combines the wedge product with the Lie bracket on
$\mathfrak g$. Thus $\partial_tF_A=d_A\partial_tA$ and
$d_A^2=[F_A,\cdot]$ where $d_A$ is the exterior covariant derivative. We will also use the Bianchi identity $d_AF_A=0$.

With our convention, a gauge transformation acts by
$k^*A=k^{-1}Ak+k^{-1}dk$, and its infinitesimal action is $d_Aq$ for
$q\in\Omega^0(E)$. The Coulomb condition $d_A^*a=0$ therefore means
that $a$ is $L^2$ orthogonal to all infinitesimal gauge directions:
$\langle a,d_Aq\rangle_{L^2}=\langle d_A^*a,q\rangle_{L^2}=0$.
This is how we remove changes that only alter the gauge representative.

Recall that $A$ is a $G_2$ instanton if its curvature $F_A\in \Omega^2_{14}$ which means that $F_A\wedge \psi=0$. Using \eqref{eq:mainconv} the instanton condition can we expressed as the vanishing of a one form.
Define $L_\psi:\Omega^2(E)\to\Omega^1(E)$ by
$L_\psi\beta=*(\beta\wedge\psi)$. Let $\pi_7$ be orthogonal projection
to the $\Lambda^2_7$ component, acting on each $E$ component. The adjoint
and the two compositions of $L_\psi$ are
\begin{align}\label{eq:Ladj}
 L_\psi^*a=*(a\wedge\psi)=a^\sharp\ip\varphi,\qquad
 L_\psi L_\psi^*=3\mathrm{id},\qquad L_\psi^*L_\psi=3\pi_7.
\end{align}

Write $\chi=L_\psi F_A$ and $\PP=L_\psi d_A$, as in
\eqref{eq:introdefect}. Thus $\chi=0$ is the instanton equation, and
$|\chi|^2=3|\pi_7F_A|^2$. In particular, the defect measures one
component of the curvature, rather than the whole curvature. For fixed
$\psi$, its linearization in a connection direction $a$ is $\PP a$. We also note, as before that
\begin{align}\label{eq:divdefect}
 d_A^*\chi=-*d_A(F_A\wedge\psi)=0.
\end{align}
For bundle valued forms, $\langle\alpha\wedge\beta\rangle$ denotes
the form obtained by pairing the $E$ components: in an orthonormal
frame of $E$, it is $\sum_\lambda\alpha^\lambda\wedge\beta^\lambda$.
Second, $\PP$ is formally self-adjoint. For $a,b\in\Omega^1(E)$,
integration by parts, using $d\psi=0$, gives
\begin{align*}
 \int_M\langle\PP a,b\rangle\,d\mu
 =\int_M\langle d_Aa\wedge b\rangle\wedge\psi
 =\int_M\langle a\wedge d_Ab\rangle\wedge\psi
 =\int_M\langle a,\PP b\rangle\,d\mu.
\end{align*}
This proves self-adjointness. In particular
\begin{align}\label{eq:graddefect}
 \PP\chi=\PP^*\chi=d_A^*L_\psi^*L_\psi F_A
                          =3d_A^*\pi_7F_A.
\end{align}
For fixed $\psi$, the derivative of $\frac12\int_M|\chi|^2d\mu$ in
the direction $a$ is $\langle\chi,\PP a\rangle_{L^2}$.
Self-adjointness moves $\PP$ onto $\chi$, so the gradient is
$\PP\chi$. This explains the term $-\PP\chi$ in the flow below.

We next compare this gradient with the Yang--Mills gradient. From \eqref{eq:mainconv} we have 
\begin{align}\label{eq:starF}
 *F_A=\varphi\wedge F_A+3*\pi_7F_A.
\end{align}
On two forms in dimension seven, $d_A^*=*d_A*$. Apply $*d_A$ to
\eqref{eq:starF}. The Bianchi identity removes the derivative of $F_A$ in the
first term, leaving $*(d\varphi\wedge F_A)=*(F_A\wedge d\varphi)$.
The last term is $3d_A^*\pi_7F_A=\PP\chi$ by
\eqref{eq:graddefect}. We, thus obtain
\begin{align}\label{eq:YMdecomp}
 d_A^*F_A=*(F_A\wedge d\varphi)+\PP\chi.
\end{align}
Therefore the torsion corrected Yang--Mills velocity
\begin{align}\label{eq:torsionYM}
 -d_A^*F_A+*(F_A\wedge d\varphi)=-\PP\chi
\end{align}
is the negative gradient of the instanton defect energy give by
\begin{align*}
 \frac12\int_M|\chi|^2\,d\mu
 =\frac12\left(\int_M|F_A|^2\,d\mu
                  -\int_M\langle F_A\wedge F_A\rangle\wedge\varphi\right).
\end{align*}

\begin{remark}
For a fixed
ambient structure, this energy is one half of the $\Omega$-Yang--Mills functional considered by Chen and
Wentworth~\cite[Section~1.1]{CW}, with $\Omega=-\varphi$ in our
conventions and with their trace pairing expressed using our fixed
invariant fibre inner product. 
\end{remark}

Now we also allow $\psi_t$ to vary, and set $\Omega^4(M) \ni Q=\partial_t\psi$ and
$b=*(F_A\wedge Q)$. When we differentiate $\chi=*(F_A\wedge\psi)$,
the Hodge star, the curvature, and the final factor $\psi$ each contribute
one term. Keeping these contributions separate gives
\begin{align}\label{eq:defectgeneral}
 \partial_t\chi=K_t\chi+\PP\partial_tA+b,
 \qquad K_t=(\partial_t*_{t})*_{t}^{-1}\big|_{\Omega^1}.
\end{align}
The differentiated star acts on six forms. Thus $K_t$ is an algebraic
operator determined by $g_t^{-1}\partial_tg_t$, with no derivatives of
$\chi$. For ordinary Yang--Mills flow, substitute
$\partial_tA=-d_A^*F_A$ and use \eqref{eq:YMdecomp} to get
\begin{align}\label{eq:forcingb}
 \partial_t\chi=K_t\chi-\PP^2\chi
            -\PP\bigl(*(F_A\wedge d\varphi)\bigr)+b.
\end{align}
At an instanton, the terms containing $\chi$ vanish. The remaining
terms come from torsion and the variation of $\psi$. Thus we see that the ordinary
Yang--Mills flow would preserve the instanton equation to first order
only if $*(F_A\wedge Q)=\PP\bigl(*(F_A\wedge d\varphi)\bigr)$.
Adding the torsion term in \eqref{eq:torsionYM} removes the torsion
error, but leaves the ambient forcing $b$.

At an instanton, the compatibility equation obtained by
differentiating the instanton condition also appears in the study
of simultaneous deformations of the ambient structure and the
connection by de la Ossa, Larfors and
Svanes~\cite[Section~5.1]{OLS16}. Their ambient deformations are in the setting of torsion-free $G_2$ structures. Here we allow coclosed backgrounds with torsion
and use the compatibility equation to choose a connection velocity
for the prescribed ambient variation $Q$.

We would therefore like to solve $\PP u=-b$. We cannot invert $\PP$
on all one forms because at an instanton, it vanishes on every gauge
direction $d_Aq$. Instead, we impose $d_A^*u=0$ and introduce a function (or a $0$-form) $f$. The enlarged operator is precisely $\DD$ from \eqref{eq:introD} and
its first component imposes the gauge condition, and its second
component is $d_Af+\PP u$. The first step is to check ellipticity of the operator. We do this below.

\medskip

\begin{lemma}\label{lem:D}
For $d\psi=0$, the operator $\DD$ in \eqref{eq:introD} is formally
self-adjoint and elliptic
on $\Omega^0(E)\oplus\Omega^1(E)$. At an instanton its kernel is
\begin{align}\label{eq:Dkernel}
 \ker\DD=\ker d_A|_{\Omega^0(E)}\oplus\mathcal H^1_{A,\psi},\qquad \text{where} \qquad
 \mathcal H^1_{A,\psi}=\{a:d_A^*a=0,\ \PP a=0\}.
\end{align}
\end{lemma}
\begin{proof}
Self adjointness follows from that of $\PP$ and the adjoint relation
between $d_A$ and $d_A^*$. With the symbol convention of
Lemma~\ref{lem:subopen}, the principal symbol at a real covector $\xi$ is
\begin{align}\label{eq:Dsymbol}
 \sigma_\xi(\DD)(f,a)=\sqrt{-1}\,(-\langle\xi,a\rangle,\ \xi f+\xi\times a).
\end{align}
The cross product on covectors is the metric dual of the cross
product on vectors and acts trivially on the $E$ factor. For ordinary
one forms $a,b$, the adjoint identity \eqref{eq:Ladj} gives
$\langle *(\xi\wedge a\wedge\psi),b\rangle
=\langle\xi\wedge a,b^\sharp\ip\varphi\rangle
=\varphi(\xi^\sharp,a^\sharp,b^\sharp)$, and applying the same identity
componentwise to $E$ valued forms proves the cross product expression
in \eqref{eq:Dsymbol}. Since
$\xi\times(\xi\times a)=-|\xi|^2a+\langle\xi,a\rangle\xi$, the square
of this symbol is $|\xi|^2\mathrm{id}$, which proves ellipticity.

At an instanton the compositions
\begin{align}\label{eq:complex}
 \PP d_Af=*([F_A,f]\wedge\psi)=0,\qquad
 d_A^*\PP a=-*([F_A,a]\wedge\psi)=0
\end{align}
vanish. The latter identity follows from $[F_A,a]\wedge\psi=[F_A\wedge\psi,a]$.
Applying $d_A^*$ to the second equation in $\DD(f,a)=0$ therefore
yields $d_A^*d_Af=0$, and integration gives $d_Af=0$ and then
$\PP a=0$, proving \eqref{eq:Dkernel}.
\end{proof}

Note that the two kernel spaces have different meanings.
The space $\ker d_A|_{\Omega^0(E)}$ consists of infinitesimal gauge
transformations fixing $A$. The space $\mathcal H^1_{A,\psi}$ consists
of infinitesimal instanton deformations orthogonal to gauge directions.
The notion of acyclicity means that both spaces vanish. This is the usual
nondegeneracy condition for isolated instantons which is discussed in detail in the work of Walpuski \cite{Wal}.

For a self-adjoint elliptic operator on a compact manifold, the
cokernel is identified with the kernel. Thus zero kernel gives both
existence and uniqueness for the linear equation. The inverse gains
one derivative, with estimates
\begin{align}\label{eq:ellipticestimate}
 \|\DD^{-1}v\|_{H^{s+1}}\leq C_s\|v\|_{H^s},\qquad s\geq0.
\end{align}
To compare nearby operators, fix the bundle
$\BB=E\oplus(T^*M\otimes E)$ and reference norms, and write
$\DD_0=\DD_{A_0,\psi_0}:H^1(\BB)\to L^2(\BB)$ for an invertible
reference operator, whose coefficients give the estimate
\begin{align*}
 \|\DD_{A,\psi}-\DD_0\|_{H^1\to L^2}
 \leq C\bigl(\|A-A_0\|_{C^0}+\|\psi-\psi_0\|_{C^1}\bigr).
\end{align*}
Although the comparison uses fixed norms, the adjoints still use the
current metric. Its first derivatives enter the zeroth-order
coefficients, which explains the $C^1$ norm of $\psi-\psi_0$.
If the perturbation is small enough that
$\|(\DD_{A,\psi}-\DD_0)\DD_0^{-1}\|_{L^2\to L^2}<1$,
\begin{align*}
 \DD_{A,\psi}^{-1}
 =\DD_0^{-1}\bigl[I+(\DD_{A,\psi}-\DD_0)\DD_0^{-1}\bigr]^{-1}
 :L^2(\BB)\longrightarrow H^1(\BB),
\end{align*}
where the bracketed inverse is a norm convergent Neumann series.
For smooth coefficients, the elliptic estimate
$\|v\|_{H^{s+1}}\leq C_s(\|\DD v\|_{H^s}+\|v\|_{L^2})$
and this base inverse prove \eqref{eq:ellipticestimate} at every integer $s\geq0$
on the same neighborhood. Smooth inversion between the fixed spaces
$H^{s+1}(\BB)$ and $H^s(\BB)$, with the inverse differentiation
formula used in Lemma~\ref{lem:subopen}, gives smooth dependence on
smooth parameter families at each order. The estimate constants may
depend on higher coefficient norms.

For the calculations below, abbreviate $u=\UU_{A,\psi}(Q)$. The expression described in~\eqref{eq:introU} is then
\begin{align}\label{eq:ellipticcorrection}
 \DD_{A,\psi}(f,u)=(0,-b).
\end{align}
Thus
\begin{align}\label{eq:correctionequations}
 d_A^*u=0,\qquad \PP u+d_Af=-b.
\end{align}
These equations define $(f,u)$ on an open neighborhood, including
connections with possibly nonzero defect. Invertibility makes the solution
unique. Every operation commutes with gauge transformations, so the
correction is gauge-equivariant. It is nonlocal because computing it
requires an elliptic solution over $M$. We retain the scalar $f$ in
this open neighborhood as it need not vanish away from instantons but as we have already see, it vanishes on an instanton locus.
%The next proposition explains why it vanishes on the instanton locus.

\begin{proposition}\label{prop:horizontal}
At an acyclic instanton, the auxiliary section $f$ in \eqref{eq:ellipticcorrection}
is zero, and $u$ is the unique solution of
\begin{align}
 d_A^*u=0,\qquad *(d_Au\wedge\psi)=-*(F_A\wedge Q).
\end{align}
Among all velocities $v$ satisfying
$d_Av\wedge\psi+F_A\wedge Q=0$, it has minimal $L^2$ norm.
\end{proposition}
\begin{proof}
We have already proved $f=0$. If $v$ is another compatible velocity, we can solve
$d_A^*d_Aq=d_A^*(v-u)$ because the covariant Laplacian $d_A^*d_A$
on $\Omega^0(M,E)$ is invertible, its kernel being $\ker d_A=\{0\}$.
By the choice of $q$, the form $v-u-d_Aq$ is Coulomb. It is also
killed by $\PP$ which can be seen from the facts that the velocities $v$ and $u$ satisfy the same
compatibility equation, and $\PP d_Aq=0$. Hence it belongs to
$\mathcal H^1_{A,\psi}$ and vanishes by acyclicity. Thus, $v=u+d_Aq$. Since $d_A^*u=0$, these summands are orthogonal, giving
$\|v\|_2^2=\|u\|_2^2+\|d_Aq\|_2^2$ and the minimum norm claim.
\end{proof}

The Coulomb condition therefore has a concrete role here: by removing
the component of the velocity tangent to the gauge orbit, it selects
the minimum norm solution of the cancellation equation, just as
orthogonality does in the coassociative construction in \textsection~\ref{sec:first}.

\medskip
Consider a smooth solution of \eqref{eq:coupledYM} on an interval where the correction is defined. By
\eqref{eq:torsionYM}, its equation is
\begin{align}
 \partial_tA=-\PP\chi+u.
\end{align}
Combining \eqref{eq:defectgeneral} and \eqref{eq:correctionequations}
gives the defect equation
\begin{align}\label{eq:exactdefect}
 \partial_t\chi=-\PP^2\chi-d_Af+K_t\chi.
\end{align}

Notice that the ambient forcing $b$ from \eqref{eq:forcingb} has been replaced by $-d_Af$. Although this term need not vanish pointwise, it disappears from the energy identity when we do
integration by parts and pair $f$ with $d_A^*\chi=0$. This uses
\eqref{eq:divdefect} for arbitrary connections, not just instantons.
We now calculate all contributions from the changing metric.

\begin{theorem}\label{thm:energy}
Let $M$ be compact without boundary, let $\psi_t$ be a smooth family of
closed positive four forms, and let $A_t$ be a smooth solution of
\eqref{eq:coupledYM}. Suppose $\DD_{A_t,\psi_t}$ is invertible throughout
the interval.
Let $h=\partial_tg_t$ and define
$\mathcal E(t)=\frac12\int_M|\chi|_{g_t}^2d\mu_t$. Then
\begin{align}\label{eq:energyidentity}
 \mathcal E'(t)=-\|\PP\chi\|_{L^2(g_t)}^2
       +\frac12\int_M h^{ij}\langle\chi_i,\chi_j\rangle\,d\mu_t
       -\frac14\int_M(\tr_{g_t}h)|\chi|^2\,d\mu_t.
\end{align}
Let $\|h\|_{\mathrm{op},\infty}$ denote the supremum of the pointwise
operator norm of $g^{-1}h$, and let $\|\tr h\|_\infty$ be the
supremum of the absolute metric trace. Then
\begin{align}\label{eq:energybound}
 \mathcal E(t)\leq\mathcal E(0)
                 \exp\left(\int_0^t \|h\|_{\mathrm{op},\infty}+\frac12\|\tr h\|_\infty\,ds\right).
\end{align}
\end{theorem}

\begin{proof}
This is a strightforward computation.
\end{proof}
If the metric is stationary, the defect energy is nonincreasing,
including when the positive four form varies within that fixed metric. This gives the following

\begin{corollary}
\label{cor:defect-decay}
Along any smooth corrected connection solution, suppose
\begin{align*}
 \|\DD_{A_t,\psi_t}v\|_{L^2}\geq\lambda(t)\|v\|_{L^2}
\end{align*}
with $\lambda(t)>0$. Let
$c_h=\|g^{-1}\partial_tg\|_{\mathrm{op},\infty}
+\frac12\|\tr_g\partial_tg\|_\infty$.
Then
\begin{align*}
 \mathcal E(t)\leq \mathcal E(0)
 \exp\left(\int_0^t(c_h(q)-2\lambda(q)^2)\,dq\right).
\end{align*}
In particular, $c_h\leq2\lambda^2-\eta$, $\eta>0$, implies
$\mathcal E(t)\leq \mathcal E(0)e^{-\eta t}$.
\end{corollary}

\begin{proof}
The Bianchi identity gives $d_A^*\chi=0$, and hence
$\DD(0,\chi)=(0,\PP\chi)$. Therefore
$\|\PP\chi\|_2^2\geq\lambda^2\|\chi\|_2^2=2\lambda^2\mathcal E$.
We insert this in the defect-energy identity \eqref{eq:energybound} above and apply Gronwall's inequality.
\end{proof}
% note:This energy estimate alone does not replace the connection and regularity estimates in Theorem~\ref{thm:moving-stability}.

We now give the proof of Theorem~\ref{thm:transport}.

\begin{proof}[Proof of Theorem~\ref{thm:transport}]
Let $A_t$ be a smooth solution of \eqref{eq:coupledYM} on an interval
where the correction is defined. The initial defect vanishes. On each compact time subinterval,
Theorem~\ref{thm:energy} bounds the nonnegative energy by its zero
initial value times a finite factor. Thus $\mathcal E(t)=0$.
The defect is smooth, so its zero $L^2$ norm makes it vanish
everywhere and hence we have $F_{A_t}\wedge\psi_t=0$. By
\eqref{eq:torsionYM} the connection equation then reduces to
$\partial_tA=\UU_{A,\psi}(Q)$, and Proposition~\ref{prop:horizontal}
gives $f=0$ together with the two equations in
\eqref{eq:horizontaltransport}, while existence and uniqueness for a
short time are proved in Theorem~\ref{thm:instanton}.
\end{proof}

\subsection{Existence for the connection equation}\label{sec:second-existence}
As in the case of coassociatives, we now construct the solutions used in the preservation argument. The basic strategy is the same as in \textsection~\ref{sec:first}.
We first prove that the correction has the required mapping properties.
We then use the nearby instanton branch to construct the flow from an
instanton. Finally, a gauge-fixed parabolic equation gives existence
from nearby connections whose defect may be nonzero.

Let $\psi_t$, $0\leq t\leq T_0$, be a prescribed smooth family of
closed positive four forms, and set $Q_t=\partial_t\psi_t$.
Pairs of zero and one forms are sections of the fixed bundle $\BB$
introduced after Lemma~\ref{lem:D}.
Choose a smooth reference connection $A_\star$ for which
$\DD_{A_\star,\psi_0}$ is invertible, hold $A_\star$ fixed, and write
nearby connections as $A=A_\star+a$, with $a\in\Omega^1(M,E)$. The
initial connection $A_{\mathrm{in}}$ in
Theorem~\ref{thm:nearby-existence} is chosen near this reference
connection, and neither connection is required to be an instanton.
Writing $\DD_\star=\DD_{A_\star,\psi_0}$, we use the Sobolev norms
from Section~\ref{sec:basics} with fixed metric $g_0=g_{\psi_0}$ and
the connection induced by $A_\star$. We write $\mathcal L(H^m,H^n)$
for bounded linear maps with the operator norm. We use $H^6$ for connections and $H^5$ for curvature, accounting for
the one derivative in $F_A$. In dimension seven, $H^5$ is an algebra
and $H^6$ controls two classical derivatives of the coefficients.
In particular, multiplication maps $H^6\times H^5$ continuously into
$H^5$. The inverse of $\DD$ will return the derivative used in forming
the curvature.

\begin{lemma}\label{lem:Umap}
On a sufficiently small $H^6$ neighborhood of $A_\star$ and a sufficiently
short ambient interval, the map
\begin{align*}
 H^6(M,T^*M\otimes E)\ni a
 \longmapsto\UU_{A_\star+a,\psi_t}(Q_t)
 \in H^6(M,T^*M\otimes E)
\end{align*}
is smooth, with local Lipschitz constants uniform in $t$.
It is smooth at every higher Sobolev order on the intersection with
that same neighborhood. For smooth $a=A-A_\star$ with
$\|a\|_{H^6}\leq K$ in this neighborhood, and every integer $r\geq6$, we have
\begin{align}\label{eq:Utame}
 \|(f,u)\|_{H^r}\leq C_{r,K}(1+\|a\|_{H^r}).
\end{align}
%The constants include the bounds for the prescribed ambient family and the uniform low order inverse bound.
On bounded subsets at any order $s\geq6$ one also has
\begin{align}\label{eq:u-lipschitz}
 \|\UU_{A_1,\psi}(Q)-\UU_{A_2,\psi}(Q)\|_{H^s}
 \leq C\|Q\|_{C^{s+1}}\|A_1-A_2\|_{H^s}.
\end{align}
For smooth connections with uniform coefficient bounds and a uniform
low-order inverse bound, the correction is locally Lipschitz in $L^2$.

\end{lemma}
\begin{proof}
We first work at the low order. At the reference pair, elliptic
invertibility gives $\DD_\star:H^6(\BB)\to H^5(\BB)$ as an
isomorphism. For fixed $\psi_t$, changing $A$ changes only the
zeroth-order coefficients of $\DD$: the difference of two covariant
derivatives is multiplication by the difference of the connections.
The product estimates therefore give
\begin{align*}
 \|\DD_{A,\psi_t}-\DD_\star\|_{\mathcal L(H^6,H^5)}
 \leq C\bigl(\|A-A_\star\|_{H^6}
                    +\|\psi_t-\psi_0\|_{C^6}\bigr).
\end{align*}
Choose the neighborhood and time interval so that
$\|(\DD_{A,\psi_t}-\DD_\star)\DD_\star^{-1}\|_{\mathcal L(H^5,H^5)}<\frac 12$.
Then
\begin{align*}
 \DD_{A,\psi_t}^{-1}
 =\DD_\star^{-1}
 [I+(\DD_{A,\psi_t}-\DD_\star)\DD_\star^{-1}]^{-1}
 :H^5\longrightarrow H^6
\end{align*}
is bounded uniformly and depends smoothly on $A$, while the curvature
\begin{align*}
 F_{A_\star+a}=F_{A_\star}+d_{A_\star}a+\tfrac12[a\wedge a]
\end{align*}
is a smooth map $H^6\to H^5$. Multiplying the curvature by the
prescribed smooth $Q_t$ preserves $H^5$, and applying the inverse
returns it to $H^6$. This proves the first assertion of the
lemma and the local Lipschitz estimate. It also explains why the
nonlocal correction does not lose a derivative of the connection.

Writing $w=(f,u)$, we observe that the difference
$\DD_{A,\psi_t}-\DD_{A_\star,\psi_t}$ is multiplication by $a$ followed
by fixed algebraic contractions, so for $r>6$ we can combine the
elliptic estimate for the smooth reference operator with the product
inequality
\begin{align*}
 \|(\DD_{A,\psi_t}-\DD_{A_\star,\psi_t})w\|_{H^{r-1}}
 \leq C_r\bigl(\|a\|_{L^\infty}\|w\|_{H^{r-1}}
                  +\|a\|_{H^{r-1}}\|w\|_{L^\infty}\bigr),
\end{align*}
to give, after interpolation,
\begin{align*}
 \|w\|_{H^r}\leq C_{r,K}\left(
 \|\DD_{A,\psi_t}w\|_{H^{r-1}}
 +(1+\|a\|_{H^{r-1}})\|w\|_{H^6}\right).
\end{align*}
The interpolation step uses
$\|w\|_{H^{r-1}}\leq\varepsilon\|w\|_{H^r}
+C_{r,\varepsilon}\|w\|_{H^6}$ for $r>6$.
The bound on $\|a\|_{H^6}$ controls $\|a\|_{L^\infty}$, so we can
choose $\varepsilon$ small enough to move the resulting multiple of
$\|w\|_{H^r}$ to the left. Also $H^6$ controls
$\|w\|_{L^\infty}$. For
our solution the base inverse estimate bounds $\|w\|_{H^6}$ in terms
of $K$, while $*(F_A\wedge Q)$ has $H^{r-1}$ norm at most
$C_{r,K}(1+\|a\|_{H^r})$, giving \eqref{eq:Utame}, with the case
$r=6$ following directly from the base estimate.

For an arbitrary source in $H^{r-1}$, we start with its $H^6$ solution
and apply the reference elliptic estimate successively to
$\DD_{A_\star,\psi_t}w=\DD_{A,\psi_t}w-
(\DD_{A,\psi_t}-\DD_{A_\star,\psi_t})w$.
Beginning with the inverse already constructed at the low order,
this argument successively gives one more derivative whenever the
source and coefficients have that derivative. It yields the inverse
$H^{r-1}\to H^r$ on the same neighborhood, with bounds uniform on
sets where $\|a\|_{H^r}$ is bounded. Smooth inversion between these
fixed Sobolev spaces proves smooth dependence at each order. We remark that we are not imposing a new smallness condition on every high derivative.
For proving \eqref{eq:u-lipschitz}, we differentiate
$\DD_{A,\psi}w=-(0,*(F_A\wedge Q))$ in a connection direction $b$ to get
\begin{align*}
 \DD_{A,\psi}\,\partial_Aw[b]
 =-(\partial_A\DD[b])w-(0,*(d_Ab\wedge Q)),
\end{align*}
and similar self-adjointness and elliptic regularity argument proves \eqref{eq:u-lipschitz}.
\end{proof}

\begin{lemma}\label{lem:IFT}
Let $A_0$ be an acyclic instanton for a closed positive $\psi_0$.
Every smooth closed positive $\psi$ sufficiently close to $\psi_0$ in
$C^{2,\alpha}$ has a unique nearby instanton in the slice
$d_{A_0}^{*_{g_0}}(A-A_0)=0$. It is acyclic and depends smoothly on
smooth finite-dimensional parameter families of $\psi$.
\end{lemma}
\begin{proof}
We use the implicit function theorem argument following
\cite[Proposition~2.17]{Wal17}. The enlarged equation allows us to
work with co-closed structures. We keep the gauge condition fixed
at $(A_0,g_0)$, so that the domain and target are fixed Banach spaces.
Consider
\begin{align*}
 \mathcal F(f,a,\psi)=
 \left(d_{A_0}^{*_{g_0}}a,
 *_\psi(F_{A_0+a}\wedge\psi)+d_{A_0+a}f\right).
\end{align*}
The unknowns lie in $C^{2,\alpha}(\BB)$, the target is
$C^{1,\alpha}(\BB)$, and the ambient parameter ranges over positive
forms in the closed subspace
$\ker d\subset C^{2,\alpha}(\Lambda^4T^*M)$. This map is smooth,
with derivative in $(f,a)$ at $(0,0,\psi_0)$ equal to
$\DD_{A_0,\psi_0}$, which is an isomorphism by acyclicity and elliptic
regularity, so the implicit function theorem supplies
$(f(\psi),a(\psi))$.

We must check that a solution of the enlarged equation is an actual
instanton. Write $A=A_0+a$ and apply $d_A^*$ to its second equation.
By $d_A^*\chi=0$, we get $d_A^*d_Af=0$ and hence $f=0$ and we get an acyclic instanton. 

For smooth $\psi$, the slice and instanton equations have smooth
principal coefficients and terms polynomial in $a$ with no derivatives
of $a$ beyond the linear first order part, so elliptic regularity makes
$a$ smooth. At the solution, the derivative in the unknowns is
\begin{align*}
 B_\psi(q,b)=\left(d_{A_0}^{*_{g_0}}b,d_Aq+\PP_{A,\psi}b\right).
\end{align*}
The implicit function theorem makes this derivative invertible at the
original low regularity, and its symbol remains elliptic nearby.
Elliptic regularity therefore upgrades its inverse to
$C^{r,\alpha}(\BB)\to C^{r+1,\alpha}(\BB)$ for every integer $r\geq1$.

We also need smooth dependence on one common parameter neighborhood.
Fix the neighborhood on which the $C^{2,\alpha}$ branch was
constructed. At each parameter value, the higher-order implicit
function theorem gives a local branch of higher regularity. That
branch solves the same equation in the original neighborhood, so
uniqueness at the low order identifies it with the original branch.
The original branch is therefore smooth at every order throughout
the fixed neighborhood. For a path, one time interval suffices;
we do not take an intersection of shrinking intervals as the
regularity order increases.
\end{proof}

We now prove the short-timme existence and uniqueness for the corrected system of flows.

\begin{theorem}\label{thm:instanton}
Let $\psi_t$ be a prescribed smooth closed positive family and let $A_0$
be an acyclic $\psi_0$ instanton. The coupled Yang-Mills system in equation~\eqref{eq:coupledYM} has a unique
smooth short-time solution with initial value $A_0$. It remains an
instanton and satisfies
\begin{align*}
 d_A^*\partial_tA=0,\qquad
 \PP_{A,\psi}\partial_tA=-*(F_A\wedge Q).
\end{align*}
\end{theorem}
\begin{proof}
We obtain a path of instantons $\widehat{A}_t$ in a fixed reference slice by Lemma~\ref{lem:IFT}. The slice fixes the representatives, but its
condition is different from requiring the velocity to be Coulomb for
the current connection and current metric. We impose the latter
condition by a time-dependent gauge transformation.

First, differentiating the instanton equation gives
\begin{align*}
 \PP_{\widehat A,\psi}\partial_t\widehat A
             =-*(F_{\widehat A}\wedge Q).
\end{align*}
The Hodge star variation term vanishes because the instanton defect is zero.
Using the current metric, we solve
\begin{align*}
 d_{\widehat A}^*d_{\widehat A}q
                  =-d_{\widehat A}^*\partial_t\widehat A.
\end{align*}
Acyclicity makes this Laplacian invertible on $E$ valued zero forms,
so the equation determines a smooth $q_t$, which we use to generate a
family of gauge transformations $k_t$. With gauge action
$k^*A=k^{-1}Ak+k^{-1}dk$, we solve $\dot k_tk_t^{-1}=q_t$,
$k_0=\mathrm{id}$, and set $A_t=k_t^*\widehat A_t$.
Then
\begin{align*}
 \partial_tA_t=k_t^{-1}
   (\partial_t\widehat A_t+d_{\widehat A_t}q_t)k_t.
\end{align*}
This velocity satisfies $d_A^*\partial_tA=0$ and the displayed
linearized instanton equation, since $\PP d_Aq=0$ at an instanton,
so Proposition~\ref{prop:horizontal} identifies it with the unique
solution $\UU_{A,\psi}(Q)$. Together with $\chi=0$ and
\eqref{eq:torsionYM}, this shows that the path solves
\eqref{eq:coupledYM}.

Every other smooth solution with the same initial data remains an
instanton by \eqref{eq:energyidentity} and hence solves
$\dot A=\UU_{A,\psi_t}(Q_t)$, for which Lemma~\ref{lem:Umap} gives
local Lipschitz dependence in $H^6$, so uniqueness follows from the
Banach space ODE theorem.
\end{proof}

\begin{theorem}\label{thm:nearby-existence}
Let $M$ be compact without boundary and let $\psi_t$, $0\leq t\leq T_0$, be a prescribed smooth family of closed positive four forms. Suppose $A_\star$ is smooth and $\DD_{A_\star,\psi_0}$ is invertible. There are $\delta>0$ and $0<T\leq T_0$ such that every smooth connection $A_{\mathrm{in}}$ with $\|A_{\mathrm{in}}-A_\star\|_{H^6}<\delta$ admits a unique smooth solution of \eqref{eq:coupledYM} on $[0,T]$ with $A(0)=A_{\mathrm{in}}$. The operator $\DD_{A_t,\psi_t}$ remains invertible. The initial connection need not be an instanton.
\end{theorem}
\begin{proof}
We first solve an equation with an additional gauge term, and then
remove that term by a time-dependent gauge transformation. This is
the construction of \cite[Section~2.1]{Hong}, with
Lemma~\ref{lem:Umap} controlling the nonlocal correction. We give the existence and uniqueness argument for our equation below. The additional point to check is that the elliptic
correction has the required mapping and difference estimates and
these are supplied by Lemma~\ref{lem:Umap}.

Write $\mathcal V_t(A)=-d_A^{*_{t}}F_A+*_t(F_A\wedge d\varphi_t)+\UU_{A,\psi_t}(Q_t)$
for the right-hand side of \eqref{eq:coupledYM}. Since each operation
defining $\UU$ commutes with a gauge transformation, the convention
$k^*A=k^{-1}Ak+k^{-1}dk$ gives
\begin{align*}
 \UU_{k^*A,\psi}(Q)=k^{-1}\UU_{A,\psi}(Q)k,
 \qquad \mathcal V_t(k^*A)=k^{-1}\mathcal V_t(A)k.
\end{align*}
Let $B=A_\star+a$ be the connection in the gauge fixed equation and set
\begin{align*}
 \zeta_t(B)=d_{A_\star}^{*_{g_t}}(B-A_\star).
\end{align*}
We solve
\begin{align}\label{eq:gaugefixed}
 \partial_tB=\mathcal V_t(B)-d_B\zeta_t(B),
 \qquad B(0)=A_{\mathrm{in}}.
\end{align}
The original equation is degenerate in gauge directions. The added
term supplies the missing second-order part in those directions and
makes the equation for $B$ strictly parabolic. We will remove it
when reconstructing $A$ below.

Write $K_a\omega=[a\wedge\omega]$, so $d_B=d_{A_\star}+K_a$.
All adjoints here still use $g_t$.
We denote by
\begin{align*}
 \Delta_{A_\star,t}=d_{A_\star}^{*_{g_t}}d_{A_\star}+d_{A_\star} d_{A_\star}^{*_{g_t}},
 \qquad F_B=F_{A_\star}+d_{A_\star} a+\tfrac12[a\wedge a].
\end{align*}
Direct substitution in \eqref{eq:gaugefixed} gives
\begin{align}\label{eq:semilinear}
 \partial_ta+\Delta_{A_\star,t}a
       =N_t(a)+\UU_{A_\star+a,\psi_t}(Q_t),
\end{align}
where
\begin{align}
 N_t(a)={}&-d_{A_\star}^{*_{g_t}}F_{A_\star}
 -\tfrac12d_{A_\star}^{*_{g_t}}[a\wedge a]
 -K_a^{*_{g_t}}F_B-[a,d_{A_\star}^{*_{g_t}}a]
 +*_{t}(F_B\wedge d\varphi_t).
\end{align}
In the bracket with a zero form, the remaining one-form index is
retained. The curvature expansion shows that every term in $N_t(a)$
contains at most one derivative of $a$. Lemma~\ref{lem:Umap} shows
that the correction preserves the Sobolev order of $a$. Thus all
second derivatives of the unknown are on the left. The coefficients
of that operator depend only on the prescribed metric and the fixed
reference connection, and its principal symbol is $|\xi|_{g_t}^2 I$.

Let $\mathsf H(t,s)$ be the evolution family for
$\partial_t+\Delta_{A_\star,t}$. Its coefficients are smooth and uniformly elliptic on a fixed compact time interval. Linear nonautonomous parabolic theory gives this family \cite[sections 6.1-6.2]{Lunardi95}. In fixed reference Sobolev norms it satisfies
\begin{align}
 \|\mathsf H(t,s)v\|_{H^6}\leq C\|v\|_{H^6},\qquad
 \|\mathsf H(t,s)v\|_{H^6}
       \leq C(t-s)^{-1/2}\|v\|_{H^5},\quad 0\leq s<t\leq T.
\end{align}
For completeness, these bounds directly follow from energy estimates. For a smooth $v$, put $u(q) = H(q,s)v$ and $y_{j}(q) = \|\mathsf u(q)\|^{2}_{H^{j}}$. Uniform elliptcity, commutation with fixed reference derivatives, and integration by parts gives
\begin{align*}
    y_{j}^{\prime}(q) + c_{j}y_{j+1}(q)\leq C_{j}y_{j}(q), \qquad
    j=5,6
\end{align*}
with constants uniform for $0\leq s\leq q\leq T$. Grownwall's inequality at order six prove sthe first bound above. At order five it gives
\begin{align*}
    y_{5}(t) + c\int^{t}_{s}y_{6}(q)dq \leq Cy_{5}(s)
\end{align*}
Multiplying the order six differential inequality by $q-s$ and dropping its negative coercive term yeilds $\frac{d}{dq}((q-s)y_{6}(q))\leq (1+C_{6}(q-s))y_{6}(q)$. After integration and the preceding order five estimate, $(t-s)y_{6}(t)\leq Cy_{5}(s)$. Taking square roots proves the second bound, and density extends it to $v\in H^{5}$.
The second estimate gains the one derivative needed for the local
nonlinearity. Its time singularity is integrable because  
$\int_0^t(t-s)^{-1/2}\,ds=2\sqrt t$.
This is the reason that the integral equation below is a contraction
for small time. The constants are uniform on a fixed short interval.
Write
$\mathcal N_t(a)=N_t(a)+\UU_{A_\star+a,\psi_t}(Q_t)$.
On a sufficiently small $H^6$ ball, products and Lemma~\ref{lem:Umap}
give
\begin{align}
 \|\mathcal N_t(a)-\mathcal N_t(b)\|_{H^5}
              \leq L\|a-b\|_{H^6},\qquad
 \sup_t\|\mathcal N_t(a)\|_{H^5}\leq M.
\end{align}
Choose a closed ball of radius $\rho$ strictly inside the inverse
neighborhood. For $a_{\mathrm{in}}=A_{\mathrm{in}}-A_\star$, the integral
equation is $a=\mathcal Ta$, where
\begin{align}
 (\mathcal Ta)(t)=\mathsf H(t,0)a_{\mathrm{in}}
        +\int_0^t\mathsf H(t,s)\mathcal N_s(a(s))\,ds.
\end{align}
On $C([0,T],H^6)$, the integral is bounded by $2CM\sqrt T$, and
\begin{align*}
 \sup_{0\leq t\leq T}\|(\mathcal Ta)(t)-(\mathcal Tb)(t)\|_{H^6}
       \leq 2CL\sqrt T\,\sup_{0\leq t\leq T}\|a(t)-b(t)\|_{H^6}.
\end{align*}
First choose $\delta$ so that $C\delta<\rho/4$, and then choose $T$
so that $2CM\sqrt T<\rho/4$ and $2CL\sqrt T<1/2$. The map then sends
the closed ball of radius $\rho$ in $C([0,T],H^6)$ into the ball of
radius $\rho/2$ and is a contraction there. Its unique fixed point is
a solution of \eqref{eq:semilinear}, continuous into $H^6$ and
differentiable into $H^4$.

To prove smoothness on this same interval, use the bound
$\|a\|_{H^6}\leq K$, product estimates and
\eqref{eq:Utame} to obtain, for $r\geq6$,
\begin{align*}
 \|\mathcal N_t(a)\|_{H^{r-1}}
                  \leq C_{r,K}(1+\|a\|_{H^r}).
\end{align*}
We differentiate \eqref{eq:semilinear} using the fixed reference
connection and integrate by parts in the second-order term. Uniform
ellipticity controls $\|a\|_{H^{r+1}}^2$. We pair the remaining term in
$H^{r-1}$ with $a$ in $H^{r+1}$ and use Young's inequality. This gives
an arbitrarily small multiple of $\|a\|_{H^{r+1}}^2$ plus
$C_{r,K}(1+\|a\|_{H^r}^2)$. Absorbing the small multiple yields
\begin{align}\label{eq:highenergy}
 \frac{d}{dt}\|a\|_{H^r}^2+c\|a\|_{H^{r+1}}^2
          \leq C_{r,K}(1+\|a\|_{H^r}^2).
\end{align}
The commutators with the fixed reference derivatives have smooth
coefficients and are absorbed in the same estimate. For each integer
$r>6$, we first solve locally in $H^r$ by the contraction argument
with heat estimates at that order, and uniqueness in $H^6$ identifies
this solution with the one already constructed wherever both are
defined. Applying Gronwall's inequality to \eqref{eq:highenergy} then
bounds its $H^r$ norm in terms of the initial $H^r$ norm and $C_{r,K}$
up to any time before $T$.

The low-order solution stays in the smaller $H^6$ ball chosen above,
so it remains a positive distance from the boundary where the inverse
was defined. On bounded $H^r$ subsets of this neighborhood,
Lemma~\ref{lem:Umap} and the inverse differentiation formula following
Lemma~\ref{lem:D} give uniform local Lipschitz bounds for
$\mathcal N_t:H^r\to H^{r-1}$. These bounds and the distance to the
boundary give a positive restart time for the contraction.

If the $H^r$ solution stopped before $T$, its norm would still be
bounded by \eqref{eq:highenergy}. Restarting sufficiently close to the
stopping time would extend it, a contradiction. Consequently
$a\in C([0,T],H^r)$ for every $r$. This proves higher regularity on
the same interval as the low-order solution. The equation then gives
$a\in C^1([0,T],H^{r-2})$ at each order, and repeated differentiation
in time supplies all time derivatives, including at $t=0$.

To recover the original equation from $B_t$, solve the pointwise ODE
for gauge transformations
\begin{align}
 \dot k_tk_t^{-1}=\zeta_t(B_t),\qquad k_0=\mathrm{id},
 \qquad A_t=k_t^*B_t.
\end{align}
The right-hand side is a smooth section of the adjoint bundle, so the
ODE defines smooth gauge transformations on $[0,T]$. Using the chosen
gauge action convention gives
\begin{align*}
 \partial_tA_t
 &=k_t^{-1}(\partial_tB_t+d_{B_t}(\dot k_tk_t^{-1}))k_t\\
 &=k_t^{-1}\mathcal V_t(B_t)k_t
 =\mathcal V_t(A_t).
\end{align*}
This gives a solution of the original equation with
$A(0)=A_{\mathrm{in}}$, and since $\DD_{B_t,\psi_t}$ is invertible
by the choice of fixed point neighborhood, its gauge conjugate
$\DD_{A_t,\psi_t}$ remains invertible as well and hence we get the existence of solutions. 

\medskip
We now prove uniqueness.
The contraction proves uniqueness of the gauge-fixed equation. We
still need uniqueness for the original equation, because two original
solutions need not come with the same gauge-fixing construction.
Let $A_1,A_2$ be smooth solutions with the same initial value and put
$v=A_1-A_2$. Abbreviate $\chi_j=\chi_{A_j,\psi_t}$ and
$\PP_j=\PP_{A_j,\psi_t}$. The curvature expansion gives
\begin{align*}
 \chi_1-\chi_2=\PP_2v+\tfrac12*_t([v\wedge v]\wedge\psi_t), \quad 
 (\PP_1-\PP_2)w=*_t([v\wedge w]\wedge\psi_t).
\end{align*}
Put $n(v)=\tfrac12*_t([v\wedge v]\wedge\psi_t)$ and
$\mathcal C(v,w)=*_t([v\wedge w]\wedge\psi_t)$. Substituting the
two expansions into \eqref{eq:torsionYM} gives the exact difference
equation
\begin{align*}
 \partial_tv=-\PP_2^2v-\PP_2n(v)
 -\mathcal C(v,\PP_2v)-\mathcal C(v,n(v))
 -\mathcal C(v,\chi_2)+\UU_{A_1,\psi_t}(Q_t)-\UU_{A_2,\psi_t}(Q_t).
\end{align*}
Pair this equation with $v$ in the current $L^2$ inner product.
Self-adjointness gives $-\|\PP_2v\|_2^2$ from the first term and
$-\langle n(v),\PP_2v\rangle_{L^2}$ from the second. On a compact
time interval both solutions are smooth, so $v$ and $\chi_2$ have
bounded pointwise norms. The algebraic estimates
$|n(v)|\leq C|v|^2$ and
$|\mathcal C(v,w)|\leq C|v||w|$, followed by Young's inequality,
bound the next four terms by
$\tfrac12\|\PP_2v\|_2^2+C\|v\|_2^2$.

For the nonlocal term, one may subtract its two defining elliptic
equations. If $w_j=(f_j,u_j)$, their difference satisfies
\begin{align*}
 \DD_1(w_1-w_2)
 =-(\DD_1-\DD_2)w_2
   -(0,*_t((F_{A_1}-F_{A_2})\wedge Q_t)).
\end{align*}
We estimate this source in $H^{-1}$ rather than in $L^2$, so that
we do not need a derivative of $v$ in the final bound. The first
term is algebraic in $v$, with bounded smooth coefficients. In the
curvature difference, $d_{A_2}v$ has $H^{-1}$ norm at most
$C\|v\|_2$, and the quadratic term is bounded by
$C\|v\|_\infty\|v\|_2$. Smoothness on the compact time interval
bounds $\|v\|_\infty$. Thus the whole source has $H^{-1}$ norm at
most $C\|v\|_2$. Applying $\DD_1^{-1}:H^{-1}\to L^2$ gives
$\|u_1-u_2\|_2\leq C\|v\|_2$. This also shows that no assumption
about invertibility along the line segment between the two
connections is needed. Finally, differentiating the metric and
volume in $\|v\|_2^2$ contributes at most $C\|v\|_2^2$. Thus
\begin{align*}
 \frac{d}{dt}\|v\|_2^2\leq C\|v\|_2^2,
\end{align*}
and Gronwall's inequality proves uniqueness for the original equation.

\end{proof}

\subsection{Continuation and transport along ambient paths}
\label{sec:connection-transport}
There are two questions about longer time intervals. We first ask
whether an instanton can be transported along a controlled ambient
$\psi_t$ evolution for all time. We then give sufficient bounds for a general
connection solution to continue past a finite time.

For transport, staying in the reference neighborhood gives an
instanton at each time through the local branch. We also need a
time-dependent gauge to make its velocity Coulomb. The size of this
gauge is controlled by the total amount of ambient motion, measured
by the time integral of the velocity. 
%Staying close to one ambient form does not by itself bound that integral.

Fix an acyclic reference pair $(A_\star,\psi_\star)$ and write
$\mathcal I(\psi)$ for the instanton in the fixed reference Coulomb
slice supplied by Lemma~\ref{lem:IFT}. We have the following result which tells us about the global transport along a small ambient path.

\begin{theorem}
\label{thm:global-instanton-transport}
Let $(A_\star,\psi_\star)$ be as in Lemma~\ref{lem:IFT}. There
are a finite integer $r$ and $\epsilon>0$ such that a smooth closed
positive family $\psi_t$, $t\geq0$, satisfying
\begin{align}
\label{eq:small-instanton-path}
 \sup_{t\geq0}\|\psi_t-\psi_\star\|_{C^r}<\epsilon,\qquad
 \int_0^\infty\|Q_t\|_{C^r}\,dt<\epsilon
\end{align}
has a global instanton lift satisfying
\eqref{eq:horizontaltransport}, starting at
$\mathcal I(\psi_0)$. Assume the analogous higher-norm bounds are
finite. If $\psi_t$ converges smoothly and $Q_t$ is integrable in
every norm, the lift converges smoothly to an acyclic instanton for
the limiting four form. Exponential bounds for the ambient tail
give exponential convergence of the lift. The lift remains uniformly
close to $A_\star$ at any prescribed fixed finite order if the
corresponding bounds in \eqref{eq:small-instanton-path} are small.
\end{theorem}

\begin{proof}
The proof is similar  in spirit to the proof of Theorem~\ref{thm:integrable-coassoc}.  Let $\bar B_t=\mathcal I(\psi_t)$. The local branch is defined for all
$t$ because the whole ambient path stays in its fixed neighborhood.
It is acyclic, with a uniform inverse bound. Differentiating its
instanton equation gives
$\PP_{\bar B,t}\partial_t\bar B=-*_t(F_{\bar B}\wedge Q_t)$.
We solve the equation
\begin{align*}
 d_{\bar B}^{*t}d_{\bar B}q_t
 =d_{\bar B}^{*t}\partial_t\bar B.
\end{align*}
The scalar Laplacian is uniformly invertible because its kernel is
zero at the reference pair and this persists in the fixed
neighborhood. We record why its solution is controlled by the
ambient velocity. The derivative of the slice condition is
$d_{A_\star}^{*_{g_\star}}\partial_t\bar B=0$, so the invertible
linearization $B_\psi$ in Lemma~\ref{lem:IFT} gives
\begin{align*}
 B_{\psi_t}(0,\partial_t\bar B_t)
       =(0,-*_t(F_{\bar B_t}\wedge Q_t)).
\end{align*}
The inverse of $B_{\psi_t}$ gains one derivative and bounds
$\partial_t\bar B_t$ by $Q_t$. The scalar equation for $q_t$ has one
derivative of $\partial_t\bar B_t$ on the right, and its second-order
inverse gains two derivatives. Combining the estimates, for each
fixed $k\geq0$ we obtain
\begin{align}\label{eq:transport-gauge-bound}
 \|\partial_t\bar B_t\|_{C^{k+1,\alpha}}
       +\|q_t\|_{C^{k+2,\alpha}}
 \leq C_k\|Q_t\|_{C^{k+2,\alpha}}.
\end{align}
The constants use bounded ambient and branch coefficients at the
corresponding higher orders and the uniform low-order inverse
bounds. These are available by the regularity argument in
Lemma~\ref{lem:IFT}. The derivative counts in
\eqref{eq:transport-gauge-bound} are deliberately more than is
needed at some orders as they allow us to use one fixed finite norm
when choosing the small neighborhood. Put
$v_t=\partial_t\bar B_t-d_{\bar B_t}q_t$. Then
\begin{align*}
 d_{\bar B_t}^{*t}v_t=0,\qquad
 \PP_{\bar B_t,t}v_t=-*_t(F_{\bar B_t}\wedge Q_t).
\end{align*}
Here $\PP_{\bar B,t}d_{\bar B}q=0$ because $\bar B$ is an instanton.
Thus $v_t=\UU_{\bar B_t,\psi_t}(Q_t)$.

Choose a gauge $g_t$ with
$(\partial_tg_t)g_t^{-1}=-q_t$ and $g_0=1$. Then
$B_t=g_t^*\bar B_t$ has velocity $\UU_{B_t,\psi_t}(Q_t)$.
At each point of $M$, the gauge equation is an ordinary differential
equation in the compact structure group, so it exists for every
finite time. To bound spatial derivatives, we differentiate in local
bundle trivializations. At each order, the highest derivative of
$g_t$ satisfies a linear equation. Its coefficients involve $q_t$ and
the remaining terms involve derivatives of $q_t$ and derivatives of
$g_t$ already controlled at lower orders. The time-integral bounds,
\eqref{eq:transport-gauge-bound}, and Gronwall's inequality therefore
bound all derivatives successively.

Integrating each differentiated equation between two later times
bounds the difference by the tails of the corresponding time
integrals. These tails tend to zero, so $g_t$ is Cauchy in every
$C^k$ norm. The pointwise limit lies in the compact structure group and
its smoothness follows from the derivative convergence. Thus the
limit is a gauge transformation. The same estimates
give
\begin{align*}
 \|g_t-1\|_{C^{k+1}}
 \leq C_k\int_0^t\|Q_q\|_{C^{k+3}}\,dq
\end{align*}
when the relevant integral is small, in local charts about the
identity. The formula for the gauge action then keeps $B_t$ close
to $A_\star$ in $C^k$. Taking the finite integer $r$ large enough
gives the first assertion. The higher-norm integrability gives
smooth convergence, and exponential bounds for the integrals
from $t$ to infinity give exponential convergence.
The limiting instanton is acyclic because the uniform reference
neighborhood has a uniform inverse bound, and conjugating by a
gauge transformation preserves the spectrum. Uniqueness follows
from Theorem~\ref{thm:nearby-existence}.
\end{proof}

We now give an analogue of Proposition~\ref{prop:coassoc-continuation} for the coupled system of the flow of $4$-forms and the Yang-Mills flow.

\begin{proposition}
\label{prop:connection-continuation}
Suppose $\psi_t$ extends smoothly across a finite $T$ and remains
positive and closed. A solution of \eqref{eq:coupledYM}
extends smoothly across $T$, modulo a time-dependent gauge, if in
such a gauge
\begin{align*}
 \sup_{t<T}\|A_t-A_{\mathrm{ref}}\|_{H^s}<\infty,
 \qquad
 \sup_{t<T}\|\DD_{A_t,\psi_t}^{-1}\|_{L^2\to L^2}<\infty
\end{align*}
for one integer $s\geq8$. Here $A_{\mathrm{ref}}$ is fixed and smooth.
Equivalently, the displayed bounds may be imposed on a gauge-fixed
representative with its gauge term included in the equation.
\end{proposition}

\begin{remark}
This is a sufficient finite-order criterion using a connection bound
in a gauge and the inverse bound.
\end{remark}

\begin{proof}
Since the proof if similar to the proof of Proposition~\ref{prop:coassoc-continuation}, we omit some details. We prove continuation by finding a local existence time that is
uniform as the starting time approaches $T$. Use the same fixed
smooth reference connection to express the gauge-fixed local
equations. The assumed $H^s$ bound controls the coefficients of
$\DD$ at low order and in the higher Sobolev norms. The uniform $L^2$ inverse bound and the
elliptic estimate successively give uniform bounds
$\DD^{-1}:H^{s-1}\to H^s$. To see that no neighborhood at every order
is required, first use the $L^2$ bound in the order-one elliptic
estimate and then apply the higher estimates inductively. Sobolev
multiplication bounds all lower-order coefficients by the given
$H^s$ bound. The resolvent identity then gives a uniform radius on which the
inverse remains bounded. The curvature expansion and the estimates
for $\UU$ give uniform bounds and Lipschitz constants for the
nonlinearity $H^s\to H^{s-1}$ on that radius. The heat estimates
used in Theorem~\ref{thm:nearby-existence} are uniform because the
ambient metric extends smoothly across $T$. Consequently the
contraction argument has an existence time bounded below by one
positive number, independently of its starting time $t_0<T$.

If the assumed bound is obtained in a time-dependent gauge, fix
its value at $t_0$ and use that single gauge transformation on the
initial connection for the new local solution. Gauge covariance
and uniqueness identify the new solution with the original one
on their common interval, modulo this fixed transformation. Choose
$t_0$ closer to $T$ than the uniform existence time. The new
solution then continues the geometric solution across $T$.
Parabolic and elliptic regularity give smoothness there. This proves
continuation of the geometric solution. If the gauge transformations
themselves extend smoothly, the original representative extends
as well.
\end{proof}

\subsection{Nonlinear stability}
\label{sec:instanton-stability}

In this section, we first show that, near an acyclic instanton, the defect controls the connection after a Coulomb gauge has been chosen. We then compare a
solution with a nearby moving instanton. A fixed positive elliptic
operator gives decay of the difference in one Sobolev norm. We use this
decay to keep the solution in the neighborhood where the correction is
defined, and then use the parabolic equation to obtain decay of every
higher derivative. The initial connection in this argument need not
be an instanton.

Fix a closed positive $\psi_\star$ and an acyclic
$\psi_\star$-instanton $A_\star$. We use an integer $s\geq8$ throughout this subsection. This is a
sufficient regularity assumption, and we do not claim that it is
optimal. In dimension seven, $H^{s-1}$ is a Banach algebra and
$H^s$ embeds continuously into $C^4$. We also remain above the
$H^6$ regularity used in
Theorem~\ref{thm:nearby-existence}.
The higher coefficient bounds in
\eqref{eq:background-small} are separate assumptions on the
prescribed background and they do not follow from this Sobolev
embedding. Use the notation
\begin{align}
\label{eq:reference-L}
 \mathscr L_\star=\PP_\star^2+d_{A_\star}d_{A_\star}^{*},\qquad
 \PP_\star=\PP_{A_\star,\psi_\star}.
\end{align}
At an instanton $\PP_\star d_{A_\star}=0$ and
$d_{A_\star}^*\PP_\star=0$, so the mixed terms in $\DD_\star^2$
vanish. Its one-form block is therefore $\mathscr L_\star$.
For a one form $a$, integration by parts gives
\begin{align*}
 \langle\mathscr L_\star a,a\rangle_{L^2}
 =\|\PP_\star a\|_{L^2}^2+\|d_{A_\star}^*a\|_{L^2}^2.
\end{align*}
The first term detects a change in the instanton equation, while the
second detects the component in the gauge directions. If both vanish,
acyclicity gives $a=0$. Thus $\mathscr L_\star$ is positive,
self-adjoint and strongly elliptic, with no kernel. Since the manifold is compact, we get that there exists a lowest eigenvalue $\nu>0$ of $\mathscr L$. For integers $k\geq0$, we use the
following fixed Sobolev norms
\begin{align}
\label{eq:spectral-norms}
 |a|_k=\|\mathscr L_\star^{k/2}a\|_{L^2(g_\star)}.
\end{align}
Elliptic regularity and positivity make these norms equivalent to the
ordinary $H^k$ norms. To see the estimate that we will use, expand $a$
in an orthonormal eigenbasis of $\mathscr L_\star$. If its eigenvalues
are $\mu_j\geq\nu$ and the coefficients of $a$ are $a_j$, then
$|a|_k^2=\sum_j\mu_j^k|a_j|^2$. Comparing these sums term by term gives
\begin{align*}
 |a|_{k+1}^2\geq\nu|a|_k^2.
\end{align*}

\begin{lemma}
\label{lem:defect-distance}
Every sufficiently close $H^s$ connection is gauge equivalent, by a
gauge transformation close to the identity, to $A_\star+a$ with
$d_{A_\star}^*a=0$. In that slice,
\begin{align*}
 \|a\|_{H^s}\leq C\|\chi_{A_\star+a,\psi_\star}\|_{H^{s-1}}.
\end{align*}
If $a$ is sufficiently small in $C^0$, one also has
$\|a\|_{H^1}\leq C\|\chi_{A_\star+a,\psi_\star}\|_{L^2}$.
\end{lemma}

\begin{proof}
We first choose the gauge. Consider
$d_{A_\star}^*((\exp q)^*(A_\star+a)-A_\star)$ as a function of the
scalar gauge parameter $q$ and the one form $a$. Its derivative with
respect to $q$ at zero is $d_{A_\star}^*d_{A_\star}$. Acyclicity gives
$\ker d_{A_\star}|_{\Omega^0(E)}=0$, so this scalar Laplacian is an
isomorphism $H^{s+1}\to H^{s-1}$. The gauge action is smooth between
these spaces for $s>\frac72$. The implicit function theorem therefore
provides a unique sufficiently small $q$ for which the transformed
connection satisfies the stated Coulomb condition.

In the slice, $\DD_\star(0,a)=(0,\PP_\star a)$, so the elliptic inverse gives
$\|a\|_{H^s}\leq C\|\PP_\star a\|_{H^{s-1}}$. Since the instanton defect satisfies 
$\chi_{A_\star+a}=\PP_\star a+\frac12*([a\wedge a]\wedge\psi_\star)$,
Sobolev multiplication gives
\begin{align*}
 \|a\|_{H^s}
 \leq C\|\chi_{A_\star+a}\|_{H^{s-1}}+C\|a\|_{H^s}^2.
\end{align*}
Choose the ball i the Sobolev space $H^s$ so that $C\|a\|_{H^s}\leq \frac 12$.
The last term is then at most $\frac12\|a\|_{H^s}$ and can be moved
to the left-hand side. This proves the first estimate. At order one,
the quadratic term is bounded in $L^2$ by
$C\|a\|_{C^0}\|a\|_{L^2}$. Since $\|a\|_{L^2}\leq\|a\|_{H^1}$,
sufficiently small $C^0$ norm gives the same absorption and proves the
second estimate. 
\end{proof}

\begin{remark}
The gauge condition is essential in the proof of Lemma~\ref{lem:defect-distance}. A nontrivial pure
gauge change can have zero defect and a nonzero connection difference. Lemma~\ref{lem:defect-distance} gives a local comparison between the instanton defect and the distance from the reference instanton after gauge fixing. We record it for this geometric interpretation.
The dynamical proof below instead estimates the connection
difference directly and does not use this lemma.

\end{remark}

For the evolution argument, the reference instanton is allowed to
move with the ambient $44$-form. We subtract its equation before estimating
the perturbation. This removes the forcing that would otherwise be
present even when the difference of the two connections is zero.

\begin{lemma}
\label{lem:nonlinear-estimates}
Suppose $B_t$ is a smooth family of $\psi_t$-instantons satisfying
$\partial_tB_t=\UU_{B_t,\psi_t}(Q_t)$. For $\widehat A=B_t+a$, add the
gauge term $-d_{\widehat A}d_{B_t}^{*t}a$. Then $a$ satisfies 
\begin{align}
\label{eq:perturbation-equation}
 \partial_ta=-\mathscr L_ta+N_t(a)
       +\UU_{B_t+a,\psi_t}(Q_t)-\UU_{B_t,\psi_t}(Q_t),
 \qquad \text{where} \quad 
 \mathscr L_t=\PP_{B_t,t}^2+d_{B_t}d_{B_t}^{*t}.
\end{align}
On a small $H^s$ ball, $s\geq8$, we have
\begin{align}
\label{eq:nonlinear-bounds}
 \|N_t(a)\|_{H^{s-1}}\leq C\|a\|_{H^s}^2.
\end{align}
Moreover, if
\begin{align}
\label{eq:background-small}
 \|\psi_t-\psi_\star\|_{C^{s+2}}
 +\|B_t-A_\star\|_{C^{s+1}}
 +\|Q_t\|_{C^{s+1}}\leq\epsilon,
\end{align}
then
\begin{align}
\label{eq:operator-small}
 \|(\mathscr L_t-\mathscr L_\star)a\|_{H^{s-1}}
 &\leq C\epsilon\|a\|_{H^{s+1}},\\
 \|\UU_{B_t+a,\psi_t}(Q_t)-\UU_{B_t,\psi_t}(Q_t)\|_{H^{s-1}}
 &\leq C\epsilon\|a\|_{H^s}.
\nonumber
\end{align}
\end{lemma}

\begin{proof}
We first expand the equation and identify the terms that are at least
quadratic in $a$. Put
$\mathcal C_t(a,v)=*_t([a\wedge v]\wedge\psi_t)$ and
$n_t(a)=\frac12*_t([a\wedge a]\wedge\psi_t)$. Since $B_t$ is an
instanton, the curvature expansion gives
$\chi_{B_t+a,\psi_t}=\PP_{B_t,t}a+n_t(a)$, and
$\PP_{B_t+a,t}=\PP_{B_t,t}+\mathcal C_t(a,\cdot)$.
The gauge term is
$-d_{B_t}d_{B_t}^{*t}a-[a,d_{B_t}^{*t}a]$.
The linear terms combine to $-\mathscr L_ta$. Subtracting
$\partial_tB_t=\UU_{B_t,\psi_t}(Q_t)$ leaves the difference of the two
corrections in \eqref{eq:perturbation-equation}. The remaining local
terms are exactly
\begin{align*}
 N_t(a)=-\PP_{B_t,t}n_t(a)
       -\mathcal C_t(a,\PP_{B_t,t}a)
       -\mathcal C_t(a,n_t(a))
       -[a,d_{B_t}^{*t}a].
\end{align*}
Every summand has either two factors of $a$ and at most one derivative,
or three factors of $a$ and no derivatives. Because $H^{s-1}$ is an
algebra, the quadratic terms are bounded by $C\|a\|_{H^s}^2$.
On the chosen small ball the cubic term obeys the same bound. This
proves \eqref{eq:nonlinear-bounds}.

We next compare the second-order operators. In fixed coordinate and
bundle charts, write
\begin{align*}
 (\mathscr L_t-\mathscr L_\star)a
 =\sum_{|\beta|\leq2}c_\beta(t,x)\partial^\beta a.
\end{align*}
The coefficients are smooth functions of the background forms,
connections, and their derivatives. Each coefficient vanishes when
$(\psi_t,B_t)=(\psi_\star,A_\star)$. The derivative orders in
\eqref{eq:background-small} therefore bound their $C^{s-1}$ norms by
$C\epsilon$. When we apply up to $s-1$ derivatives to
$c_\beta\partial^\beta a$, Leibniz's rule puts at most $s-1$ derivatives
on the coefficient and at most $s+1$ derivatives on $a$.
Summing these product bounds gives $C\epsilon\|a\|_{H^{s+1}}$,
which is the first estimate in \eqref{eq:operator-small}.

For the correction, put $w_A=-\DD_{A,\psi}^{-1}
(0,*(F_A\wedge Q))$. Subtracting the equations at $B+a$ and $B$
gives
\begin{align}\label{eq:correction-difference}
 \DD_{B+a,\psi}(w_{B+a}-w_B)
 ={}-(\DD_{B+a,\psi}-\DD_{B,\psi})w_B-(0,*((d_Ba+\tfrac12[a\wedge a])\wedge Q)),
\end{align}
and estimates of Lemma~\ref{lem:Umap} give the second line of
\eqref{eq:operator-small}.
\end{proof}

We will also need estimates at orders above $s$ without assuming
smallness of the higher background norms. For every fixed $r\geq s$,
the same product estimates and \eqref{eq:correction-difference} give
\begin{align}\label{eq:higher-remainder-bound}
 \|N_t(a)\|_{H^{r-1}}
 +\|\UU_{B_t+a,\psi_t}(Q_t)-\UU_{B_t,\psi_t}(Q_t)\|_{H^{r-1}}
 \leq C_r\|a\|_{H^r}
\end{align}
on the fixed small $H^s$ ball. Here the constants may depend on bounded
higher derivatives of the reference family; these derivatives need not
be small. We explain why the highest norm of $a$ occurs only linearly.
The low-order inverse first bounds $w_{B+a}-w_B$ in terms of
$\|a\|_{H^s}$. Moving the zeroth-order difference
$\DD_{B+a}-\DD_B$ to the right lets us apply higher elliptic estimates
successively, as in the proof of \eqref{eq:Utame}. In the product
terms, a highest derivative of $a$ is multiplied by a controlled
low-order factor. The other terms have their highest derivative on
the smooth reference coefficients or on a lower-order unknown already
controlled by the induction. The same product estimates apply to
$N_t(a)$. This proves \eqref{eq:higher-remainder-bound} with constants
depending on the fixed low-order ball and the higher coefficient
bounds. Every term still contains a factor of $a$, so no constant
forcing term appears in this estimate.

\begin{theorem}[Attraction to a moving acyclic instanton]
\label{thm:moving-stability}
Let $\psi_t$ and $B_t$ be smooth for all $t\geq0$, with $\psi_t$
closed and positive, $B_t$ an instanton, and
$\partial_tB_t=\UU_{B_t,\psi_t}(Q_t)$. Assume that $\psi_t$, $B_t-A_\star$ and $Q_t$ are uniformly bounded
in every spatial $C^k$ norm. For some sufficiently small $\epsilon>0$, suppose
\eqref{eq:background-small} holds for all $t$. If the smooth
initial connection $A_0$ is sufficiently close to $B_0$ in $H^s$,
$s\geq8$, then \eqref{eq:coupledYM} has a unique global
smooth solution. Moreover, there are smooth gauge transformations $g_t$, with
$g_0=1$, such that
\begin{align}
\label{eq:tracking}
 \|g_t^{-1*}A_t-B_t\|_{H^s}\leq Ce^{-\gamma t}\|A_0-B_0\|_{H^s}
\end{align}
for some $\gamma>0$. The difference decays exponentially in every
$C^k$ norm after positive time, and $g_t$ converges smoothly to a
gauge transformation $g_\infty$. If $B_t$ converges smoothly and
exponentially to $B_\infty$, then $A_t$ converges smoothly and
exponentially to $g_\infty^*B_\infty$.
\end{theorem}

\begin{proof}
We first solve \eqref{eq:perturbation-equation} for a short time with
$a(0)=A_0-B_0$. The parabolic construction in
Theorem~\ref{thm:nearby-existence} applies as replacing a fixed reference
connection by the prescribed smooth $B_t$ changes only the coefficients. We now estimate this solution while it remains in a small $H^s$ ball, and then prove that it cannot leave that ball.

We use the fixed norm \eqref{eq:spectral-norms}, so differentiating the
norm produces no metric-variation terms. The useful derivative pairing
is
\begin{align*}
 |\langle a,b\rangle_s|
 =|\langle \mathscr L_\star^{(s+1)/2}a,\mathscr L_\star^{(s-1)/2}b\rangle_{L^2}|
 \leq |a|_{s+1}|b|_{s-1}
\end{align*}
This identity follows by moving powers of the fixed self-adjoint
operator $\mathscr L_\star$ between the two factors. It allows us to
pair a term controlled in $H^{s-1}$ with the dissipative quantity
$|a|_{s+1}$. In particular, the second-order operator difference is
bounded by $C\epsilon|a|_{s+1}^2$, the local nonlinearity by
$C|a|_s^2|a|_{s+1}$, and the correction difference by
$C\epsilon|a|_s|a|_{s+1}$. Thus, we obtain from Lemma~\ref{lem:nonlinear-estimates}
\begin{align}
\label{eq:stability-energy}
 \frac12\frac{d}{dt}|a|_s^2
 &\leq-|a|_{s+1}^2
       +C\epsilon|a|_{s+1}^2
       +C|a|_s^2|a|_{s+1}
       +C\epsilon|a|_s|a|_{s+1}\\
 &\leq-\big(1-C_1\epsilon-C_2|a|_s\big)|a|_{s+1}^2,
\nonumber
\end{align}
where in the last line we used $|a|_s\leq\nu^{-1/2}|a|_{s+1}$.
Choose a radius $\rho$ inside the neighborhood of definition of the
correction and choose $\epsilon$ so that
$C_1\epsilon+C_2\rho<\frac 12$. Take $|a(0)|_s<\rho$.
As long as $|a(t)|_s<\rho$, the preceding inequality implies
\begin{align*}
 \frac12\frac{d}{dt}|a|_s^2
 \leq-\frac12|a|_{s+1}^2
 \leq-\frac\nu2|a|_s^2
\end{align*}
Integration then gives $|a(t)|_s\leq e^{-\nu t/2}|a(0)|_s<\rho$.
If a first exit time existed, continuity would give $|a|_s=\rho$ there,
which contradicts the bound $|a(t)|_s < \rho$. Thus the small ball is preserved throughout
the existence interval. This argument proves the required smallness.

Together with the assumed bound on $B_t-A_\star$, this keeps the
gauge-fixed connection in a fixed neighborhood of $A_\star$ where
$\DD$ has a uniformly bounded inverse. The local parabolic argument
then has a uniform positive restart time on each finite time interval.
If the maximal time were finite, we could restart sufficiently close
to it and continue beyond it. This contradiction proves global
existence. The same energy estimate gives exponential decay in $H^s$.

The higher-order estimates are obtained in the same as we obtained \eqref{eq:highenergy} and we get 
\begin{align}\label{eq:higher-stability-energy}
 y_r'(t)+c_r y_{r+1}(t)\leq C_r y_r(t),\qquad r\geq s,
\end{align}
where $c_r>0$ and the constants are independent of $t$.

The right-hand side of this inequality is not negative, so applying
Gronwall at order $r$ alone does not prove decay. We instead use the
negative term to gain one derivative at a time on shorter time
intervals, starting from the decay already proved in $H^s$.

Fix an interval $[u,u+1]$ and let
$M_u=\sup_{u\leq q\leq u+1}y_s(q)$. Integration of
\eqref{eq:higher-stability-energy} at $r=s$ gives
$\int_u^{u+1}y_{s+1}(q)\,dq\leq C M_u$.
Choose a smooth nondecreasing time cutoff $\eta$ which vanishes at
$u$ and is one on $[u+1/2,u+1]$. Multiplication of the inequality
at order $s+1$ by $\eta$ gives
\begin{align*}
 (\eta y_{s+1})'+c_{s+1}\eta y_{s+2}
 \leq (C_{s+1}\eta+\eta')y_{s+1}.
\end{align*}
Integrating this inequality first from $u$ to any time in $[u+1/2,u+1]$ allows us to bound $y_{s+1}$ there, and then we can integrate to $u+1$ to bound the integral of
$y_{s+2}$ on that final half interval. The initial term vanishes
because $\eta(u)=0$. Both terms on the right are bounded by the
preceding integral estimate, so
\begin{align*}
 \sup_{u+1/2\leq q\leq u+1}y_{s+1}(q)
 +\int_{u+1/2}^{u+1}y_{s+2}(q)\,dq\leq C M_u.
\end{align*}
For the next derivative order, choose a cutoff supported in the
interval where we have just controlled the integral of $y_{s+2}$.
Continue with cutoffs which become one at
$u+1-2^{-j}$, for $j=2,3,\ldots$. Their derivatives have larger
bounds as $j$ increases, but only finitely many steps are needed
for any fixed $r$. At each step the integral obtained at the
preceding order bounds the new cutoff term. We therefore obtain
control of the $r$th norm at the final endpoint from the supremum
of the $s$th norm on the full interval. Taking $u=t-1$ gives
\begin{align}\label{eq:smoothing-decay}
 \|a(t)\|_{H^r}
 \leq C_r\sup_{t-1\leq q\leq t}\|a(q)\|_{H^s}
 \leq C_r e^{-\gamma t}\|a(0)\|_{H^s},\qquad t\geq1.
\end{align}
Because the energy norms are fixed, this argument does not differentiate
the high-order background coefficients in time. Sobolev embedding now
gives exponential decay in every spatial $C^k$ norm. On the initial
bounded time interval, smooth initial data and
\eqref{eq:higher-stability-energy} give finite higher norms. Thus
we have one smooth solution for all time, even though only the
initial $H^s$ norm was required to be small.

It remains to remove the gauge term and identify the limit in the
original equation. Set $\widehat A=B_t+a$ and
$\xi_t=d_{B_t}^{*t}a$. Let $g_t$ solve
$\dot g_tg_t^{-1}=\xi_t$, $g_0=1$, and put
$A_t=g_t^*\widehat A_t$. The derivative of this time-dependent gauge
transformation adds precisely the term that was subtracted to make
the equation parabolic
\begin{align*}
 \partial_tA_t
 =g_t^{-1}\bigl(\partial_t\widehat A_t+d_{\widehat A_t}\xi_t\bigr)g_t
 =\mathcal V_t(A_t),
\end{align*}
so this recovers the original equation with the same initial value.
By \eqref{eq:smoothing-decay}, $\xi_t=d_{B_t}^{*t}a$ decays
exponentially in every spatial norm. The gauge equation is global
\begin{align*}
 \|g_t-g_\infty\|_{C^k}
 \leq C_k\int_t^\infty\|\xi_q\|_{C^k}\,dq
 \leq C_k e^{-\gamma t},\qquad t\geq1,
\end{align*}
giving exponential convergence. Since
$g_t^{-1*}A_t=B_t+a_t$, the estimate for $a_t$ proves
\eqref{eq:tracking}. If $B_t$ converges exponentially to $B_\infty$,
continuity of the gauge action in the corresponding higher norms gives
$A_t\to g_\infty^*B_\infty$ with the same type of decay. The limiting
representative need not be $B_\infty$ itself. Uniqueness follows from
the smooth uniqueness for the original equation already proved.
\end{proof}

We can finally give a proof of our Theorem~\ref{thm:static-stability}.

	\begin{proof}[Proof of Theorem~\ref{thm:static-stability}]
		We apply Theorem~\ref{thm:moving-stability} to the constant
		background
		\begin{align*}
			\psi_t=\psi_\star,\qquad B_t=A_\star
			\qquad\text{for all }t\geq0.
		\end{align*}
		Then $Q_t=\partial_t\psi_t=0$. Since $A_\star$ is an acyclic
		$\psi_\star$-instanton, each $B_t$ is an acyclic
		$\psi_t$-instanton. Moreover,
		\begin{align*}
			\partial_tB_t=0=\UU_{A_\star,\psi_\star}(0),
		\end{align*}
		so the required evolution equation for the reference
		connection is satisfied. The left-hand side of
		\eqref{eq:background-small} is zero, and all the required
		spatial bounds are uniform because the background is smooth
		and independent of time.
		
		Since $\UU_{A,\psi_\star}(0)=0$, equation
		\eqref{eq:coupledYM} reduces to
		\begin{align*}
			\partial_tA=-d_A^*F_A+*(F_A\wedge d\varphi_\star).
		\end{align*}
		Theorem~\ref{thm:moving-stability} therefore gives a unique
		global smooth solution for every sufficiently small smooth
		$H^s$ perturbation of $A_\star$. Since $B_t=A_\star$ for all
		time, its limiting connection is $B_\infty=A_\star$.
		The convergence conclusion of that theorem gives smooth
		exponential convergence of $A_t$ to $g_\infty^*A_\star$.
		\end{proof}
	
	\begin{remark}
Only $d\psi_\star=0$ is required in the above proof we have not assumed
$d\varphi_\star=0$. Thus the background need not be
torsion-free.
	\end{remark}

\subsection{Modified Laplacian coflow}\label{sec:coflow}
We now take the ambient family to solve Grigorian's modified
Laplacian coflow \cite{Gri}. Since its equation is independent of
the connection, we first solve the ambient problem and then solve
the connection equation along that solution. We explain the
conversion of conventions and the ambient estimates used below.

We first note that \cite[eq. (2.2)]{Gri} uses
$(u\ip\varphi)\wedge(v\ip\varphi)\wedge\varphi
=6g(u,v)\vol_+$, whereas \eqref{eq:mainconv} has the opposite sign.
Keeping $\varphi$ and $g$ fixed therefore gives
$\vol_+=-\vol_g$, $*_+=-*$ and $\psi_+=-\psi$.
Equation (2.10) of \cite{Gri} defines the full torsion by
$\nabla_a\varphi_{bcd}=(T_+)_a{}^e(\psi_+)_{ebcd}$,
with no extra numerical factor. Our definition consequently gives
$T_+=-T$, and hence $\tr T_+=-\tr T$. Also, \cite[eq. (1.4)]{Gri}
\begin{align*}
 \partial_t\psi_+=\Delta\psi_+
                         +2d((c-\tr T_+)\varphi),\qquad c>0,
\end{align*}
where $c$ is his positive constant $A$. We negate the equation and
substitute $\psi_+=-\psi$ and $\tr T_+=-\tr T$.
The Hodge Laplacian is unchanged by orientation reversal, since the
two changes of sign in its star operators cancel. We therefore obtain
exactly \eqref{eq:modifiedcoflow}. The same algebraic conversion applies when $c=0$.

\begin{proof}[Proof of Corollary~\ref{cor:coflow}]
Grigorian introduces a positive parameter $A$ (which we denote by $c$ here) in the modification of the Laplacian coflow, but the local existence argument also applies to every fixed $c\in\mathbb R$. Indeed, the parameter occurs only
in the first-order term $-2c\,d\varphi$, so it does not change
the second-order symbol of the gauge-fixed equation, see also the
gauge-fixing argument in \cite[Theorems~5.1 and~6.9]{Gri} and \cite[Introduction]{BVStab20}.
Thus the short-time theory used here includes both $c=0$ and
negative values of $c$.

The velocity of the flow is
$Q=\QQ_c(\psi)$, which is exact, since both terms defining $\QQ_c$ are
exterior derivatives on closed forms. In particular, $dQ=0$.
Theorem~\ref{thm:transport} therefore gives preservation of the
instanton equation for every existing smooth coupled solution
while the correction is defined.

To construct the coupled solution, first solve the ambient
equation and then apply Theorem~\ref{thm:instanton} to this
prescribed family, starting from the given acyclic instanton.
We shorten the common interval if necessary so that $\DD$ remains
invertible. Uniqueness of the ambient solution and uniqueness of
the connection evolution give uniqueness of the coupled solution.
\end{proof}

For a prescribed smooth solution of unmodified coflow, the same
preservation calculation applies with $Q=\Delta_\psi\psi$. However, the well-posedness of the unmodified Laplacian coflow is still an open problem. For the global result below, we use modified coflow with $c=0$. The next lemma gives the quantitative estimates needed for the global coupled argument. It uses the gauge-fixed equation of \cite{BVStab20} and controls the solution and its total motion
using one finite-order initial norm.

\begin{lemma}
	\label{lemma-ambient-coflow}
	Let $\bar\psi$ be a torsion-free structure on compact $M$. Fix an
	integer $r\geq0$ and $0<\alpha<1$, and put $N=r+8$. There are
	constants $\delta,C,\sigma>0$ with the following property. If $\psi_0$
	is smooth, closed and positive, $[\psi_0]=[\bar\psi]$, and
	$\|\psi_0-\bar\psi\|_{C^{N,\alpha}}<\delta$, then parameter-zero
	modified coflow has a unique global smooth solution and
	\begin{align}
	\sup_{t\geq0}\|\psi_t-\bar\psi\|_{C^{r,\alpha}}
		+\sup_{t\geq0}e^{\sigma t}\|\partial_t\psi_t\|_{C^{r,\alpha}}
		+\int_0^\infty\|\partial_t\psi_t\|_{C^{r,\alpha}}\,dt
		\leq C\|\psi_0-\bar\psi\|_{C^{N,\alpha}}.
		\label{eq-ambient-coflow-bounds}
	\end{align}
	The solution converges smoothly and exponentially to
	$F_\infty^*\bar\psi$ for a diffeomorphism $F_\infty$ isotopic to the
	identity. Every higher spatial norm of the solution is uniformly
	bounded, and every higher spatial norm of its velocity has an
	exponentially decaying bound. 
    %The constants at these higher orders may depend on higher norms of the smooth initial data. No smallness at these higher orders is assumed and all norms use the fixed metric of $\bar\psi$.
\end{lemma}

\begin{proof}
	We use the gauge fixing in \cite[Lemma~3.3]{BVStab20}, translated to
	our orientation convention as above, and take the background
	connection in its DeTurck vector field to be the Levi--Civita
	connection of $\bar\psi$. Write $\bar\nabla$ and $\bar\Delta$ for this
	connection and its nonnegative Hodge Laplacian. If
	$\widetilde\psi=\bar\psi+h$ is closed, the gauge-fixed equation is
	\begin{align*}
		\partial_t h&=-\bar\Delta h+d\bigl(\mathcal R(h)\bar\nabla h\bigr),
		&\mathcal R(0)&=0.
	\end{align*}
	Here $\mathcal R(h)$ is a smooth bundle homomorphism depending pointwise on
	$h$. We explain why the remainder has this form, since a pointwise
	quadratic bound alone would not suffice for the estimates below.
	On closed four forms, the gauge-fixed velocity is the exterior
	derivative of
	\begin{align*}
		*_{\psi}d\varphi
		-2(\operatorname{tr}T_{\psi})\varphi
		+V(\varphi)\mathbin{\lrcorner}\psi,
		\qquad \varphi=*_{\psi}\psi.
	\end{align*}
	Each summand is linear in $\bar\nabla h$, with coefficients depending
	smoothly on $h$: the reference structure is parallel, the Hodge
	star and $\varphi$ depend algebraically on $\psi$, and $V$ is linear
	in the difference of the two Levi--Civita connections. Subtract the
	linearization of this first-order expression at $h=0$. The remainder
	is $\mathcal R(h)\bar\nabla h$, with $\mathcal R(0)=0$. The exterior derivative of its linear part is $-\bar\Delta h$ on closed variations by the  linearization formula. This also explains the nonlinear expansion
	used in \cite[proof of Theorem~3.4, Step~2]{BVStab20}.
	
	For local existence we may regard the above equation as an
	equation for all four forms. If $h$ is small in $C^0$, its principal
	symbol is a small perturbation of that of $-\bar\Delta$, so this
	extension is uniformly parabolic. The usual local existence theorem
	for smooth quasilinear parabolic systems gives a smooth solution.
	It stays closed because
	\begin{align*}
		\partial_t(dh)&=-\bar\Delta(dh), & dh(0)&=0.
	\end{align*}
	Its harmonic part is constant, since both $\bar\Delta h$ and the
	exact remainder are orthogonal to harmonic forms. The initial
	perturbation is exact, so the solution remains exact. Consequently
	this extended equation agrees with the gauge-fixed coflow along
	the solution we have constructed.
	
	Choose $m=r+8$. We first use the fixed Sobolev norms
	\begin{align*}
		\|u\|_{H^j}&=\|(1+\bar\Delta)^{j/2}u\|_{L^2}.
	\end{align*}
	The same notation is used on three and four forms, with the Hodge
	Laplacian of the appropriate degree. These norms are equivalent to
	the usual Sobolev norms. Let $\lambda>0$ be the first eigenvalue of
	$\bar\Delta$ on exact four forms. Since $d$ commutes with the Hodge
	Laplacian, exactness gives
	\begin{align*}
		\|d^*(1+\bar\Delta)^{m/2}h\|_{L^2}^2
		&\geq\lambda\|h\|_{H^m}^2,\\
		\|h\|_{H^{m+1}}^2
		&=\|h\|_{H^m}^2+
		\|d^*(1+\bar\Delta)^{m/2}h\|_{L^2}^2.
	\end{align*}
	All adjoints in this proof use the fixed background metric.
	
	Because $m>7/2+1$, Sobolev multiplication and the chain rule give,
	on a fixed small $H^m$ ball,
	\begin{align*}
		\|\mathcal R(h)\bar\nabla h\|_{H^m}
		&\leq C\bigl(
		\|h\|_{C^0}\|h\|_{H^{m+1}}
		+\|h\|_{H^m}\|\bar\nabla h\|_{C^0}\bigr)\\
		&\leq C\|h\|_{H^m}\|h\|_{H^{m+1}}.
	\end{align*}
	Pair the equation with $(1+\bar\Delta)^m h$ and integrate $d$ by
	parts. The preceding estimates and exactly the same argument as in the beginning of the proof of Theorem~\ref{thm:moving-stability} proves exponential decay of $h(t)$ and thus the solution cannot leave the
	ball during its existence interval.
    
		We next establish higher estimates using only this one smallness
	assumption. For this calculation define the fixed derivative norms
	\begin{align*}
		\|h\|_{j,\bar\nabla}^2
		&=\sum_{\ell=0}^j\|\bar\nabla^\ell h\|_{L^2}^2.
	\end{align*}
	They are equivalent to the spectral Sobolev norms used above.
	We use these derivative norms to keep track of the coefficient of
	the highest derivatives in the energy estimate.
	In differentiating the remainder, separate the term with all
	derivatives on $\bar\nabla h$:
	\begin{align*}
		\bar\nabla^j(\mathcal R(h)\bar\nabla h)
		&=\mathcal R(h)\bar\nabla^{j+1}h
		+[\bar\nabla^j,\mathcal R(h)]\bar\nabla h.
	\end{align*}
	After integration by parts, the first term and its lower-order
	counterparts contribute at most
	\begin{align*}
		C\|\mathcal R(h)\|_{C^0}
		\sum_{\ell=0}^j\|\bar\nabla^{\ell+1}h\|_{L^2}^2.
	\end{align*}
	The constant here is independent of $j$. We absorb this sum directly
	into the same derivative sum supplied by the leading Laplacian
	term, before using any equivalence with spectral Sobolev norms.
	Thus the absorption uses the single fixed $C^0$ smallness chosen
	above.
	Leibniz's rule and the Sobolev product estimates give
	\begin{align*}
		\|[\bar\nabla^j,\mathcal R(h)]\bar\nabla h\|_{L^2}
		&\leq C_j\bigl(
		\|\bar\nabla \mathcal R(h)\|_{C^0}\|\bar\nabla h\|_{H^{j-1}}
		+\|\mathcal R(h)\|_{H^j}\|\bar\nabla h\|_{C^0}\bigr)\\
		&\leq C_j\|h\|_{C^1}\|h\|_{j,\bar\nabla}.
	\end{align*}
Again the same proof as for \eqref{eq:highenergy} gives
    \begin{align*}
		\frac{d}{dt}\|h\|_{j,\bar\nabla}^2
		+c_j\|h\|_{j+1,\bar\nabla}^2
		&\leq C_j\|h\|_{j,\bar\nabla}^2,
		\qquad j\geq m,
	\end{align*}
	where $c_j>0$.
	
	On every finite time interval this inequality bounds each higher
	norm by its initial value and Gronwall's inequality. If a finite
	maximal time existed, the equation would also bound every spatial
	norm of $\partial_t h$ on that interval and hence $h$ would have a
	smooth limit at the maximal time. Positivity and uniform
	parabolicity persist by the low-order bound, so local existence
	from this limit would extend the solution. This proves global
	existence.
	
	The same inequalities give exponential higher-order decay in the same manner as such proof before and we get 
    \begin{align*}
		\|h(t)\|_{j,\bar\nabla}^2
		&\leq e^{-\gamma t}
		\left(\|h(0)\|_{j,\bar\nabla}^2+
		\frac{C'_j}{\lambda-\gamma}\|h(0)\|_{H^m}^2\right).
	\end{align*}
	Thus every higher norm decays, although its initial size need not
	be small.
 
 We finally remove the DeTurck term. Let $F_t$ solve
	\begin{align*}
		\partial_tF_t&=-V(\widetilde\varphi_t)\circ F_t,
		&F_0&=\operatorname{id},
		&\widetilde\varphi_t&=*_{\widetilde\psi_t}\widetilde\psi_t.
	\end{align*}
	The vector field vanishes at the reference structure and uses only
	one derivative of $h$. Since $m=r+8$, Sobolev embedding and the
	low-order decay imply
	\begin{align*}
		\|h(t)\|_{C^{r+3,\alpha}}
		+\|V(\widetilde\varphi_t)\|_{C^{r+2,\alpha}}
		&\leq C e^{-\lambda t/2}\|h(0)\|_{H^m}.
	\end{align*}
	The vector field is therefore integrable in these norms, with a
	small integral. The same idea of differentiating the equation for $F_t$ and using Gronwall's inequality and then integrating it proves that it is a diffeomorphism and is isotopic to the identity.
    
	Put $\psi_t=F_t^*\widetilde\psi_t$. The pullback derivative formula
	cancels the DeTurck term, so
	\begin{align*}
		\partial_t\psi_t
		&=F_t^*\QQ_0(\widetilde\psi_t).
	\end{align*}
	The operator $\QQ_0$ uses two derivatives and vanishes at
	$\bar\psi$. The preceding estimates and the uniform pullback
	bounds therefore give
	\begin{align*}
		\|\partial_t\psi_t\|_{C^{r,\alpha}}
		&\leq C e^{-\lambda t/2}\|h(0)\|_{H^m},\\
		\|\psi_t-\bar\psi\|_{C^{r,\alpha}}
		&\leq C\|h(0)\|_{H^m}.
	\end{align*}
	The second inequality includes both $F_t^*h(t)$ and
	$F_t^*\bar\psi-\bar\psi$. Integrating the first inequality proves
	the required time-integral bound. Taking $N=m$ and reducing
	$\delta$ supplies the initial $H^m$ smallness and proves
	\eqref{eq-ambient-coflow-bounds}. The higher estimates give the
	remaining assertions, with a smaller common exponential rate if
	necessary. Smooth uniqueness for the ordinary modified coflow is proved by Grigorian.
\end{proof}

\begin{proof}[Proof of Theorem~\ref{thm:coflow-global}]
	We use Theorem~\ref{thm:moving-stability} with $s=8$ and
	$\psi_\star=\bar\psi$. We first choose the finite order at which
	the ambient path must be controlled. Apply the estimates in the
	proof of Theorem~\ref{thm:global-instanton-transport} at an order
	which controls the transported connection in $C^9$. Let $r_0$ be
	a sufficient ambient order for these estimates, and fix
	$r\geq\max\{r_0,10\}$. All these choices are made before the
	initial neighbourhood is chosen.
	
	Set $N=r+8$. Lemma~\ref{lemma-ambient-coflow}, at the order $r$
chosen above, gives a global ordinary modified coflow $\psi_t$
with parameter zero whenever $\|\psi_0-\bar\psi\|_{C^{N,\alpha}}$
is sufficiently small. Put $Q_t=\partial_t\psi_t$. Its estimates give
	\begin{align*}
		&\sup_{t\geq0}\|\psi_t-\bar\psi\|_{C^r}
		+\sup_{t\geq0}\|Q_t\|_{C^r}
		+\int_0^\infty\|Q_t\|_{C^r}\,dt
		\leq C\|\psi_0-\bar\psi\|_{C^{N,\alpha}}.
	\end{align*}
	The higher spatial norms of $\psi_t$ are uniformly bounded, and
	$Q_t$ decays exponentially in every spatial norm. In particular,
	the higher-order time integrals are finite.
	
	For sufficiently small initial distance, these estimates imply
	\eqref{eq:small-instanton-path}. Thus
	Theorem~\ref{thm:global-instanton-transport} supplies a global
	instanton family $B_t$, with
	\begin{align*}
		B_0=\mathcal I(\psi_0),\qquad
		\partial_tB_t=\UU_{B_t,\psi_t}(Q_t).
	\end{align*}
	We spell out the quantitative bound needed for the next theorem.
	The local branch from Lemma~\ref{lem:IFT} satisfies
	$\mathcal I(\bar\psi)=A_\star$. Its smooth dependence at the
	chosen finite order bounds
	$\mathcal I(\psi_t)-A_\star$ by $\psi_t-\bar\psi$.
	The estimate \eqref{eq:transport-gauge-bound} and the differentiated
	gauge equation in the proof of
	Theorem~\ref{thm:global-instanton-transport} bound the additional
	change by the time integral of $Q_t$. Our choice of $r$ therefore
	gives
	\begin{align*}
		\sup_{t\geq0}\|B_t-A_\star\|_{C^9}
		&\leq C\left(
		\sup_{t\geq0}\|\psi_t-\bar\psi\|_{C^r}
		+\int_0^\infty\|Q_t\|_{C^r}\,dt\right)\\
		&\leq C\|\psi_0-\bar\psi\|_{C^{N,\alpha}}.
	\end{align*}
	The constants are uniform on the fixed finite-order neighbourhood.
	The higher-order versions of the transport estimates give uniform
	higher bounds for $B_t$ and exponential convergence to an instanton
	$B_\infty$ for the limiting form $\psi_\infty$.
	
	Since $r\geq10$, the preceding estimates imply
	\begin{align*}
		\sup_{t\geq0}\left(
		\|\psi_t-\bar\psi\|_{C^{10}}
		+\|B_t-A_\star\|_{C^9}
		+\|Q_t\|_{C^9}\right)
		\leq C\|\psi_0-\bar\psi\|_{C^{N,\alpha}}.
	\end{align*}
	Reduce the initial neighbourhood so that this is below the
	$\epsilon$ in \eqref{eq:background-small}. This checks that
	condition with $s=8$. All the uniform higher spatial bounds in
	Theorem~\ref{thm:moving-stability} have also been verified.
	
	Finally,
	\begin{align*}
		\|A_0-B_0\|_{H^8}
		&\leq \|A_0-A_\star\|_{H^8}
		+\|\mathcal I(\psi_0)-A_\star\|_{H^8}\\
		&\leq \|A_0-A_\star\|_{H^8}
		+C\|\psi_0-\bar\psi\|_{C^{N,\alpha}}.
	\end{align*}
	We choose the initial neighbourhood small enough that this is
	below the initial radius in Theorem~\ref{thm:moving-stability}.
	That theorem now gives the global smooth connection evolution and
	smooth exponential convergence to $k_\infty^*B_\infty$.
	Here $k_\infty$ is the limiting gauge transformation in that
	theorem. Lemma~\ref{lemma-ambient-coflow} gives smooth exponential
	convergence of $\psi_t$ to the torsion-free form $\psi_\infty$.
	
	The gauge-fixed connection remains in the neighbourhood where the
	deformation operator has a uniformly bounded inverse. Gauge
	conjugation preserves invertibility, so the correction remains
	defined for the original connection throughout the solution.
	The limiting instanton is acyclic for the same reason. If $A_0$
	is initially an instanton, Theorem~\ref{thm:energy} preserves its
	zero defect. Uniqueness follows from ambient uniqueness and
	Theorem~\ref{thm:nearby-existence}.
\end{proof}

\medskip

We relate our results to another family of flows in $G_2$ geometry which concern the \emph{heterotic} $G_2$ system. Geometric flows related to the heterotic $G_2$ system have been studied by Ashmore, Minasian and Proto~\cite{AMP24}.
Garcia-Fernandez, Moreno, Payne and Streets~\cite{GMPS25} introduce
a parabolic flow for the large volume heterotic $G_2$ system and
prove short-time existence and Shi-type estimates. Related
stationary coupled equations for the geometry and a bundle
connection are studied by da Silva Jr., Garcia-Fernandez, Lotay
and S\'a Earp~\cite{SGLS26}.

Karigiannis, Picard and Suan ask in
\cite[Remark~4.13]{KPS} for a natural connection flow towards a
$G_2$ instanton in the anomaly-flow setting. Our
Theorem~\ref{thm:moving-stability} gives a partial result in this
direction. Near an acyclic instanton, it gives exponential
attraction for the corrected connection flow along prescribed
coclosed backgrounds satisfying its hypotheses. The correction is
nonlocal, since it is obtained by an elliptic inverse. The ambient
evolution used here does not depend on the evolving connection.
The full $\alpha'$-coupled anomaly flow, where the curvature of the
connection also enters the ambient equation, is not covered by
this argument.

\subsection{The coflow parameter}\label{sec:coflow-parameter}
The parameter $c$ multiplies a first-order term. It does not
change the principal symbol, but it can change the sign of a
low-frequency eigenvalue. The following example shows that, for
each nonzero parameter, stability can fail at the linear level
on a suitable flat background.

\begin{proposition}\label{prop:positive-c-instability}
For every $c\in\mathbb R\setminus\{0\}$, there is a flat
torsion-free $G_2$ structure on a flat seven-torus for which the
linearization of $\partial_t\psi=Q_c(\psi)$ has a positive
eigenvalue on exact closed four forms. Its eigenmode is
$L^2$ orthogonal to all infinitesimal diffeomorphism directions.
\end{proposition}

\begin{proof}
Write the torus as a product of a circle of length $L$ and a flat
six-torus, with an adapted $SU(3)$ structure on the six-plane.
Let $s$ be arclength on the circle. On real primitive three forms
of type $(2,1)+(1,2)$ on the six-plane, set
\begin{align*}
 J\eta=*(ds\wedge\eta).
\end{align*}
On this space $J^2=-I$. Choose a nonzero constant form $\eta_0$
in this space, take $k=2\pi/L$, and let $\epsilon\in\{1,-1\}$.
Define
\begin{align*}
 \eta_\epsilon(s)
 &=\cos(ks)\eta_0-\epsilon\sin(ks)J\eta_0,
 &\zeta_\epsilon&=*\eta_\epsilon.
\end{align*}
Differentiation and $J^2=-I$ give
\begin{align*}
 \partial_s\eta_\epsilon&=-\epsilon kJ\eta_\epsilon,\\
 *d\eta_\epsilon&=J\partial_s\eta_\epsilon
                    =\epsilon k\eta_\epsilon,\\
 d\eta_\epsilon&=\epsilon k\zeta_\epsilon,
 &\Delta\zeta_\epsilon&=k^2\zeta_\epsilon.
\end{align*}
Thus $\zeta_\epsilon=(\epsilon/k)d\eta_\epsilon$ is exact and
closed. The $SU(3)$ decomposition gives
$\eta_\epsilon\in\Omega^3_{27}$ and
$\zeta_\epsilon\in\Omega^4_{27}$.

We compute the derivative of each term in $Q_c$ at the flat
torsion-free form $\bar\psi$. By \eqref{eq:dual-variation}, the
variation of the associated three form in the direction
$\zeta_\epsilon$ is $\dot\varphi=-\eta_\epsilon$.
For closed four forms, $\Delta_\psi\psi=d*d\varphi$.
Since $d\bar\varphi=0$, differentiating the Hodge star in this
expression contributes zero. Hence the Laplacian term gives
\begin{align*}
 d*d\dot\varphi=-d*d\eta_\epsilon=-k^2\zeta_\epsilon,
\end{align*}
and
$\operatorname{tr}T=\tfrac14\langle d\varphi,\psi\rangle$.
At the torsion-free reference form,
\begin{align*}
 D(\operatorname{tr}T)[\zeta_\epsilon]
 &=\tfrac14\langle d\dot\varphi,\bar\psi\rangle\\
 &=-\tfrac{\epsilon k}{4}
       \langle\zeta_\epsilon,\bar\psi\rangle=0.
\end{align*}
The last equality uses $\Omega^4_{27}\perp\Omega^4_1$.
The derivative of the torsion-trace term therefore vanishes.
Finally, the derivative of $-2c\,d\varphi$ is
$2c\,d\eta_\epsilon$. Combining these terms gives
\begin{align*}
 DQ_c|_{\bar\psi}\,\zeta_\epsilon
 =(-k^2+2\epsilon ck)\zeta_\epsilon.
\end{align*}
Choose $\epsilon=1$ when $c>0$ and $\epsilon=-1$ when $c<0$.
The eigenvalue is then $-k^2+2|c|k$, which is positive if
$0<k<2|c|$. It is enough to take $L>\pi/|c|$.

We finally check that this mode is not an infinitesimal
diffeomorphism. Since $d\bar\psi=0$, every such direction is
$\mathcal L_Y\bar\psi=d(Y\mathbin{\lrcorner}\bar\psi)$.
On four forms in dimension seven, $\delta=*d*$, so
\begin{align*}
 \delta\zeta_\epsilon
 =*d\eta_\epsilon=\epsilon k\eta_\epsilon\in\Omega^3_{27}.
\end{align*}
On the other hand,
$Y\mathbin{\lrcorner}\bar\psi\in\Omega^3_7$. Integration by
parts and orthogonality of these summands give
\begin{align*}
 \langle\zeta_\epsilon,\mathcal L_Y\bar\psi\rangle_{L^2}
 &=\langle\delta\zeta_\epsilon,
                   Y\mathbin{\lrcorner}\bar\psi\rangle_{L^2}
 =0.
\end{align*}
Since $\zeta_\epsilon\ne0$, it cannot itself be an
infinitesimal diffeomorphism direction.
\end{proof}

\begin{remark}
For each nonzero parameter the proposition produces a suitable
background with a growing linear mode. The background, in
particular the circle length, is allowed to depend on the
parameter. This does not prove instability for every nonzero
parameter on a fixed background, and it is not a nonlinear
instability theorem. On a fixed torus the displayed modes have
eigenvalues $-k^2\pm2ck$; their signs depend on the allowed
frequencies as well as on $c$. The example therefore explains why
an extension of the global stability theorem to nonzero parameters
requires information about the spectrum of the chosen background.
\end{remark}

\subsection{Examples with nonzero correction}\label{sec:connection-examples}
In this section, we give examples in which the correction $\mathcal U_{A, \psi}(Q)$ described in previous sections is nonzero, so the connection must move to remain an instanton as the ambient $\psi_t$ evolves. We first compute the correction for an ambient family obtained by diffeomorphisms, and then construct examples with nonzero initial correction along the modified Laplacian coflow.

\begin{proposition}
\label{prop:diffeo-example}
Let $A$ be an acyclic instanton for a closed positive $\psi$ and
let $Y$ be an ambient vector field. For $Q=\Lie_Y\psi$, put
$v=Y\ip F_A$. Then
\begin{align}
\label{eq:diffeo-u}
 \UU_{A,\psi}(Q)
 =v-d_A(d_A^*d_A)^{-1}d_A^*v.
\end{align}
%Whenever $*(F_A\wedge\Lie_Y\psi)\neq0$, this correction is nonzero. 
If $F_A$ is not identically zero, there is a compactly
supported $Y$ in a coordinate ball for which it is nonzero.
\end{proposition}

\begin{proof}
We move both the ambient $4$-form and the connection by the flow of the vector field $Y$.
The instanton equation is preserved by this simultaneous pullback.
In a local trivialization, the covariant Cartan formula writes the
connection velocity as $Y\ip F_A$ plus an infinitesimal gauge term.
At an instanton $\PP_Ad_A=0$, so that gauge term does not affect the
linearized instanton equation. Differentiating the equation therefore
gives
\begin{align*}
 \PP_A(Y\ip F_A)=-*(F_A\wedge\Lie_Y\psi).
\end{align*}
The velocity $v=Y\ip F_A$ need not be Coulomb. We solve
$d_A^*d_Aq=d_A^*v$ and replace $v$ by $v-d_Aq$.
Acyclicity makes the scalar Laplacian invertible, so $q$ is determined
uniquely. The new velocity satisfies $d_A^*(v-d_Aq)=0$, and its
image under $\PP_A$ is unchanged because $\PP_Ad_Aq=0$.
The uniqueness in Proposition~\ref{prop:horizontal} now gives
\eqref{eq:diffeo-u}. If the forcing in the differentiated equation
is nonzero, its solution cannot be zero.

To arrange nonzero forcing, choose a point $p$ with $F_A(p)\neq0$
and a nonzero ordinary two-form component $\beta$ of the curvature
in a frame of $E$. Choose a nonzero coefficient of $\beta$ in a
coframe and wedge with four of the other five basis covectors.
The resulting four form $q$ has $\beta\wedge q\neq0$.
Because positive four forms form an open $\mathrm{GL}(7)$ orbit (or are stable forms), the derivative of the orbit map is onto the space of four forms.
We can therefore prescribe $Y(p)=0$ and the first derivative of $Y$
at $p$ so that $(\Lie_Y\psi)(p)=q$. Multiplying it 
by a cutoff equal to one near $p$ preserves these data and gives
a smooth vector field supported in a coordinate ball. The forcing
has a nonzero component at $p$ by construction.

The ambient evolution of the $4$-form obtained by pulling back $\psi$ remains closed and positive. Thus the example fits the prescribed-family coupling as
well as the instantaneous compatibility equation.
\end{proof}

\begin{proposition}
\label{prop:actual-coflow-example}
Suppose $(M,\bar\psi)$ is compact and torsion-free and carries a
nonflat acyclic instanton $A_\star$. Fix a constant $c\in \mathbb R$ and set
\begin{align*}
 \QQ_c(\psi)=\Delta_\psi\psi-2d\big((c+\tr T_\psi)\varphi\big).
\end{align*}
Then there are closed positive forms $\psi_\varepsilon$, cohomologous
and arbitrarily close to $\bar\psi$, and acyclic
$\psi_\varepsilon$-instantons $A_\varepsilon$ for which
\begin{align*}
 *_{\psi_\varepsilon}
   (F_{A_\varepsilon}\wedge \QQ_c(\psi_\varepsilon))\neq0.
\end{align*}
Consequently their correction is initially nonzero along the
actual modified-coflow coupling. When $c=0$ and the perturbation
is sufficiently small, these examples belong to the global
convergent solutions of Theorem~\ref{thm:coflow-global}.
\end{proposition}

\begin{proof}
Choose $p$ where a scalar component $\beta$ of $F_{A_\star}$ is nonzero.
It is an element of $\Lambda^2_{14}=\mathfrak g_2$. the Lie algebra of the Lie group $G_2$. We conjugate
inside $G_2$ and put it in a maximal torus, so that in an adapted
coframe on the six-plane perpendicular to $ds$,
\begin{align*}
 \omega=e^{12}+e^{34}+e^{56},\qquad
 \beta=a e^{12}+b e^{34}+d e^{56},\qquad a+b+d=0.
\end{align*}
Not all $a,b,d$ vanish. If, for example,
$a\neq b$, take $\theta=(e^{12}-e^{34})\wedge e^5$.
Then $\theta$ is a real primitive form of type
$(2,1)+(1,2)$ and
$\beta\wedge\theta=(b-a)e^{12345}\neq0$.
%For a different unequal pair, permute the three complex planes.
Put $\widehat\zeta=ds\wedge\theta$. The $\mathrm{SU}(3)$ decomposition gives $\widehat\zeta\in\Lambda^4_{27}$. Thus $F_{A_\star}(p)\wedge\widehat\zeta\neq0$.

We next choose an exact perturbation whose coflow velocity at $p$
contains this four-form direction. Since coflow is second order,
we can do this by prescribing the second derivatives of the
perturbation while making its value and first derivatives vanish. In fact, the same computation as in Proposition~\ref{prop:positive-c-instability} shows that the symbol satisfies
\begin{align}
\label{eq:coflow-symbol}
 \sigma_{ds}(\mathrm D\QQ_c|_{\bar\psi})\widehat\zeta=-\widehat\zeta.
\end{align}

Now we use normal coordinates at $p$ and extend $\theta$ with constant
coordinate coefficients. Let $\rho$ be a cutoff equal to one near
$p$, and define
\begin{align*}
 \zeta=d\left(\frac{\rho s^3\theta}{6}\right).
\end{align*}
Near $p$, the cutoff is one and $\theta$ has constant coordinate
coefficients, so $\zeta=\frac12 s^2\,ds\wedge\theta$. At $p$, where
$s=0$, both $\zeta$ and all its first derivatives vanish. Its only
nonzero second coordinate derivative is
$\partial_s^2\zeta(p)=\widehat\zeta$. Interpreting
\eqref{eq:coflow-symbol} as the coefficient of this coordinate
derivative gives
\begin{align*}
 (\mathrm D\QQ_c|_{\bar\psi})\zeta(p)=\widehat\zeta.
\end{align*}
The form $\zeta$ need not lie in $\Omega^4_{27}$ away from $p$.
Only the algebraic projection at $p$ enters the second-order
coefficient applied to its nonzero second derivative
$\partial_s^2\zeta(p)=\widehat\zeta$. Every term in which a
derivative hits a projection or another coefficient contains
$\zeta(p)$ or a first derivative of $\zeta$ at $p$. These vanish
by construction. All lower-order terms vanish for the same
reason, including the first-order term containing $c$.

Set $\psi_\varepsilon=\bar\psi+\varepsilon\zeta$.
It is closed, cohomologous to $\bar\psi$, and positive for
sufficiently small $\varepsilon$. By
Lemma~\ref{lem:IFT} there is a smooth acyclic
instanton family $A_\varepsilon$ with $A_0=A_\star$. Since
$\QQ_c(\bar\psi)=0$,
\begin{align*}
 *_{\psi_\varepsilon}
  (F_{A_\varepsilon}\wedge \QQ_c(\psi_\varepsilon))(p)
 =\varepsilon*_{\bar\psi}
  (F_{A_\star}\wedge\widehat\zeta)(p)+O(\varepsilon^2).
\end{align*}
When differentiating this expression at $\varepsilon=0$, the
variations of the Hodge star and of the curvature multiply
$\QQ_c(\bar\psi)=0$ and hence disappear. The remaining first-order
term is the one displayed above, and it is nonzero by our choice of
$\widehat\zeta$. Thus the forcing is nonzero for all sufficiently
small nonzero $\varepsilon$. At an instanton the defining equations
give $\PP\UU=-*(F_A\wedge Q_c)$, so the correction cannot vanish.

These are initial data for actual solutions. The modified coflow has
smooth short-time existence for them \cite{Gri,BVStab20}, and
Theorem~\ref{thm:instanton} supplies the connection lift. For $c=0$,
we may make the perturbation small enough to apply
Theorem~\ref{thm:coflow-global}, which proves the global assertion.
\end{proof}

\begin{remark}
The hypothesis in Proposition~\ref{prop:actual-coflow-example}
has compact examples. The gluing theorem of \cite{Wal}
produces acyclic instantons on compact torsion-free $G_2$
manifolds, and choosing nontrivial instanton gluing data gives
nonflat curvature. The proposition uses those examples as initial
backgrounds; it does not require a new compact instanton gluing
construction.
\end{remark}

\subsection{The obstruction without acyclicity}\label{sec:nonacyclic}

In the preceding sections, acyclicity made the augmented operator
invertible and supplied the spectral gap used in the stability
proof. When its kernel is nonzero, a compatible connection velocity
may fail to exist. The relation between ambient deformations and
instanton deformations was studied by de la Ossa, Larfors and
Svanes~\cite[Sections~5.3--5.4]{OLS16}, using a curvature map which
determines which infinitesimal ambient deformations preserve the
instanton condition. Here we state the solvability condition for
a prescribed closed variation $Q$ by applying the Fredholm
alternative to the self-adjoint operator $\DD_{A,\psi}$. Thus the
source must be orthogonal to its kernel. This is a linear
compatibility statement and it does not give nonlinear continuation
through a change in the kernel dimension.

\begin{proposition}
\label{prop:nonacyclic}
Let $A$ be any instanton for a closed positive $\psi$ on compact
$M$, and let $Q$ be a closed four form. Put
$b=*(F_A\wedge Q)$ and
\begin{align*}
 \mathcal H^1_{A,\psi}=\{a:d_A^*a=0,\ \PP_Aa=0\}.
\end{align*}
There is a compatible velocity $v$ satisfying
\begin{align}
\label{eq:compatible-general}
 d_A^*v=0,\qquad \PP_Av=-b
\end{align}
if and only if
\begin{align}
\label{eq:obstruction}
 \langle b,a\rangle_{L^2}=0\quad
 \text{for every }a\in \mathcal H^1_{A,\psi}.
\end{align}
When this holds, there is a unique such velocity orthogonal to
$\mathcal H^1_{A,\psi}$. It has minimum $L^2$ norm among all velocities
satisfying the differentiated instanton equation, including those
not in Coulomb gauge. On smooth families with constant kernel
dimension, these minimum norm solutions depend smoothly provided
\eqref{eq:obstruction} holds.
\end{proposition}

\begin{proof}
At an instanton the kernel of the self-adjoint elliptic operator
$\DD$ is $\ker d_A|_{\Omega^0}\oplus\mathcal H^1_{A,\psi}$.
The range is its $L^2$ orthogonal complement. The source $(0,-b)$
is already orthogonal to the scalar part of the kernel, and its
orthogonality to the one-form part is exactly
\eqref{eq:obstruction}. This proves that the augmented equation
$\DD(f,v)=(0,-b)$ is solvable if and only if that condition holds.

When it holds, choose the unique solution orthogonal to $\ker\DD$.
Apply $d_A^*$ to its second equation and as before we get $f=0$. The remaining component $v$ satisfies\eqref{eq:compatible-general} and is orthogonal to
$\mathcal H^1_{A,\psi}$.

To compare with any other compatible velocity $w$, including one
that is not Coulomb, we follow the same strategy as in the proof of Proposition~\ref{prop:horizontal} and obtain the minimum norm property and its uniqueness

For a smooth family with constant kernel dimension, identify the
varying $L^2$ spaces by the positive square roots of their metric
and density changes. On this fixed Hilbert space choose, at any reference parameter,
a circle $|z|=\rho$ which encloses zero but no other eigenvalue.
The resolvent remains defined on that circle for nearby parameters,
and the projection is
\begin{align*}
 \Pi=\frac{1}{2\pi\sqrt{-1}}
       \int_{|z|=\rho}(z-\DD)^{-1}\,dz.
\end{align*}
Smooth dependence of the resolvent gives smooth dependence of this
projection. Its rank is locally constant. Since the kernel dimension
is also constant by hypothesis, the enclosed spectral subspace cannot
contain any additional nonzero eigenvectors nearby and hence it is exactly
$\ker\DD$. On that kernel, $\DD+\Pi$ acts as the identity. On the
orthogonal complement it acts as the invertible restriction of
$\DD$. Therefore
\begin{align*}
 \DD^\dagger=(\DD+\Pi)^{-1}(I-\Pi)
\end{align*}
vanishes on the kernel and inverts $\DD$ on its orthogonal
complement. This formula shows its smooth dependence on the
parameters, and elliptic regularity gives the corresponding maps
between higher Sobolev spaces. Applying it to $(0,-b)$ gives the
minimum norm solution smoothly. The construction is local in the
parameter. The locally constructed solutions agree wherever their
domains overlap, by uniqueness, so they define smooth solutions
along the given family. 
\end{proof}
\begin{remark}
In the proof above, uniform bounds over the whole family would
require uniform estimates in addition to the local argument provided.
\end{remark}
\begin{remark}\label{rem:remacyc}
This proposition gives the exact first-order obstruction to
removing acyclicity. Constant kernel dimension on a chosen
family is an additional hypothesis, and is generally not an
open condition on all connections. Thus the proposition does
not define a smooth correction on an arbitrary open
neighborhood of a nonacyclic instanton. If an eigenvalue tends
to zero, the inverse on its eigenspace can diverge unless the
corresponding forcing component vanishes at a sufficient
rate. A general nonlinear transport theorem across such a
crossing would need further information about the moduli
space and these obstruction components.
\end{remark}

\section{Further questions}\label{sec:questions}
We saw many unifying themes in constructing coupled system of flows of the ambient $G_2$ structures (via $\varphi$ and $\psi$) and two naturally related objects in the form of coassociative submanifolds moving by the mean curvature flow and instantons moving by the Yang-Mills flow. There were many unifying themes between the ideas behind these objects. Both constructions choose a velocity by an $L^2$ minimization
condition, but the quantities being minimized are different. For
coassociatives, we minimize the correction to mean curvature. With
ordinary Laplacian flow in the ambient space, this gives the normal
velocity $\Pi_iH-R_i d_C(i^*\tau)$. It is the compatible velocity
closest to $H$, and its kernel component is $\Pi_iH$. It need not be
the compatible velocity of smallest norm. For instantons, we remove
the gauge component and choose the compatible velocity of smallest
$L^2$ norm. In both cases the correction is nonlocal: an elliptic
inverse is applied over the whole submanifold or ambient manifold.
The preceding theorems show that these choices give actual evolutions
and convergence under the stated neighborhood and stability assumptions. All these points give rise to several interesting questions for further research. We mention some of them below.

\medskip
\noindent\textbf{1. Can the embedding flow start from an
	embedding that is not coassociative?}
Our correction is defined for embeddings close to a coassociative
embedding, even when their pullback is not zero. Does the corrected
embedding equation have a smooth solution for a short time from
every such initial embedding?

Our existence proof does not answer this question, since it uses
the finite-dimensional family of embeddings whose pullback is zero.
An existence proof for initial embeddings with nonzero pullback
would therefore require a different argument. A separate question
is whether the correction can be changed so that a small nonzero
pullback tends to zero.

\medskip
\noindent\textbf{2. Can the coassociative embeddings develop
	a singularity while the ambient flow remains smooth?}
More precisely, can the tubular radius tend to zero or the second
fundamental form become unbounded in finite time, even though the
ambient structure remains smooth? If this can happen, what does
the degeneration look like? If it cannot happen under suitable
additional assumptions, what estimates prevent it?

Proposition~\ref{prop:coassoc-continuation} gives continuation under
its compactness and tubular bounds. These assumptions allow us to
keep writing nearby embeddings as normal graphs. The proposition
does not explain what happens when these bounds fail. Answering
the questions above would help us understand whether the embedding
can break down before the ambient flow does.

\medskip
\noindent\textbf{3. Can the connection flow continue when
	the kernel dimension changes?}
Suppose that the deformation operator develops a larger kernel
along a family of instantons. Under what additional conditions can
the connection flow be continued smoothly through this point? See also Remark~\ref{rem:remacyc}.

\medskip\noindent\textbf{4. For which positive coflow parameters
	does stability hold?}
Fix a torsion-free background. Is there a range of nonzero
parameters for which sufficiently small perturbations remain
close to the background and converge? Can such a result give an
analogue of Theorem~\ref{thm:coflow-global}?

Our global coupled theorem uses parameter zero.
Proposition~\ref{prop:positive-c-instability} produces, for every
nonzero parameter, a flat background with a growing linear mode.
The background is allowed to vary with the parameter. The
proposition therefore does not answer the stability question for
a fixed background. One would first need to identify parameter
values for which the linearized operator has the required decay
on the fixed cohomology class, modulo diffeomorphisms, and then
control the nonlinear terms.

\medskip
\noindent\textbf{5. Can attraction to an instanton be proved
	for larger initial perturbations?}
What additional geometric bounds would allow us to prove global
existence and convergence when the initial connection is not
sufficiently close to an acyclic instanton?

Corollary~\ref{cor:defect-decay} gives decay of the defect under its
spectral and metric assumptions but we need further estimates to prevent curvature concentration, where the curvature becomes large in increasingly
small regions. We would also need to control the ambient structure
and address any loss of invertibility of the deformation operator.
The question is which additional assumptions provide these bounds
and allow the attraction theorem to extend beyond small initial
perturbations.

For ordinary Yang--Mills flow on a fixed compact torsion-free
$G_2$ manifold, Oliveira and Waldron~\cite[Theorem~1.1]{OW} show
that a uniform pointwise bound on the $\Lambda^2_7$ component of
the curvature prevents a finite-time singularity. The $L^2$ defect
estimate in Corollary~\ref{cor:defect-decay} does not provide this
pointwise bound. Can a corresponding continuation criterion be
proved for our corrected connection flow along a moving coclosed
background?

\medskip
\noindent\textbf{6. Are there similar coupled flows in Spin(7)-geometry?} There is a closely related exceptional geometry on $8$-dimensional manifolds where the special geometric structure is described by an admissible $4$-form $\Phi$. The submanifolds corresponding to coassociatives are called Cayley submanifolds and there is an analogous notion of Spin(7)-instantons. So can we define similar coupled flows as in this paper for those objects? Natural flows of Spin(7)-structures are discussed in \cite{Dwi-RHF} and \cite{Dwi-spin}. Probably the first step would be use the results in \cite{BM} and do a computation similar to Lemma~\ref{lem:H} and Proposition~\ref{prop:unmodified-pullback} by using particulars of Spin(7)-geometry.

\printbibliography

\noindent
(SD:) Fachbereich Mathematik, Universität Hamburg, Bundesstraße 55, 20146 Hamburg, Germany.\\
\href{shubham.dwivedi@uni-hamburg.de}{shubham.dwivedi@uni-hamburg.de}

\medskip
\noindent
(SP:) \href{siva17aug@gmail.com}{siva17aug@gmail.com}

\end{document}